%% file: burgers.tex
\documentclass[12pt]{article}
\usepackage[margin=1.1in]{geometry}
\usepackage{authblk}
\usepackage{amsmath}
\usepackage{amssymb}
\usepackage{amsthm}
\usepackage{mathrsfs}
\usepackage{bbm}
\usepackage{empheq}
\usepackage[tight]{subfigure}
\usepackage{epsfig}
\usepackage{graphicx}
\usepackage{psfrag}
\usepackage[usenames,dvipsnames]{pstricks}
\usepackage{pst-plot}
\usepackage[colorlinks]{hyperref}
\usepackage{tabls}
\usepackage{paralist}
\usepackage{hyperref}
\usepackage{caption}

\DeclareSymbolFontAlphabet{\Bbb}{AMSb}

\newlength{\fixboxwidth}
\newcommand{\COMMENT}[1]{}

\newcommand{\bs}{\boldsymbol}

\newcommand{\defeq}{\mathrel{\mathord{:}\mathord{=}}}

\newcommand{\N}{\mathbb{N}}

\newcommand{\E}{\mathbb{E}}

\renewcommand{\P}{\mathbb{P}}

\newcommand{\R}{\mathbb{R}}

\newcommand{\quark}{\setbox0\hbox{$x$}\hbox to\wd0{\hss$\cdot$\hss}}

\newtheorem{thm}{Theorem}[section]

\theoremstyle{definition}

\hypersetup{
    linkcolor=blue,
}

\renewcommand{\thefigure}{\arabic{figure}}

\makeatletter
\renewcommand{\p@subfigure}{\thefigure}
\makeatother

\newcounter{mycount}

\makeatletter
\def\blfootnote{\gdef\@thefnmark{}\@footnotetext}
\makeatother

\begin{document}

\title{Optimal Bounds on Nonlinear Partial Differential Equations in Model Certification, Validation, and Experiment Design}

\author[1]{M.\ McKerns}
\author[2]{F.\ J.\ Alexander}
\author[3]{K.\ S.\  Hickmann}
\author[4]{T.\ J.\ Sullivan}
\author[5]{D.\ E.\ Vaughan}

\affil[1]{Information Sciences\\
Los Alamos National Laboratory\\
Los Alamos, New Mexico 87545, USA\\
E-mail: mmckerns@lanl.gov}
\affil[2]{Computational Science Initiative\\
Brookhaven National Laboratory\\
Upton, New York 11973, USA}
\affil[4]{Institute of Mathematics\\
Free University of Berlin\\
14195 Berlin, Germany}
\affil[3,5]{Verification and Analysis\\
Los Alamos National Laboratory\\
Los Alamos, New Mexico 87545, USA}

\date{March 24, 2019}%\today}

\maketitle
%\\
%%%%%%%%%%%%%%%%%%%%%%%%%%%%%
%\newpage
%%%%%%%%%%%%%%%%%%%%%%%%%%%%%

%%% Mike's comments:
% #1: Run OUQ P(z_{\ast} > \bar{z_{\ast}}) code on LANL resources to get a set of timings
% #2: Run OUQ P(z_{\ast} > \bar{z_{\ast}}) code with parallel ensemble optimizer
% #3: Get MC bounds for different distributions of \delta (non-uniform)
% #4: Explain how to use OUQ to calculate the optimal estimator
% #5: Explain how to incorporate data as prior information
%%% Kyle's comments:
% #6: Explain extending this to PDE with QOI derived from numerical simulation
% #7: Address model-form uncertainty [see (4)]
%%% Diane's comments:
% #8: Present MC results and MC bounds for different distributions [see (3)]
%%%

\section{Introduction}
A common goal in areas of science and engineering that rely on making
accurate assessments of performance and risk (e.g.\ aerospace
engineering, finance, geophysics, operations research) for complex
systems is to be able to guarantee the quality of the assessments
being made.  Very often, the knowledge of the system is incomplete or
contains some form of uncertainty.  There can be uncertainty in the
form of the governing equations, in information about the parameters,
in the collected data and measurements, and in the values of the input
variables and their bounds.  One of the most common cases is that
initial conditions and/or boundary conditions are known only to a
certain level of accuracy.  Even in the case where the dynamics of a
system is known exactly at a fine-grained level, computationally more
tractable coarse-grained models of the system often have to be derived
under approximation, and thus contain uncertainty.  One way to refine
a model of a system under uncertainty is to perform experiments to
help solidify what is known about the parameters and/or the form of
the governing equations.  Sampling methods (such as Monte Carlo) can
be used to determine the predictive capacity of the models under the
resulting uncertainty.  Unfortunately, this determination can be
computationally costly, inaccurate, and in many cases impractical.

One promising approach to dealing with this challenge is Optimal
Uncertainty Quantification (OUQ) developed by Owhadi et al
\cite{OSSMO:2011, Sullivan:2013, MOSSO:2010}.  OUQ integrates the
knowledge available for both mathematical models and any knowledge
that constrains outcomes of the system, and then casts the problem as
a constrained global optimization problem in a space of probability
measures; this optimization is made tractable by reducing the problem
to a finite-dimensional effective search space of discrete probability
distributions, parameterized by positions and weights.  Much of the
early work with OUQ has been to provide rigorous certification for the
behavior of engineered systems such as structures under applied
stresses or metal targets under impact by projectiles
\cite{OSSMO:2011, Kamgaouq:2014}.

\input ./oed.tex % \section{OUQ for Model Validation and Experiment Design}

% Outline of Paper
\section{Outline}
In the remainder of this manuscript, we demonstrate utility of the OUQ
approach to
understanding the behavior of a system that is governed by a partial
differential equation (and more specifically, by Burgers' equation).
In particular, we solve the problem of predicting shock location when
we only know bounds on viscosity and on the initial conditions.
We will calculate the bounds on the probability that the shock
location occurs at a distance greater than some selected target
distance, given there is uncertainty in the location of the left
boundary wall.  Through this example, we demonstrate the potential to
apply OUQ to complex physical systems, such as systems governed by
coupled partial differential equations.  We compare our results to
those obtained using a standard Monte Carlo approach, and show that
OUQ provides more accurate bounds at a lower computational cost.  As
OUQ can take advantage of solution-constraining information that Monte
Carlo cannot, and requires fewer assumptions on the form of the
inputs, the predicted bounds from OUQ are more rigorous than those
obtained with Monte Carlo.
% OUQ uses numerical optimization to discover rigorous lower and upper bounds on some quantity of interest, and (as evidenced in Figure \ref{fig:punchline}) is much more capable than standard sampling approaches in describing behavior governed by rare events.
We conclude with a brief discussion of how to extend this approach to more complex systems,
and note how to integrate our approach into a more ambitions program
of optimal experimental design.

\section{Burgers' Equation with a Perturbed Boundary}
\noindent
For a viscous Burgers' equation:
\[    
  u_t + u u_x - v u_{xx} = 0
\]
with $x \in [-1,1]$ and viscosity $v > 0$, we have:
\begin{gather*}
  u(-1) = 1 + \delta \\
  u(1) = -1
\end{gather*}
where $\delta > 0$ is a perturbation to the left boundary condition.
The solution has a transition layer, which is a region of rapid variation, and has location $z$, defined as the zero of the solution profile $u(z) = 0$ which varies with time and at steady state is extremely sensitive to the boundary perturbation
\cite{Xiu:2004}.

The exact solution of Burgers' equation with the small \emph{deterministic} perturbed boundary condition, at a steady state, is:
\[    
  %\label{eq:burgers_solution}
  u(x) = -A \: \tanh (\frac{A}{2v} (x - z_{\ast})).
\]
Here, $z_{\ast}$ is defined as the location of the transition layer where $u(z_{\ast}) = 0$, and the slope at $z_{\ast}$ is $-A = \partial u/\partial x|_{x=z_{\ast}}$.
We can solve for the two unknowns $z_{\ast}$ and $A$ with:
\begin{equation}
  \begin{gathered}
  \label{eq:burgers_steady}
  A \: \tanh (\frac{A}{2v} (1-z_{\ast})) - 1 = 0 \\
  A \: \tanh (\frac{A}{2v} (1+z_{\ast})) - 1 - \delta = 0
  \end{gathered}
\end{equation}

\section{Solving for Location of the Transition Layer}

The equations in \eqref{eq:burgers_steady} are considered difficult to solve exactly due to the high nonlinearity of the solution space.
%  and are typically solved by approximately
%numerical methods (such as \emph{asymptotic analysis})
% and $generalized$ $polynomial$ $chaos$)
%\cite{Xiu:2004}.
Solving the viscous Burgers' equation with approximate numerical
methods can still be difficult, as the aforementioned supersensitivity
of the solution can lead to a large variance in solved $z_{\ast}$ for
a distribution of input $\delta$ and $v$. We will use the
\texttt{mystic} optimization framework %\cite{McKernsHungAivazis:2009, MSSFA:2011}
to solve for $z_{\ast}$ and $A$ directly, and will
compare to approximate numerical results obtained by Xiu
\cite{Xiu:2004}.

%We use the Differential Evolution solver \texttt{mystic.solvers.diffev2}
We use an ensemble (\texttt{mystic.solvers.lattice}) of Nelder-Mead optimizers
 to solve \texttt{objective} for $z_{\ast}$ and $A$, for the given \texttt{bounds}, $v$, and $\delta$. The \texttt{bounds} and \texttt{objective} are defined as in Figure \ref{code:exactly}, bounding $z_{\ast}$ in $[0,1]$, and with an \texttt{objective} to find the minimum of $\sum_{i=0}^{1} {abs(eqns_{i})}$. Values for solved $z_{\ast}$ are presented in Table \ref{tab:exactly}, and reproduce exactly the results of numerical solution of the exact formula found by Xiu. %XXX: Xiu solves exactly, however uses a local optimizer which requires an inital trial solution close to the global minimum.
Xiu notes that iterative methods were used to solve the system of nonlinear equations, and that the convergence is very sensitive to the choice of initial guess.
%With \texttt{mystic.solvers.diffev2},
With \texttt{mystic.solvers.lattice},
the solution is found reliably for $z_{\ast}$ in $[0,1]$ and $A$ unbounded, where convergence with \texttt{mystic} is between
% $2.5e^{-12}$ and $2.8e^{-15}$
$1.1e^{-16}$ and $8.6e^{-18}$
for all given combinations of $v$ and $\delta$. On a 2.7 GHz Macbook with an Intel Core i7 processor and 16 GB 1600 MHz DDR3 memory, running python
% 2.7.12 and \texttt{mystic} 0.3.0,
2.7.14 and \texttt{mystic} 0.3.2,
solutions to \eqref{eq:burgers_steady} are found in $\sim$
%0.5 seconds.
0.15 seconds.

% \begin{compactitem}
 %\item Discuss code for solving Burgers' Equation.
 %\item Compare results from \texttt{mystic} to results in Xiu.
 %\item Discuss results validity and computational expense.
% \end{compactitem}

%TABLE: Location of transition layer found by solving Burgers' equation exactly.
%       Contents: \delta versus z, for two values of v  (Table I and II in Xiu)
%       Also provide A and fit?  (Output of exact_supersensitive.py)
\begin{table}[htbp!]
  \centering
  \begin{tabular}{c c c}
    \hline
   %\multicolumn{1}{c}{$ $} & \multicolumn{2}{c}{$z$} \\
   %\cline{2-3}
    $\delta$ & $z_{\ast}$ with $v=0.1$ & $z_{\ast}$ with $v=0.05$ \\ % A A 
    \hline
   %1e-1 & 0.723225247335678000 & 0.8616126164655098 \\
    1e-1 & 0.72322525 & 0.86161262 \\
% -1.1000000128879028 -1.0999999999977252
   %1e-2 & 0.474927411637441950 & 0.7374601530661407 \\
    1e-2 & 0.47492741 & 0.73746015 \\
% -1.0100006851146353 -1.010000000000002
   %1e-3 & 0.241423605711876460 & 0.6203095737001850 \\
    1e-3 & 0.24142361 & 0.62030957 \\
% 1.0010080277991904 1.0010000000000392
   %1e-4 & 0.052669616178245525 & 0.5048726352785534 \\
    1e-4 & 0.05266962 & 0.50487264 \\
% 1.0001535423026533 1.0001000000001408
   %1e-5 & 0.005508556194265022 & 0.3897022288168669 \\
    1e-5 & 0.00550856 & 0.38970223 \\
% -1.000095855127284 1.0000100000016647
   %0.0 &  0.000000003689826620 & 0.0000014898442216 \\
    0.0 &  3.689827e-9 & 1.489844e-6 \\
% 1.0000907216405073 1.000000004122051
    \hline
  \end{tabular}
  \caption{Location of the transition layer $z_{\ast}$ found by solving Burgers' equation directly with \texttt{mystic} for select $v$ and $\delta$. Solving Burger's equation for a given $\delta$ and $v$ also produces $A$, where $A \simeq \pm(1.0 + \delta)$.}\label{tab:exactly}
\end{table}

%\END COMMENT

%   result = diffev2(objective, args=(v, delta), x0=bounds, bounds=bounds,\
%                    npop=40, ftol=1e-8, gtol=50, disp=False, full_output=True)
\begin{figure}[htbp!]
  \centering
  \begin{minipage}{1.00\textwidth}
% \begin{\outputtextsize}
\begin{verbatim}
### exact_supersensitive.py :: solve Burgers' equation for v and delta ###
from math import tanh
from mystic import reduced
from mystic.solvers import lattice

bounds = [(0,1),(None,None)]

@reduced(lambda x,y: abs(x)+abs(y))
def objective(x, v, delta):
  z,A = x
  return [A * tanh(0.5 * (A/v) * (1. - z)) - 1.,
          A * tanh(0.5 * (A/v) * (1. + z)) - 1. - delta]

def solve(v, delta):
  "solve for (z,A) in analytical solution to burgers equation"
  result = lattice(objective, ndim=2, nbins=4, args=(v, delta),\
                   bounds=bounds, ftol=1e-8, disp=False, full_output=True)
  return delta,result[0],result[1]

if __name__ == '__main__':

  for v in [0.1, 0.05]:
    for delta in [1e-1, 1e-2, 1e-3, 1e-4, 1e-5, 0.]:
      _, result, fit = solve(v, delta)
      print "v, delta:", (v, delta)
      print "z, A:", result.tolist()
      print fit, "\n"

\end{verbatim}
% \end{\outputtextsize}
  \end{minipage}
  \caption{Solving Burgers' equation directly with \texttt{mystic} for several $v$ and $\delta$. The \texttt{lattice} solver starts \texttt{nbins} Nelder-Mead optimizers at the center-points of a uniform grid instead of requiring a single initial value for \texttt{(z,A)}. A \texttt{lambda} function is used to reduce the output of the \texttt{objective} to a single value. Results are summarized in Table \ref{tab:exactly}.}\label{code:exactly}
\end{figure}

% \newpage

\section{Monte Carlo Sampling of Transition Layer Futures}

\begin{figure}[htbp!]
  \centering
  \begin{minipage}{1.00\textwidth}
  \subfigure[$\epsilon = 0.1$, $v = 0.1$]{
    \label{fig:v01_eps01_cdf}
    \includegraphics[width=0.30\textwidth]{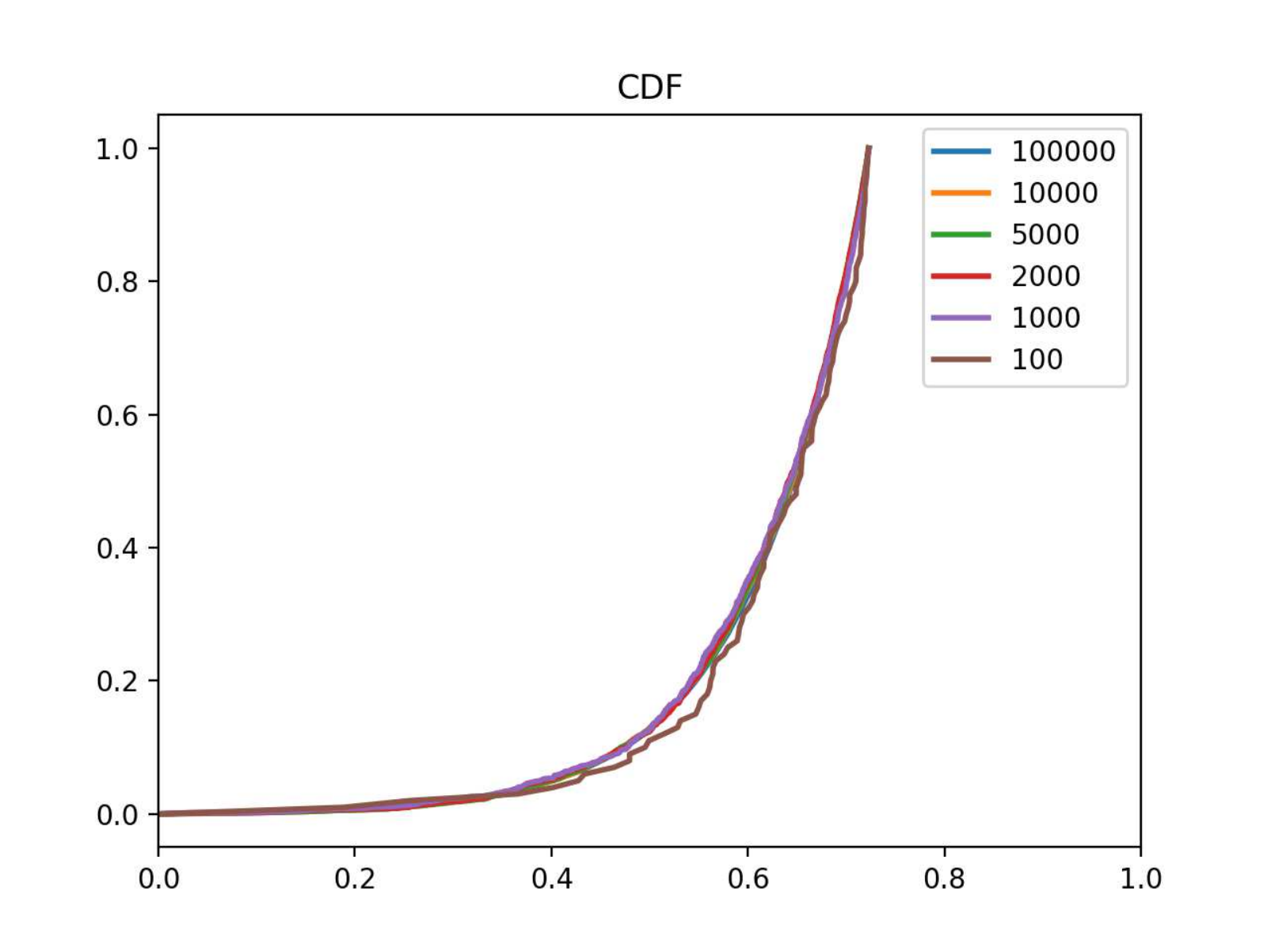}
  }
  ~
  \subfigure[$\epsilon = 0.01$, $v = 0.1$]{
    \label{fig:v01_eps001_cdf}
    \includegraphics[width=0.30\textwidth]{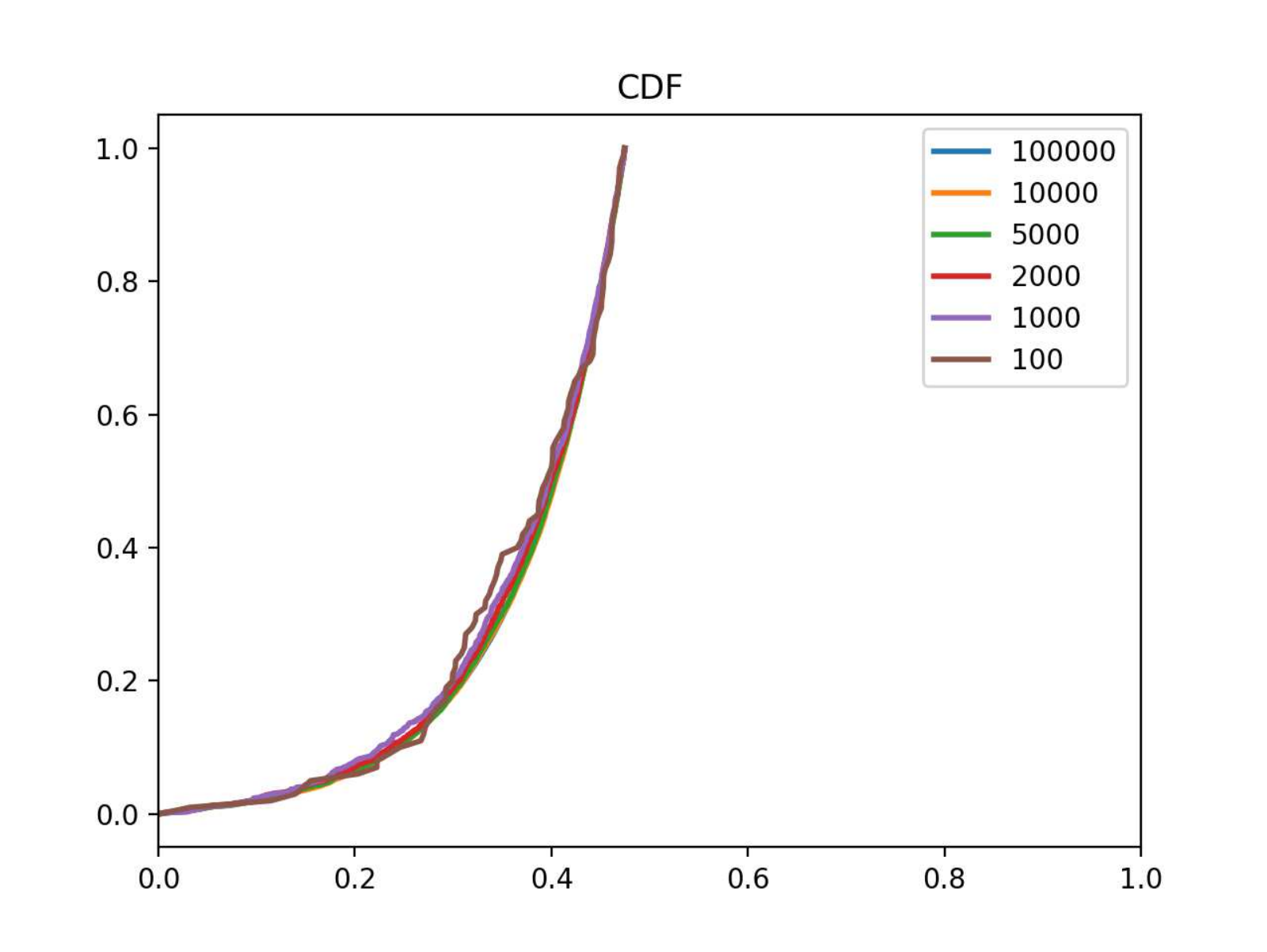}
  }
  ~
  \subfigure[$\epsilon = 0.001$, $v = 0.1$]{
    \label{fig:v01_eps0001_cdf}
    \includegraphics[width=0.30\textwidth]{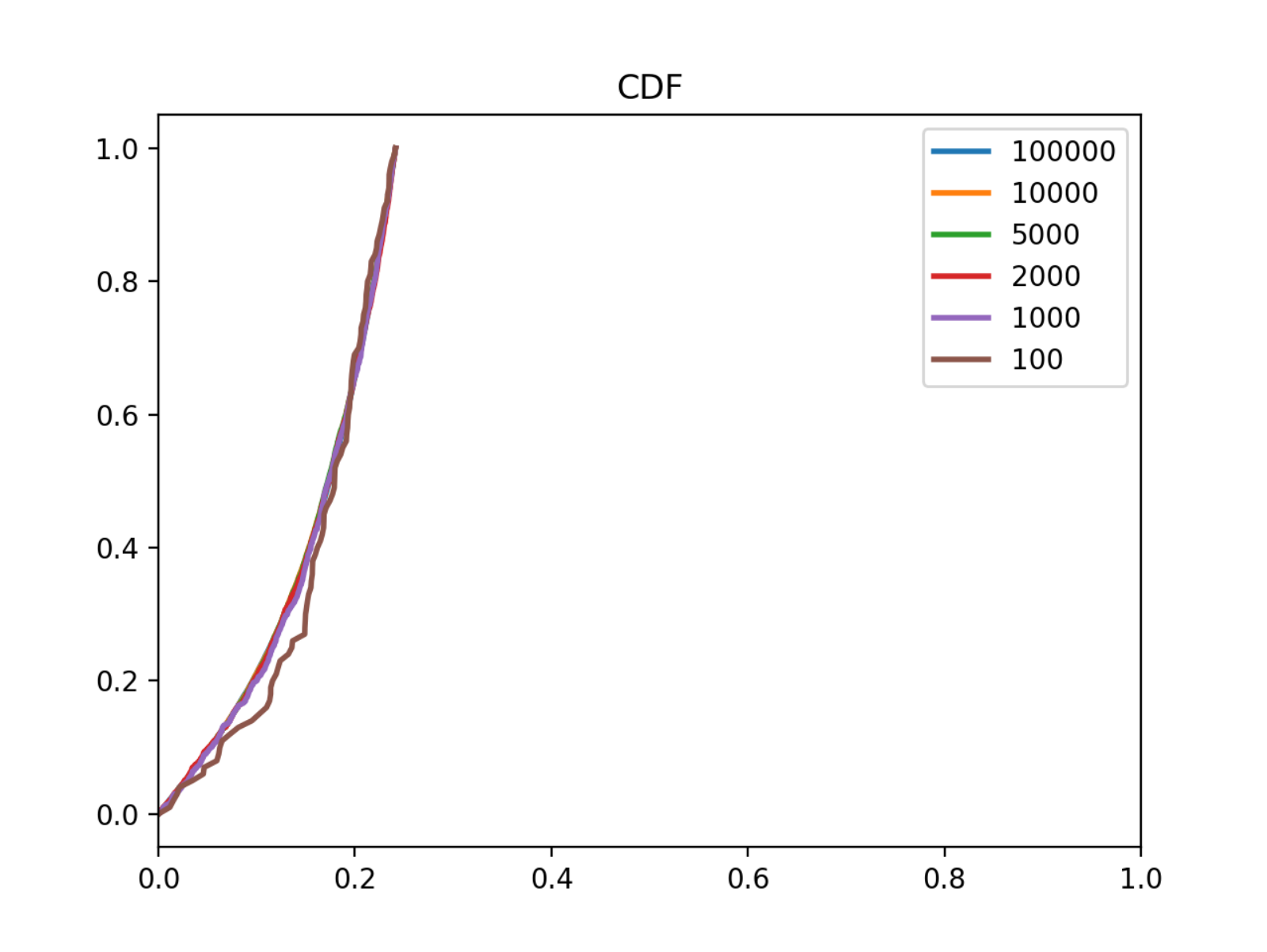}
  }
  \subfigure[$\epsilon = 0.1$, $v = 0.05$]{
    \label{fig:v005_eps01_cdf}
    \includegraphics[width=0.30\textwidth]{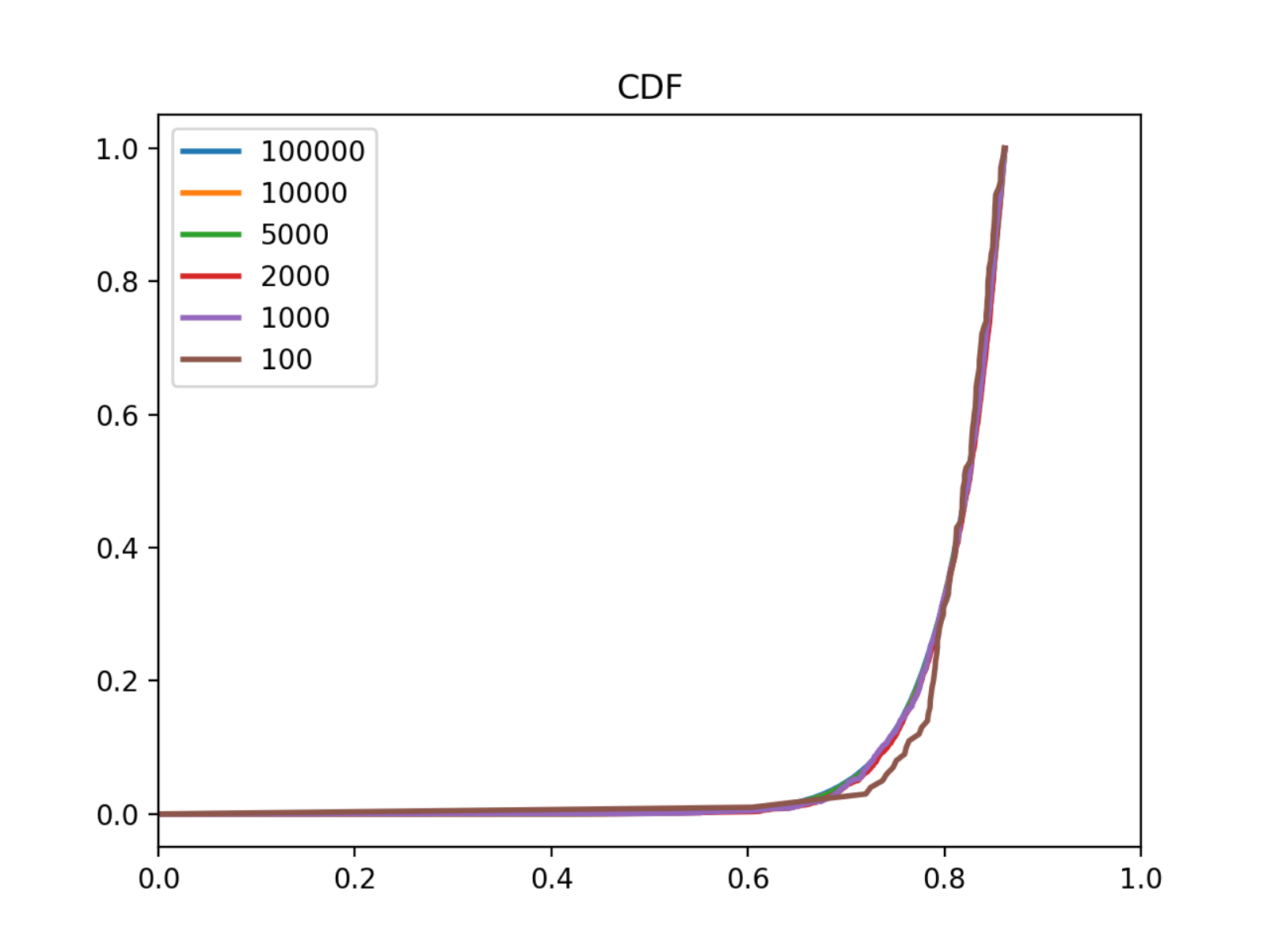}
  }
  ~
  \subfigure[$\epsilon = 0.01$, $v = 0.05$]{
    \label{fig:v005_eps001_cdf}
    \includegraphics[width=0.30\textwidth]{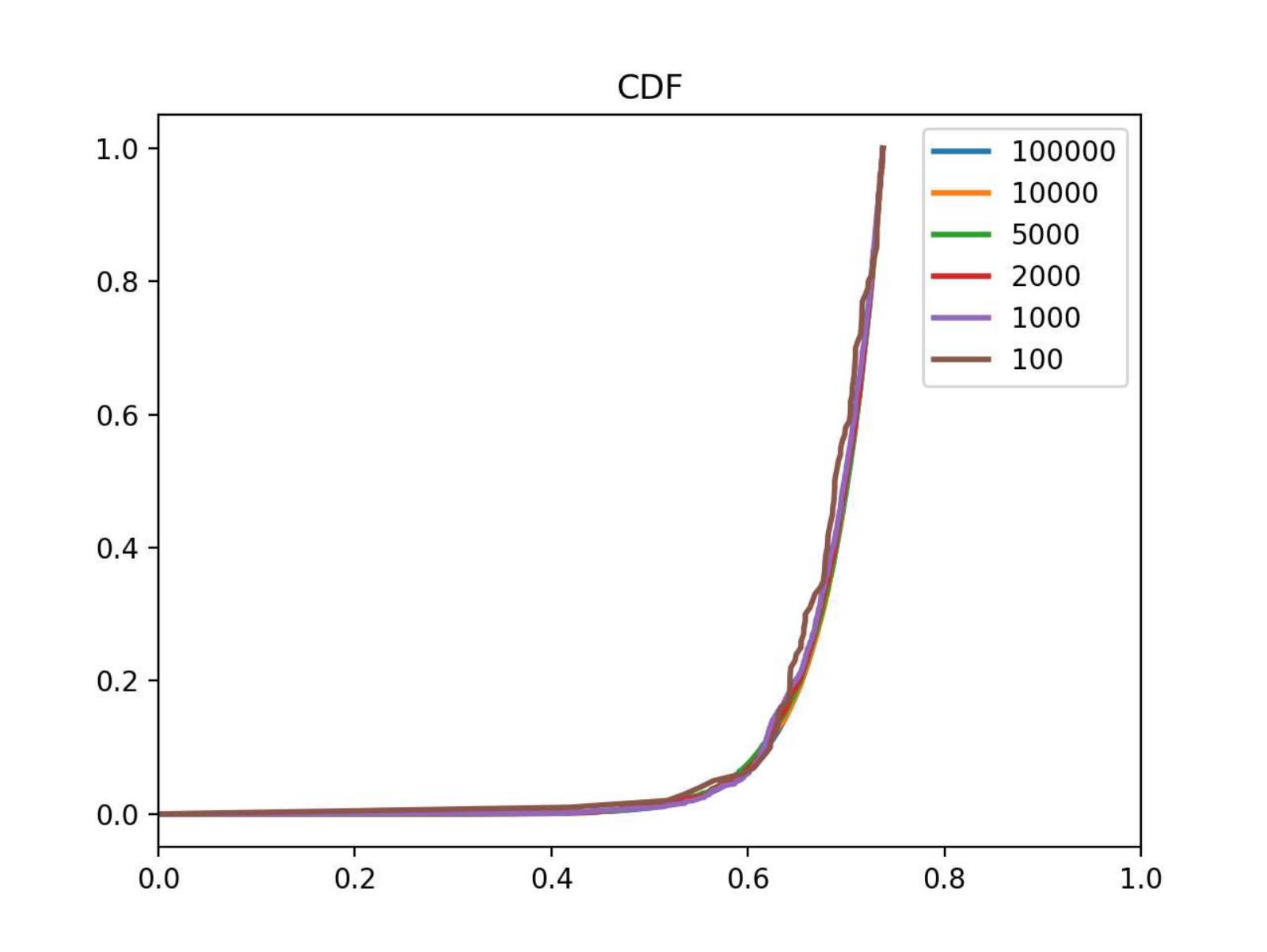}
  }
  ~
  \subfigure[$\epsilon = 0.001$, $v = 0.05$]{
    \label{fig:v005_eps0001_cdf}
    \includegraphics[width=0.30\textwidth]{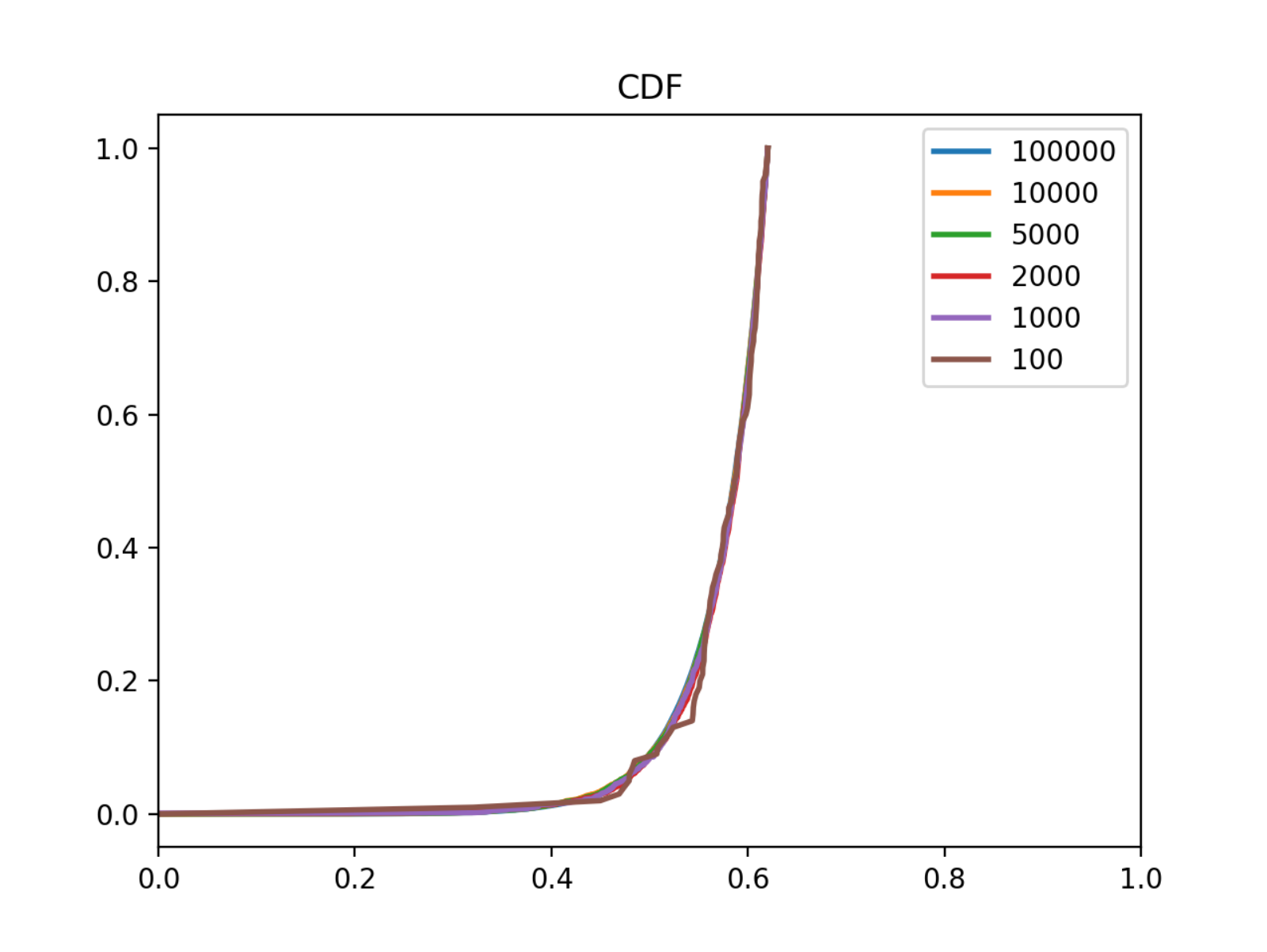}
  }
  \caption{Cumulative density functions of $z_{\ast}$ for selected $\delta \sim U(0,\epsilon)$ and $v$, calculated at different numbers of realizations $N$ of $\delta$. Shown are $N = 100,1000,2000,5000,10000,100000$, with $N=100$ in brown and $N=100000$ in blue. The plots show change in $v$ and $\epsilon$ effects the moments of $z_{\ast}$. More quantitative results are found in Table \ref{tab:multiepsilon}.}\label{fig:cdf}
  \end{minipage}
\end{figure}

\begin{figure}[htbp!]
  \centering
  \begin{minipage}{1.00\textwidth}
  \subfigure[$\epsilon = 0.1$, $v = 0.1$]{
    \label{fig:v01_eps01_icdf}
    \includegraphics[width=0.30\textwidth]{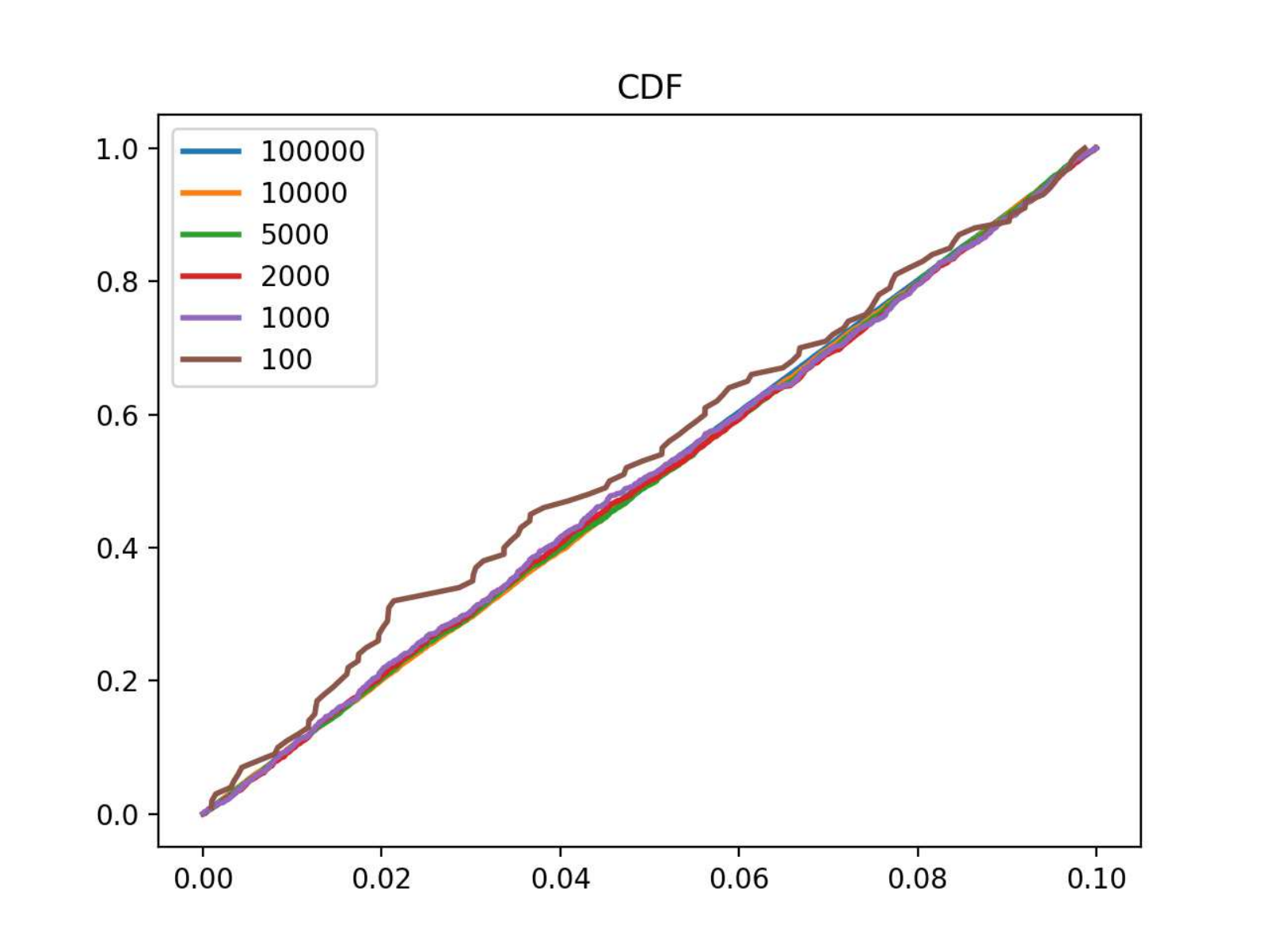}
  }
  ~
  \subfigure[$\epsilon = 0.01$, $v = 0.1$]{
    \label{fig:v01_eps001_icdf}
    \includegraphics[width=0.30\textwidth]{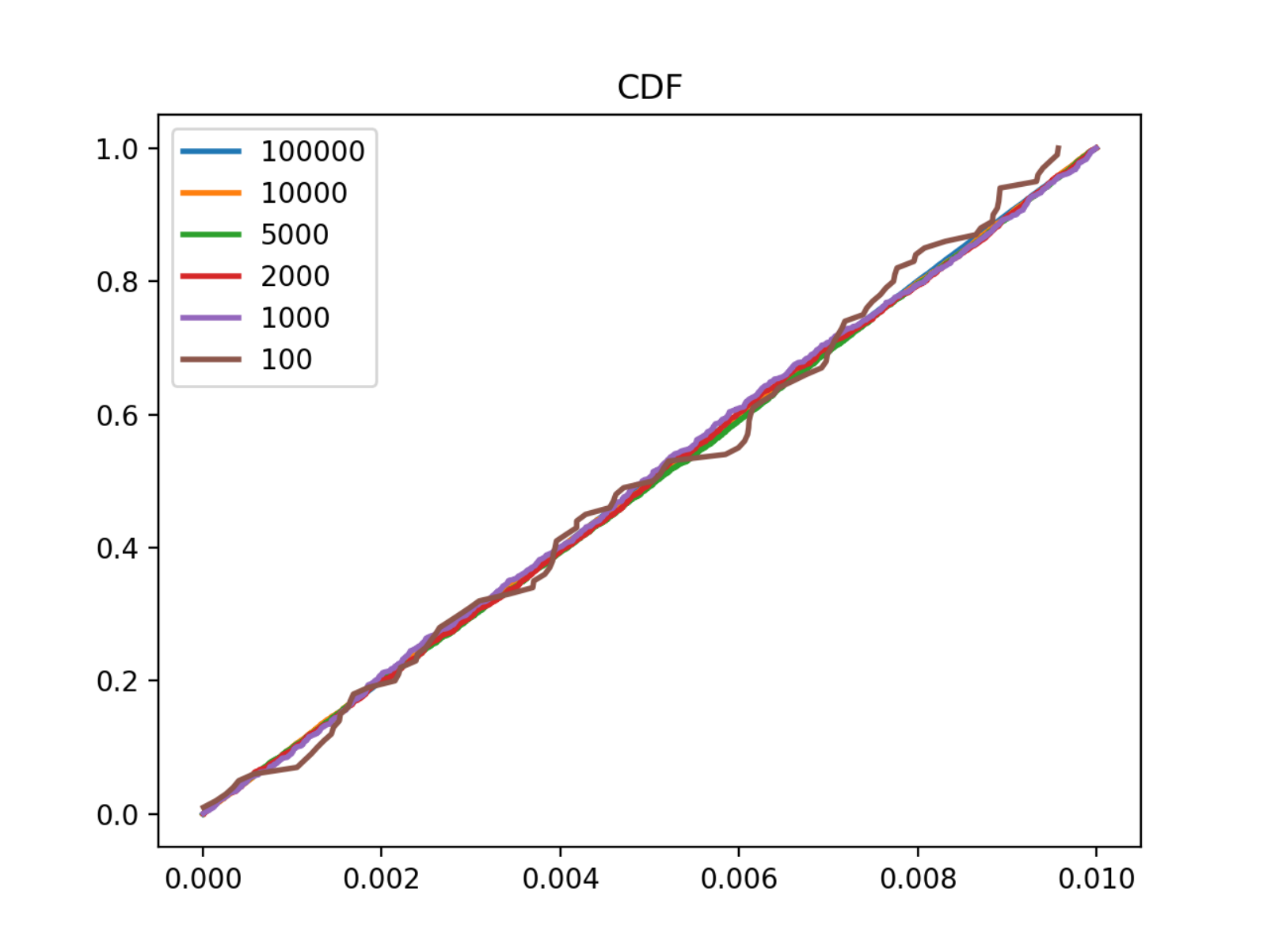}
  }
  ~
  \subfigure[$\epsilon = 0.001$, $v = 0.1$]{
    \label{fig:v01_eps0001_icdf}
    \includegraphics[width=0.30\textwidth]{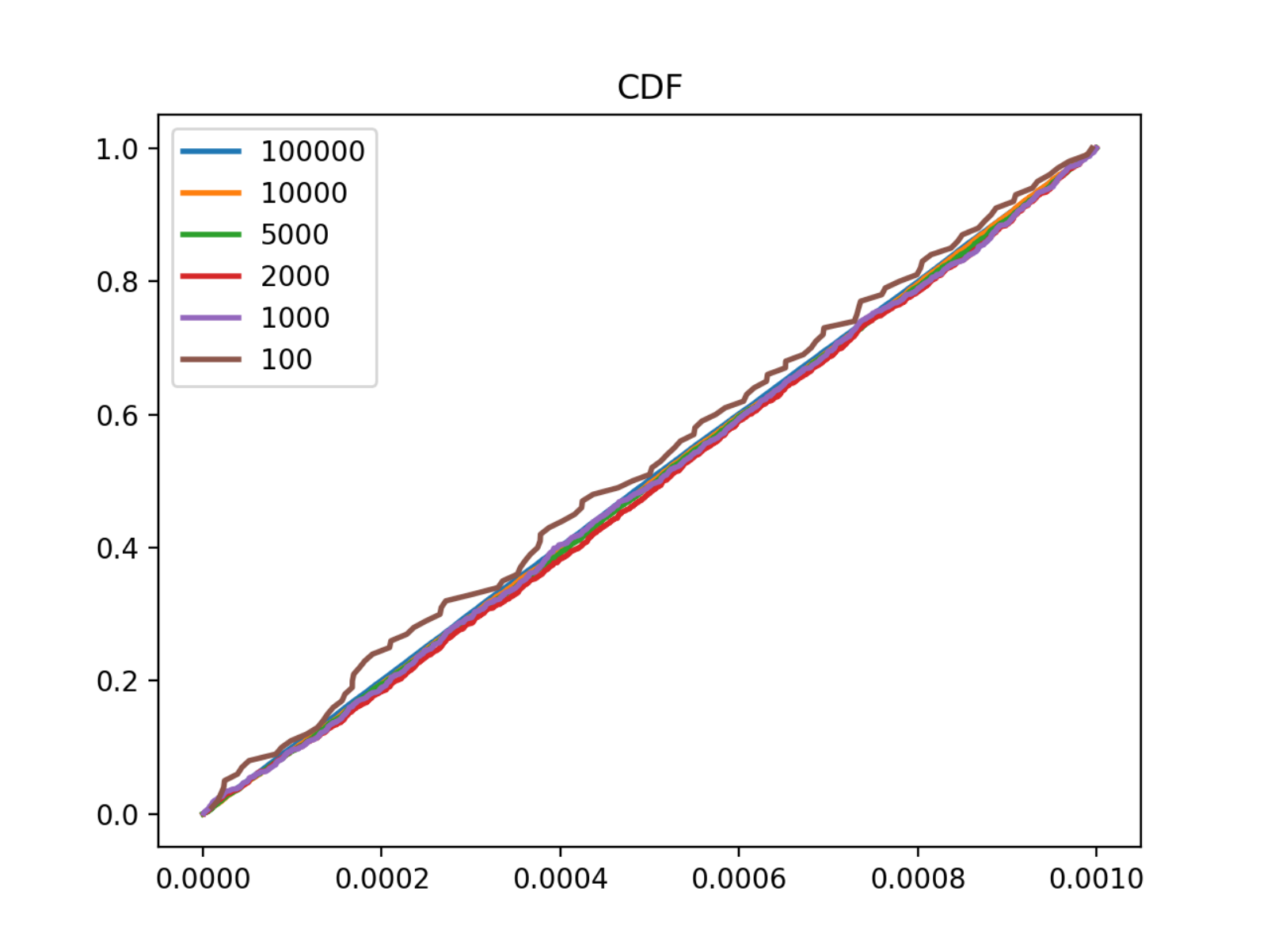}
  }
  \subfigure[$\epsilon = 0.1$, $v = 0.05$]{
    \label{fig:v005_eps01_icdf}
    \includegraphics[width=0.30\textwidth]{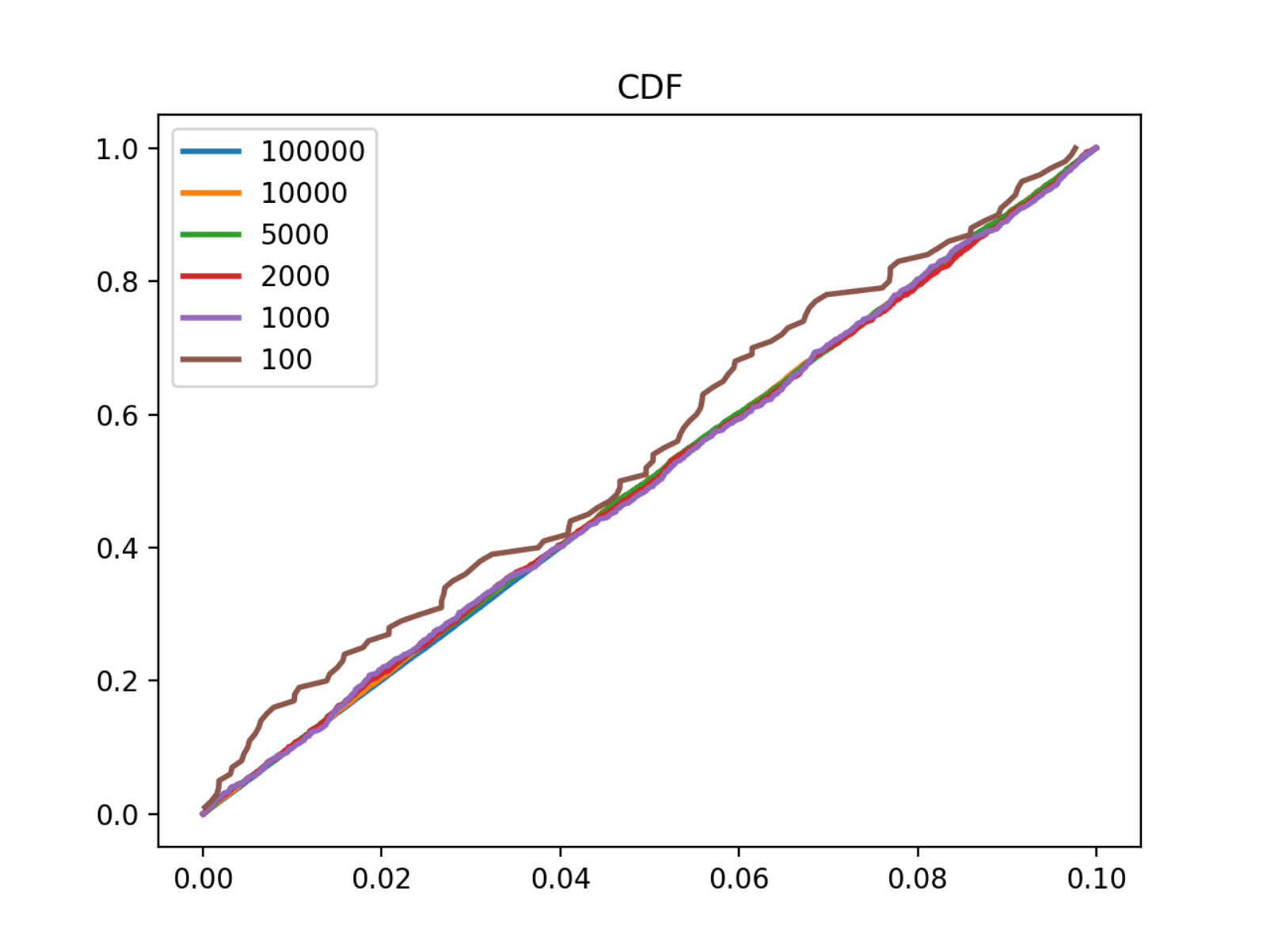}
  }
  ~
  \subfigure[$\epsilon = 0.01$, $v = 0.05$]{
    \label{fig:v005_eps001_icdf}
    \includegraphics[width=0.30\textwidth]{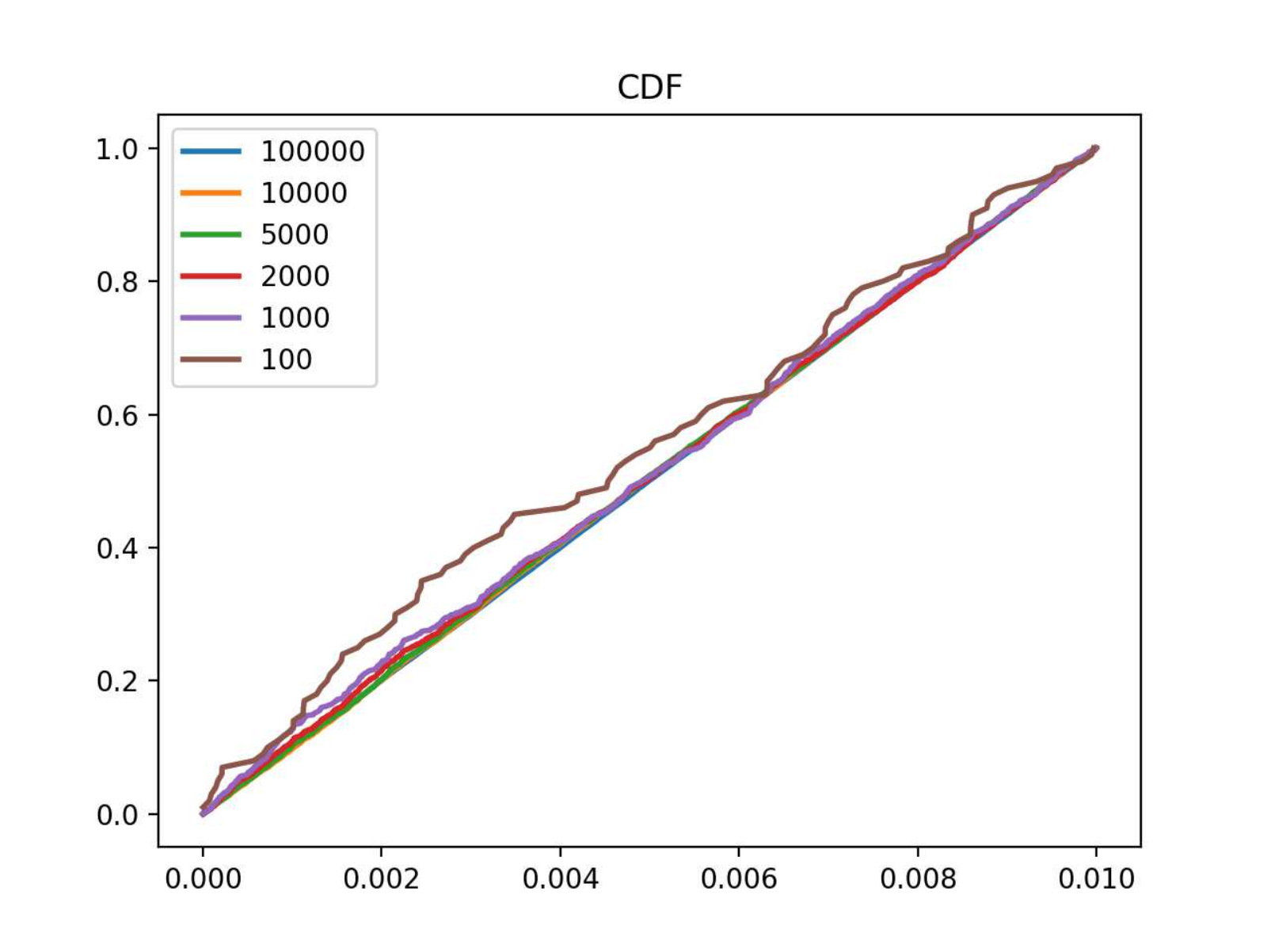}
  }
  ~
  \subfigure[$\epsilon = 0.001$, $v = 0.05$]{
    \label{fig:v005_eps0001_icdf}
    \includegraphics[width=0.30\textwidth]{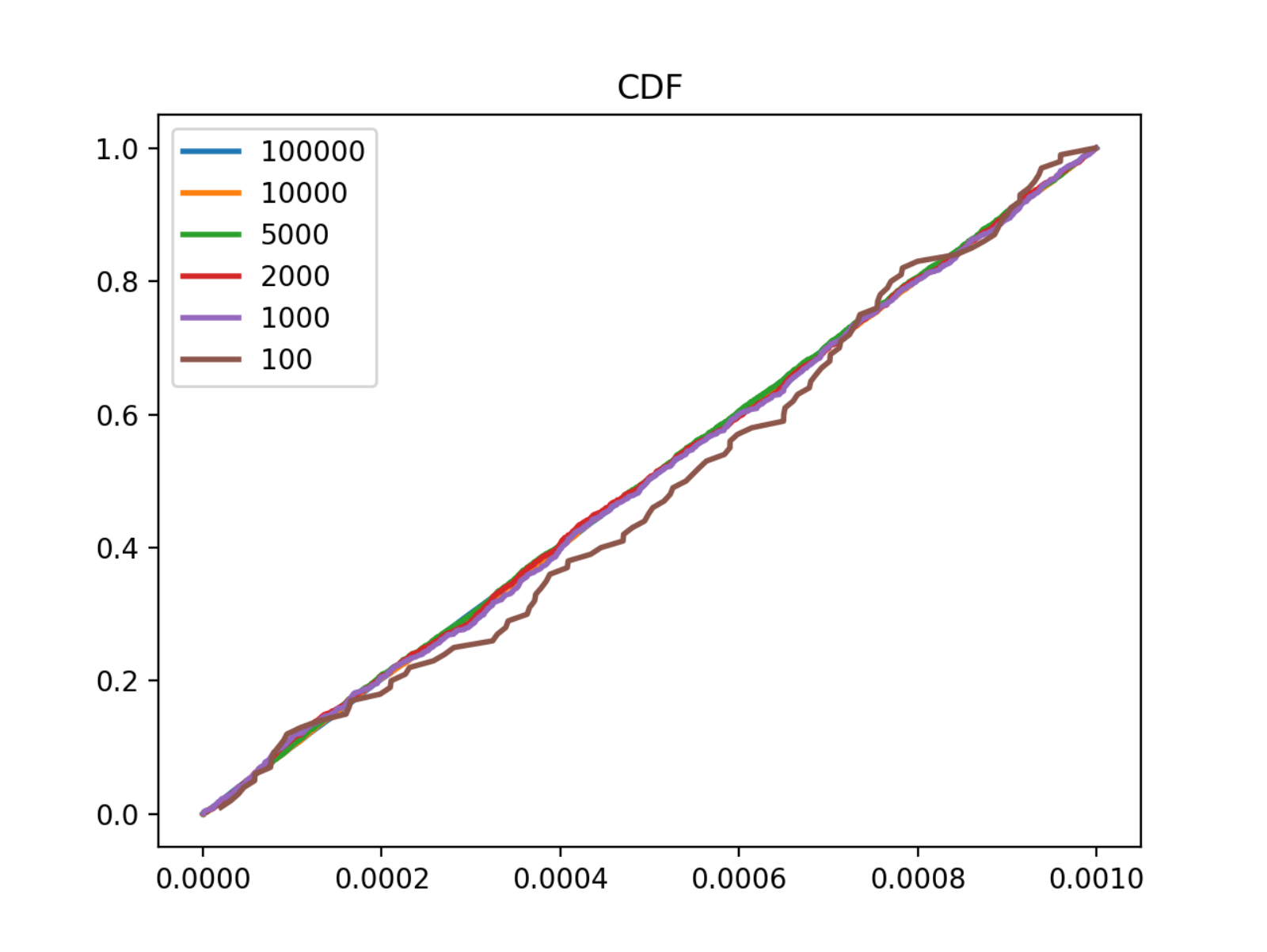}
  }
  \caption{Cumulative density functions of $\delta \sim U(0,\epsilon)$ for selected $\epsilon$ and $v$, calculated at different numbers of realizations $N$ of $\delta$. Shown are $N = 100,1000,2000,5000,10000,100000$. Sampling larger $N$ provides distributions that better approximate the theoretical \texttt{cdf}.}\label{fig:icdf}
  \end{minipage}
\end{figure}

If we assume the location of the left wall $\delta$ is a uniform random variable in $(0,\epsilon)$:
\begin{equation*}
  \delta \sim U(0,\epsilon), 
\end{equation*}
we can then repeatedly solve for location of the transition layer $z_{\ast}$ where $u(z_{\ast}) = 0$, using Monte Carlo sampling. By inspecting the cumulative distribution function \texttt{cdf} of solved $z_{\ast}$ (Figure \ref{fig:cdf}), we can get an idea of the likelihood of $z_{\ast}$ occurring at location $z$.  A more quantitative measure of likelihood can be obtained by calculating the mean location $\bar{z_{\ast}}$ and standard deviation $\sigma_{z_{\ast}}$ of the futures of $z$.
Plots of the \texttt{cdf} of $\delta$ (Figure \ref{fig:icdf}) can be used to validate that our inputs were sampled from a distribution closely resembling the theoretical \texttt{cdf}, and also helps ensure we select an appropriate number of futures $N$ in the calculation of likelihood of success, $P(z_{\ast} > \frac{x\bar{z_{\ast}}}{100})$, where we define ``success'' as when the transition occurs at a $z_{\ast}$ greater than $\frac{x\bar{z_{\ast}}}{100}$ for a given scaling $x$.

The code in Figure \ref{code:sampling} details how we 
sample $M=110000$ instances of $\delta$ from $U(0,\epsilon)$ for each combination of $v$ and $\epsilon$. The $M$ samplings of $\delta$ are passed into \texttt{solve}, and produce $N$ futures of $z$. An unordered iterative map \texttt{uimap} from \texttt{pathos.pools.ProcessPool} provides embarrassingly-parallel invocations of \texttt{solve}.  Although we start with $M$ samplings, solutions of \texttt{solve} with convergence less than \texttt{tol} = $1e^{-9}$ are discarded, and results are calculated until we have $N$ solutions of \eqref{eq:burgers_steady} within tolerance \texttt{tol}.  The $N$ futures of $z$, and corresponding $\delta$, for each $(v,\epsilon)$, are then sorted and saved to a file. 
Discarding solutions that fail to converge within \texttt{tol} does not appear to bias the results, and \texttt{solve} fails to achieve the desired tolerance \texttt{tol} in very few of the $M$ attempts.
For example, with $N=100000$, $v = 0.05$, and $\epsilon = 0.01$,
%only $668$ of $M=110000$ ($\sim 0.6\%$)
only $3$ of $M=110000$ ($\sim 2.7e^{-3}\%$)
possible invocations of \texttt{solve} fail to achieve the desired convergence \texttt{tol}.
In going to much larger $N$, we find that the normalized cumulative distribution of discarded $\delta$ approximates the sampling distribution for all combinations of $(v,\epsilon)$, thus discarding the poorly converged solutions does not appear have any bias on the results.
%In going to much larger $N$, we find that the distribution of the discarded $\delta$ is roughly uniform in $[0,\epsilon]$, thus no bias is caused by dropping solutions that do not converge within \texttt{tol}. %(see Figure \ref{fig:misses_u}).
For the above case, results for $N$ futures of $z$ were obtained in
%\texttt{02:39:10.132864}
\texttt{2:09:42.320102}
using 8-way parallel on the above-specified Macbook.
Timings and number of discarded $\delta$ were similar for all combinations of $v$ and $\epsilon$ (as seen in Table \ref{tab:misses}).

%%% to be used in calculating $\bar{z_{\ast}}$, $\sigma_{z_{\ast}}$, $\bar{\delta}$, $\sigma_{\delta}$, \texttt{cdf}($z_{\ast}$), and \texttt{cdf}($\delta$). The futures of $z$ will also be used in calculations of the probability of success.

The choice of the sampling distribution $U(0,\epsilon)$ is somewhat arbitrary.
For example, we can instead assume $\delta$ is a random variable with a truncated Gaussian distribution:
\begin{equation*}
  \delta \sim G(0,\epsilon,\epsilon/2,\sqrt{\epsilon^{2}/12}),
\end{equation*}
where $\epsilon/2$ and $\sqrt{\epsilon^{2}/12}$ are, respectively, the theoretical mean and standard deviation of $\delta$. As seen in the code in Figure \ref{code:rtnorm}, we will choose the range, mean, and standard deviation of the truncated Gaussian distribution to be the same as those of the uniform distribution $U(0,\epsilon)$. Additionally, choosing identical moments will better facilitate more direct comparisons of results between the two distributions.
We again discard solutions that fail to converge within \texttt{tol}, and again note that the discarded $\delta$ do not appear to bias the results.
%and we again see (in Figure \ref{fig:misses_n}) that the discarded $\delta$ do not appear to bias the results.
Sampling from $\delta \sim G(0,\epsilon,\epsilon/2,\sqrt{\epsilon^{2}/12})$, with $N=100000$, $v = 0.1$, and $\epsilon = 0.01$, using the above-specified Macbook, provides results in
%\texttt{02:42:27.513162},
\texttt{2:02:04.456247},
where again only
%$668$
$2$
of $M=110000$ invocations of \texttt{solve} fail to meet the desired convergence.  As with $U(0,\epsilon)$, timings and the number of discarded $\delta$ are generally similar for all combinations of $v$ and $\epsilon$.

\begin{figure}[htbp!]
  \centering
  \begin{minipage}{1.00\textwidth}
  \subfigure[$\epsilon = 0.1$, $v = 0.1$]{
    \label{fig:v01_eps01_cdf_rt}
    \includegraphics[width=0.30\textwidth]{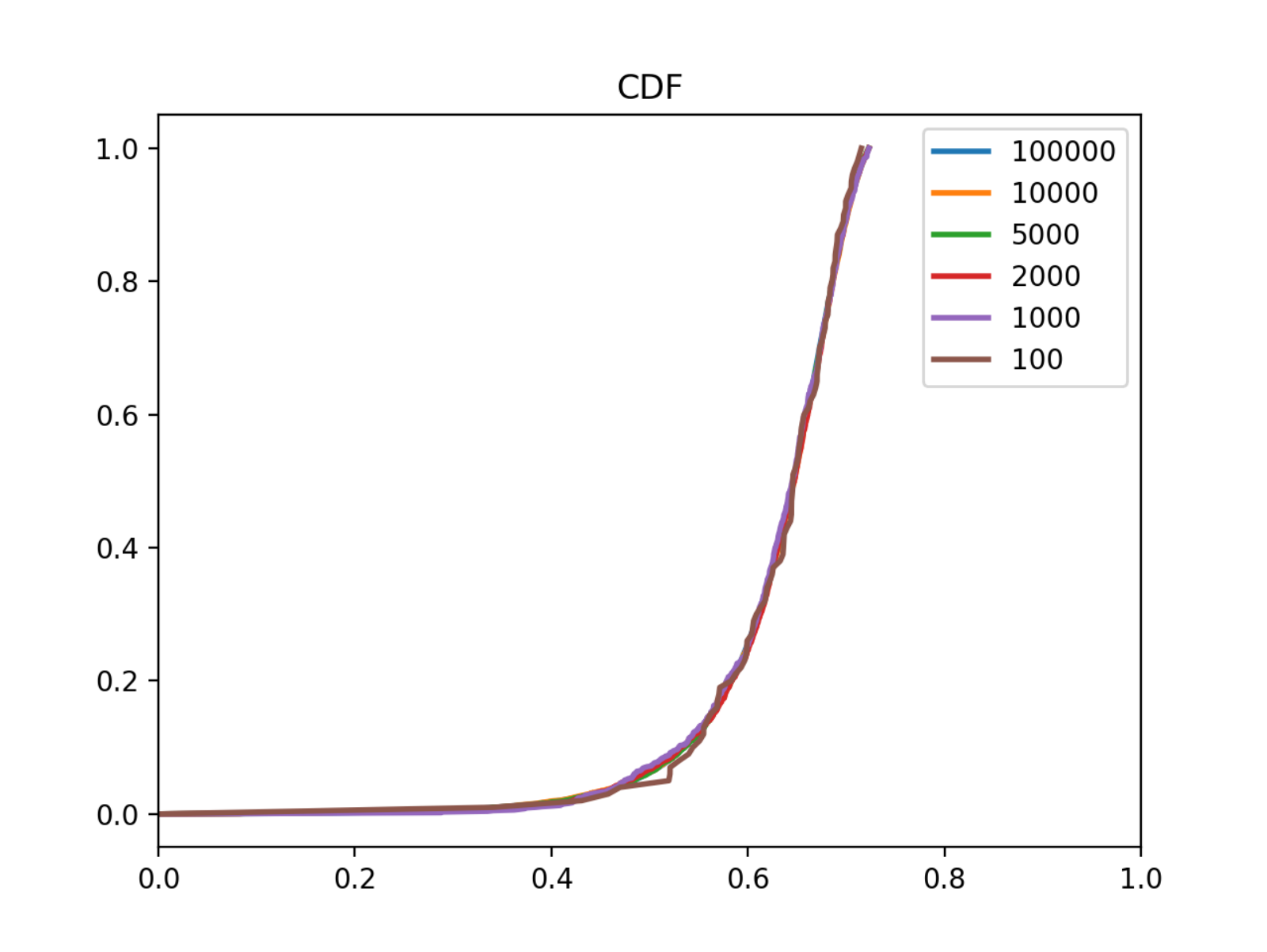}
  }
  ~
  \subfigure[$\epsilon = 0.01$, $v = 0.1$]{
    \label{fig:v01_eps001_cdf_rt}
    \includegraphics[width=0.30\textwidth]{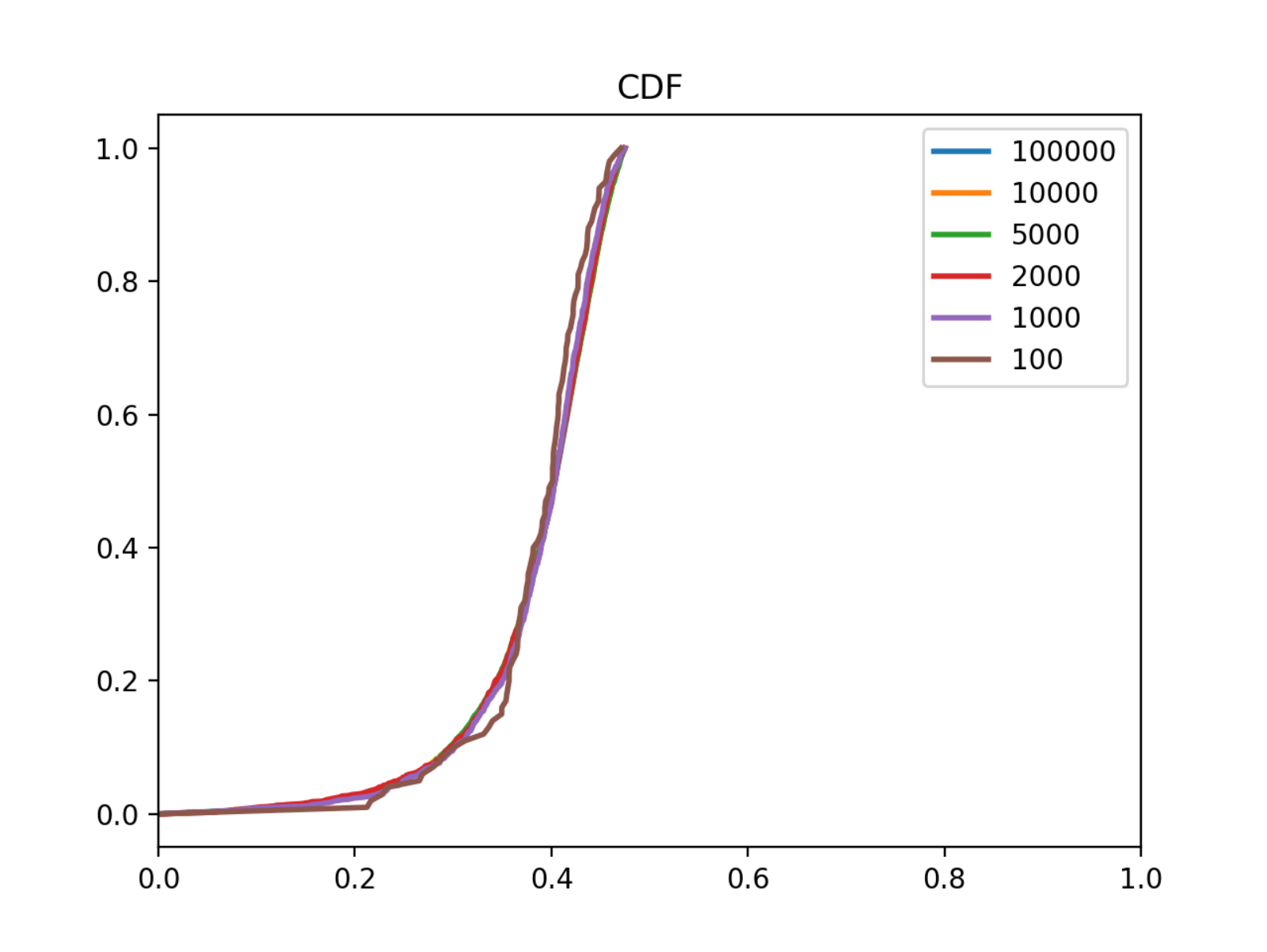}
  }
  ~
  \subfigure[$\epsilon = 0.001$, $v = 0.1$]{
    \label{fig:v01_eps0001_cdf_rt}
    \includegraphics[width=0.30\textwidth]{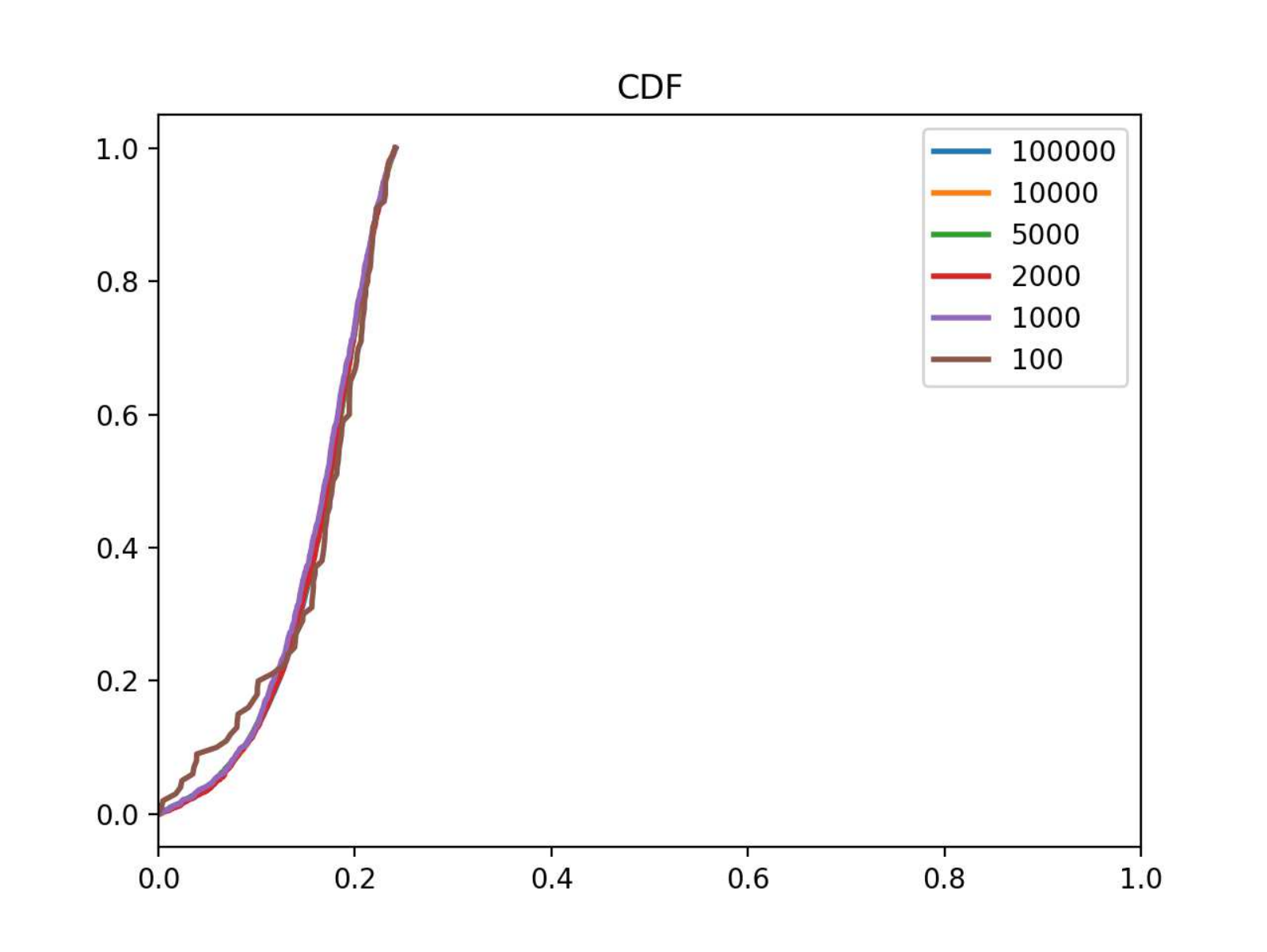}
  }
  \subfigure[$\epsilon = 0.1$, $v = 0.05$]{
    \label{fig:v005_eps01_cdf_rt}
    \includegraphics[width=0.30\textwidth]{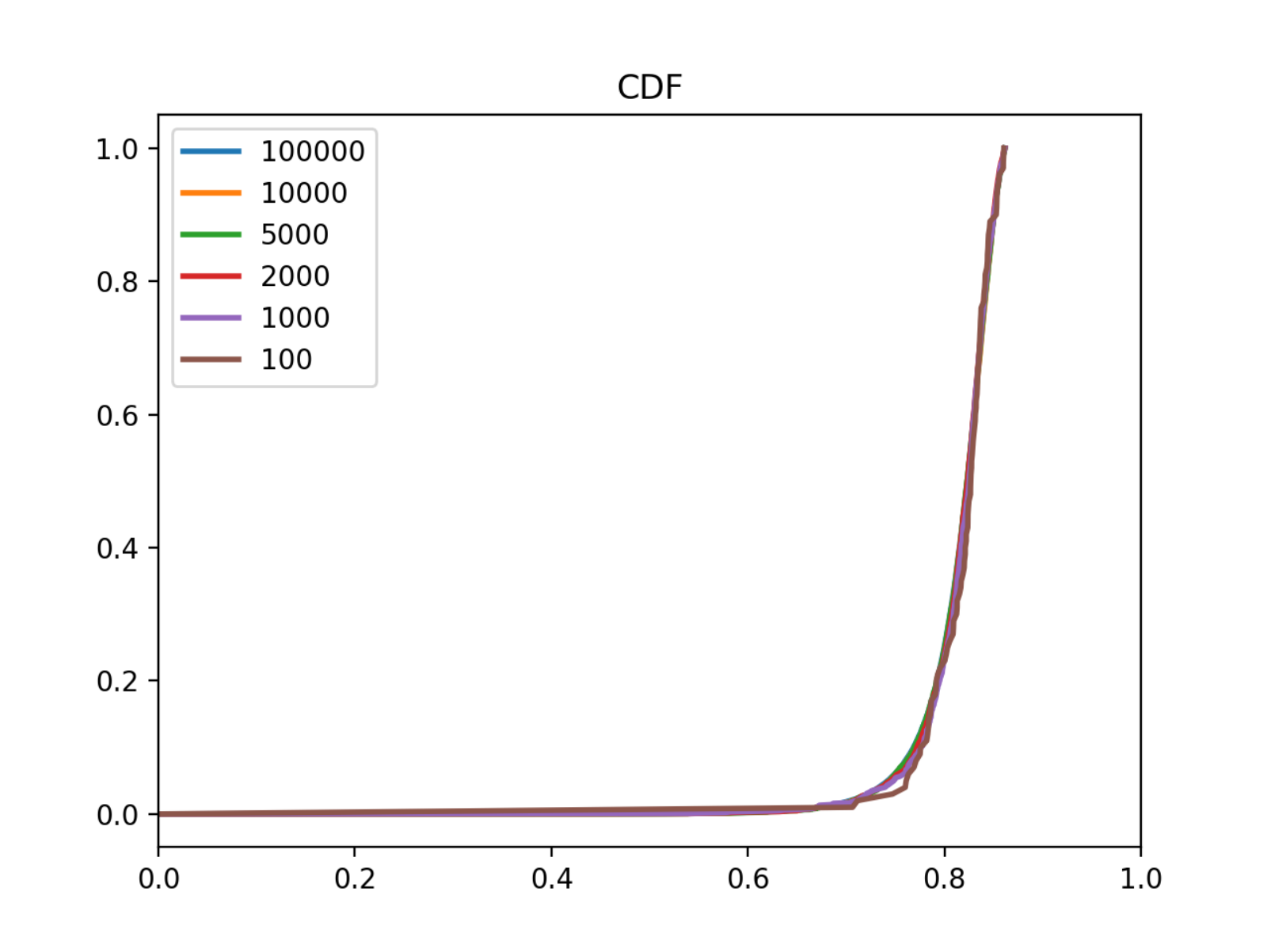}
  }
  ~
  \subfigure[$\epsilon = 0.01$, $v = 0.05$]{
    \label{fig:v005_eps001_cdf_rt}
    \includegraphics[width=0.30\textwidth]{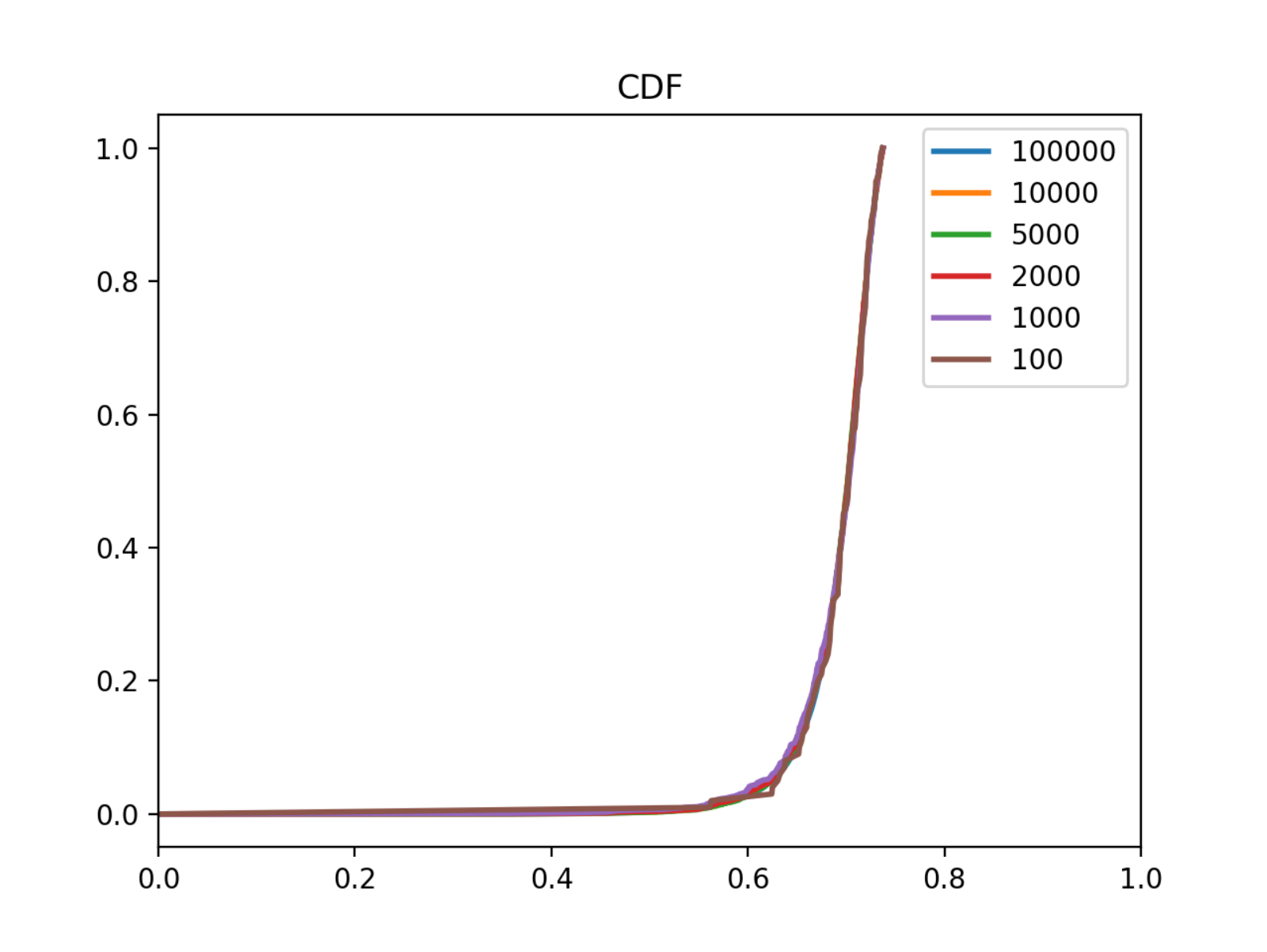}
  }
  ~
  \subfigure[$\epsilon = 0.001$, $v = 0.05$]{
    \label{fig:v005_eps0001_cdf_rt}
    \includegraphics[width=0.30\textwidth]{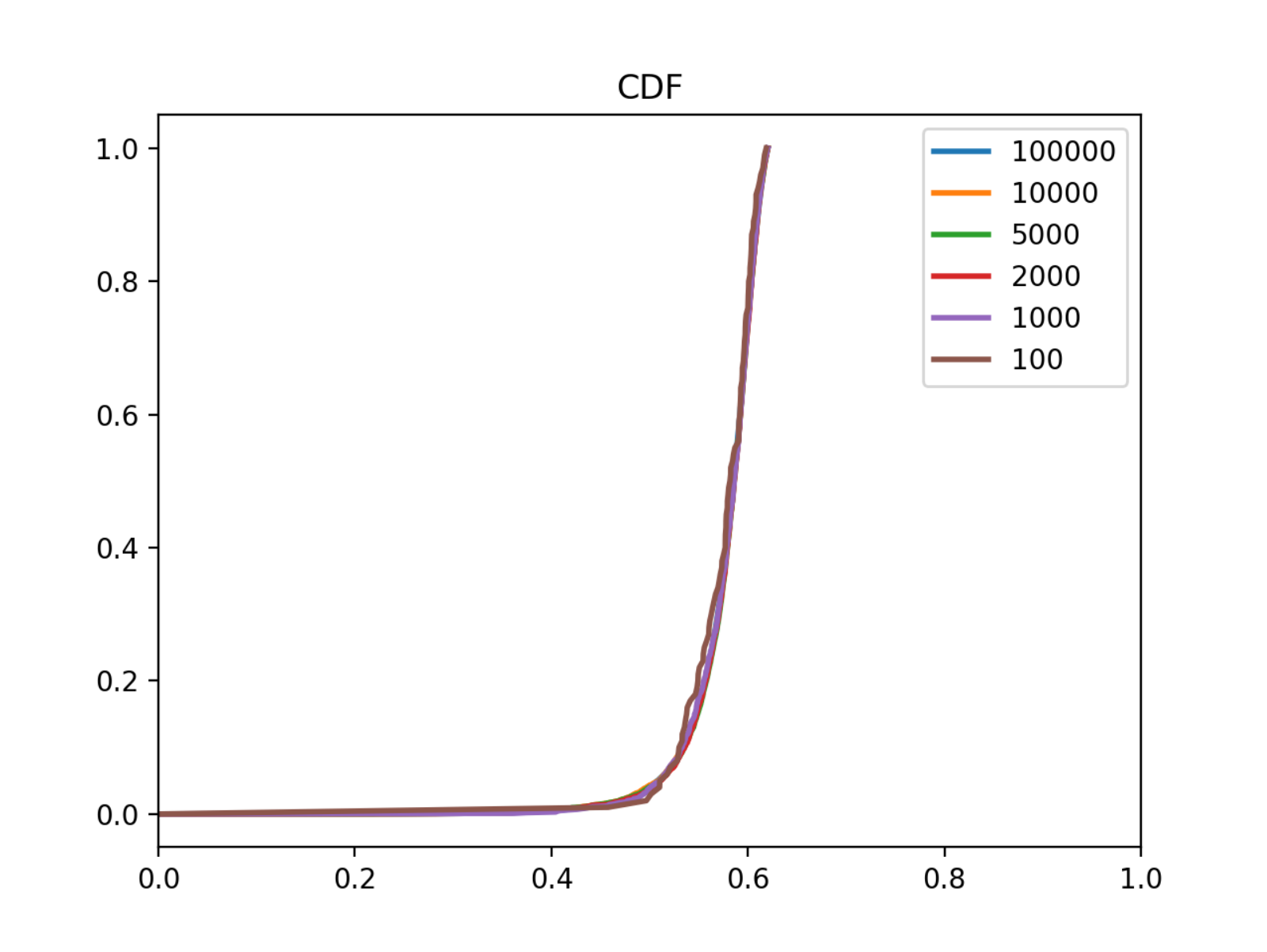}
  }
  \caption{Cumulative density functions of $z_{\ast}$ for selected $\delta \sim G(0,\epsilon,\epsilon/2,\sqrt{\epsilon^{2}/12})$ and $v$, calculated at different numbers of realizations $N$ of $\delta$. Shown are $N = 100,1000,2000,5000,10000,100000$. The plots show change in $v$ and $\epsilon$ effects the moments of $z_{\ast}$. More quantitative results are found in Table \ref{tab:multiepsilon}.}\label{fig:cdf_n}
  \end{minipage}
\end{figure}

\begin{figure}[htbp!]
  \centering
  \begin{minipage}{1.00\textwidth}
  \subfigure[$\epsilon = 0.1$, $v = 0.1$]{
    \label{fig:v01_eps01_icdf_rt}
    \includegraphics[width=0.30\textwidth]{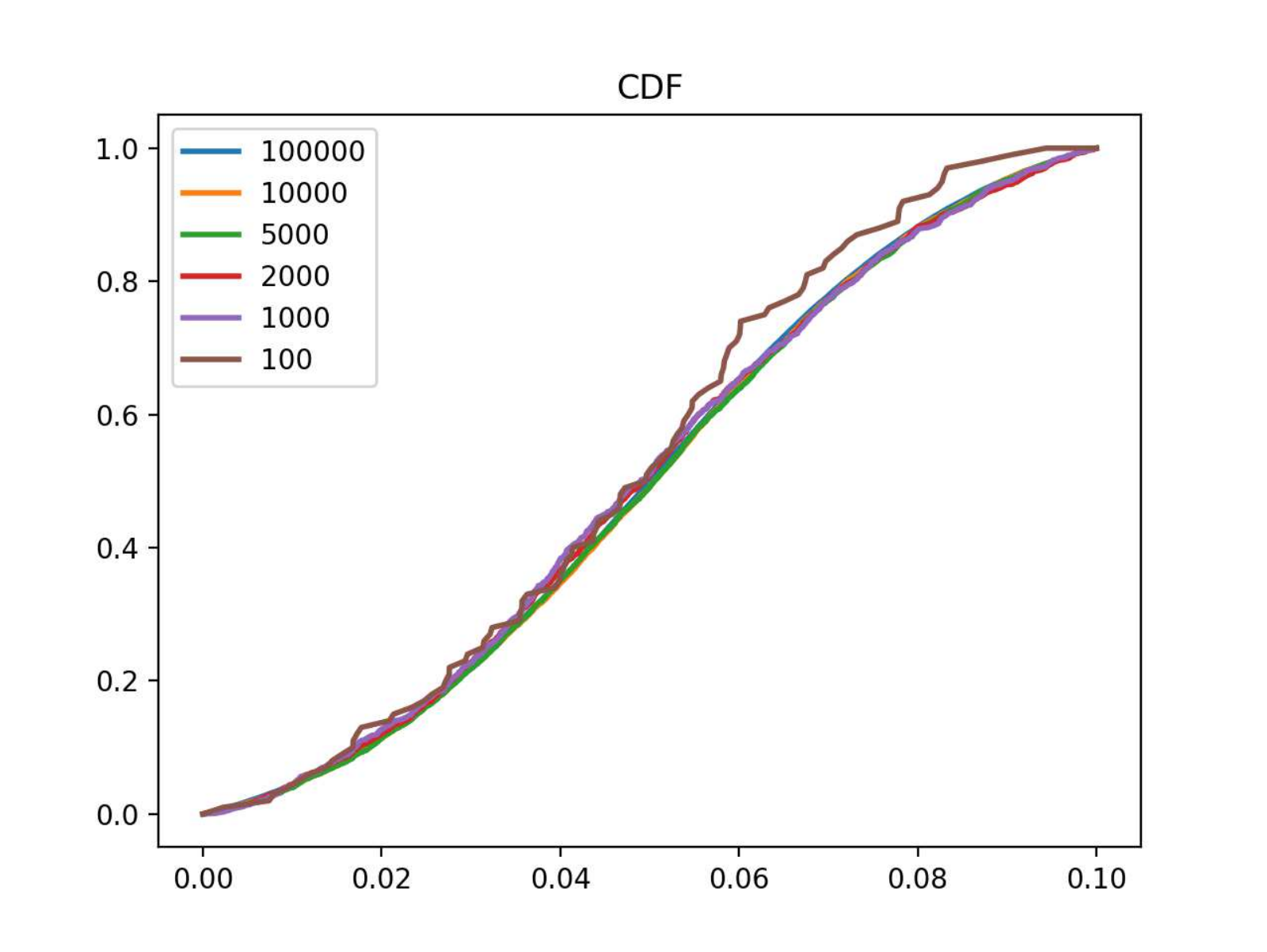}
  }
  ~
  \subfigure[$\epsilon = 0.01$, $v = 0.1$]{
    \label{fig:v01_eps001_icdf_rt}
    \includegraphics[width=0.30\textwidth]{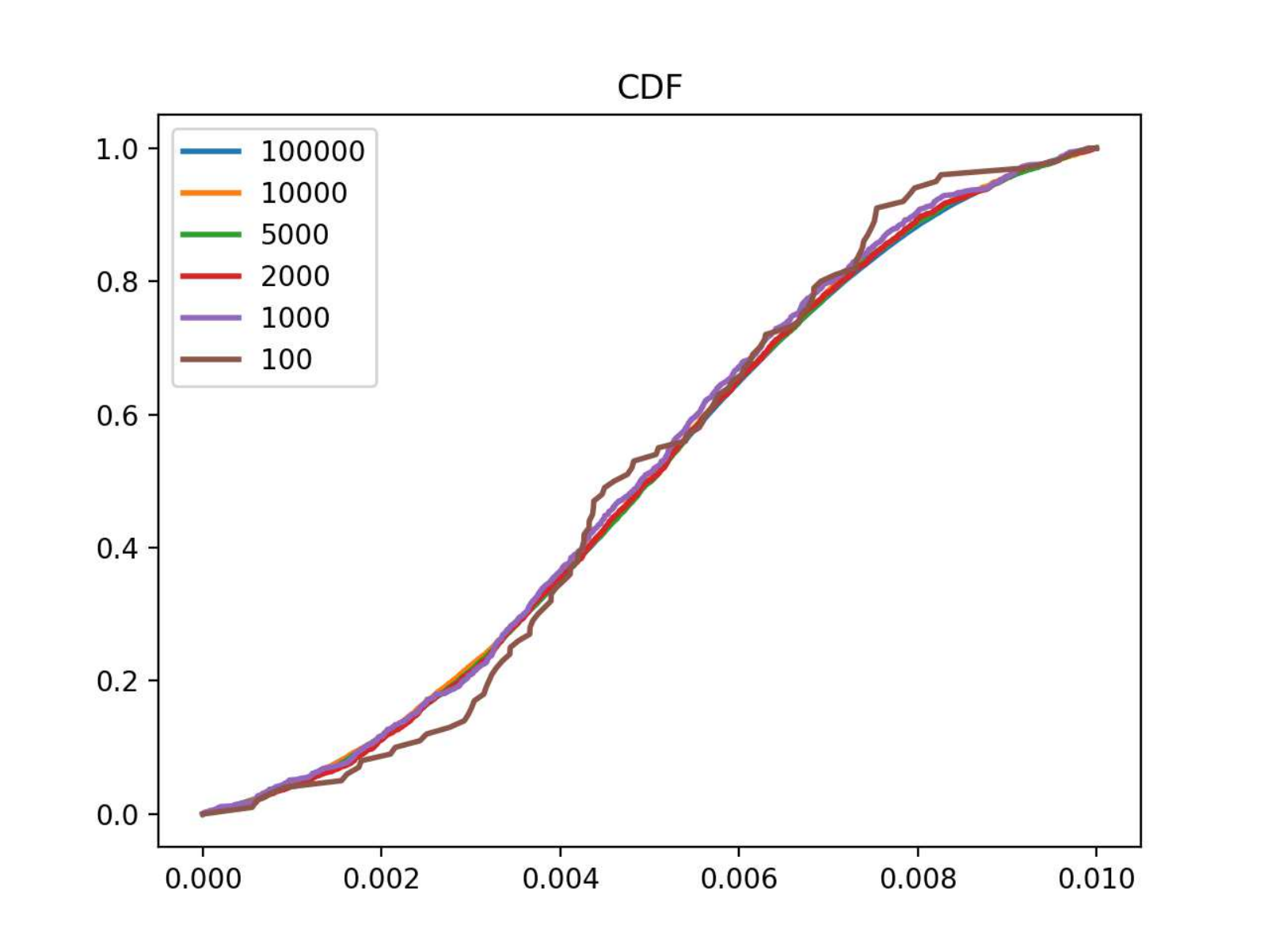}
  }
  ~
  \subfigure[$\epsilon = 0.001$, $v = 0.1$]{
    \label{fig:v01_eps0001_icdf_rt}
    \includegraphics[width=0.30\textwidth]{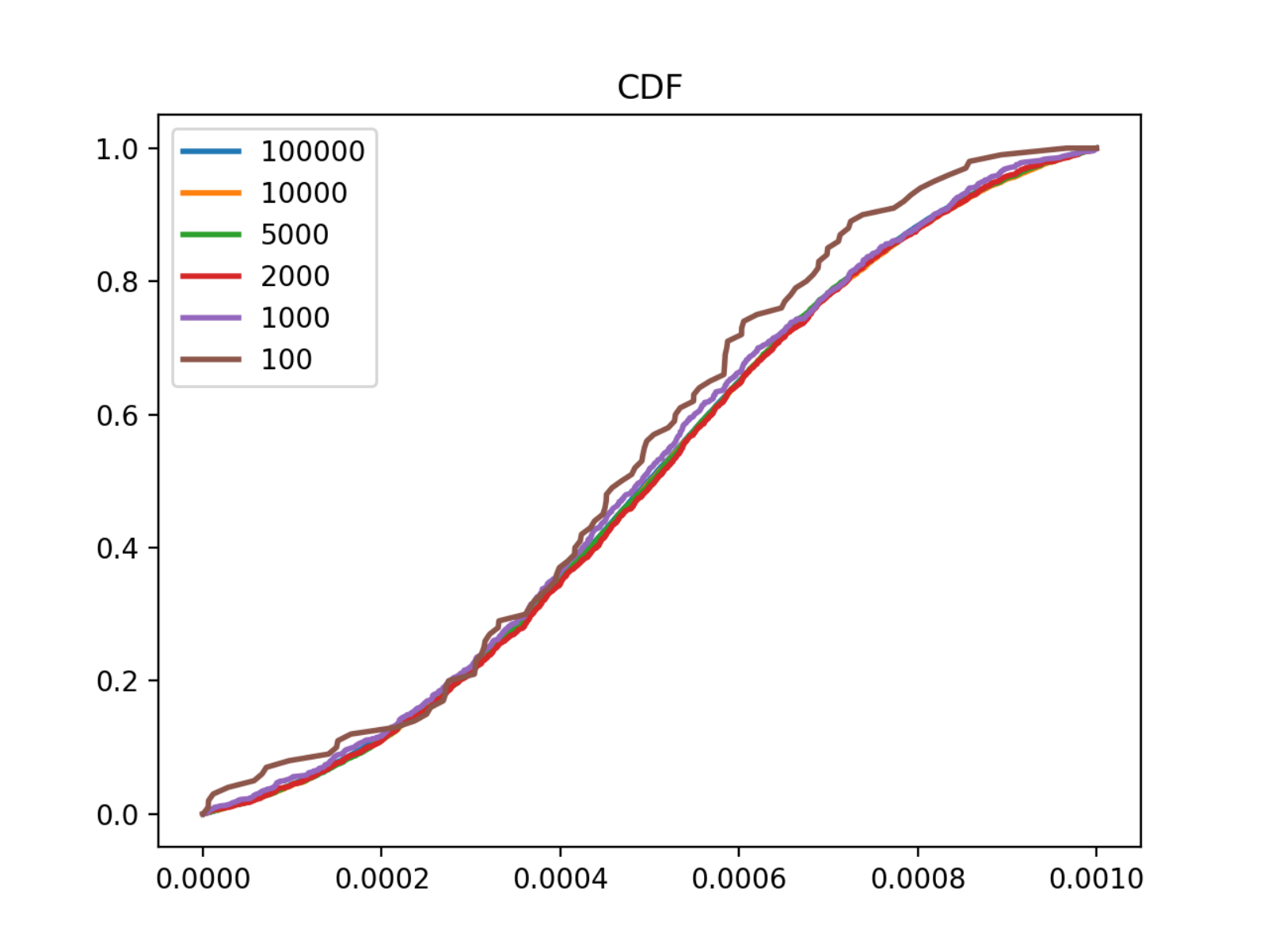}
  }
  \subfigure[$\epsilon = 0.1$, $v = 0.05$]{
    \label{fig:v005_eps01_icdf_rt}
    \includegraphics[width=0.30\textwidth]{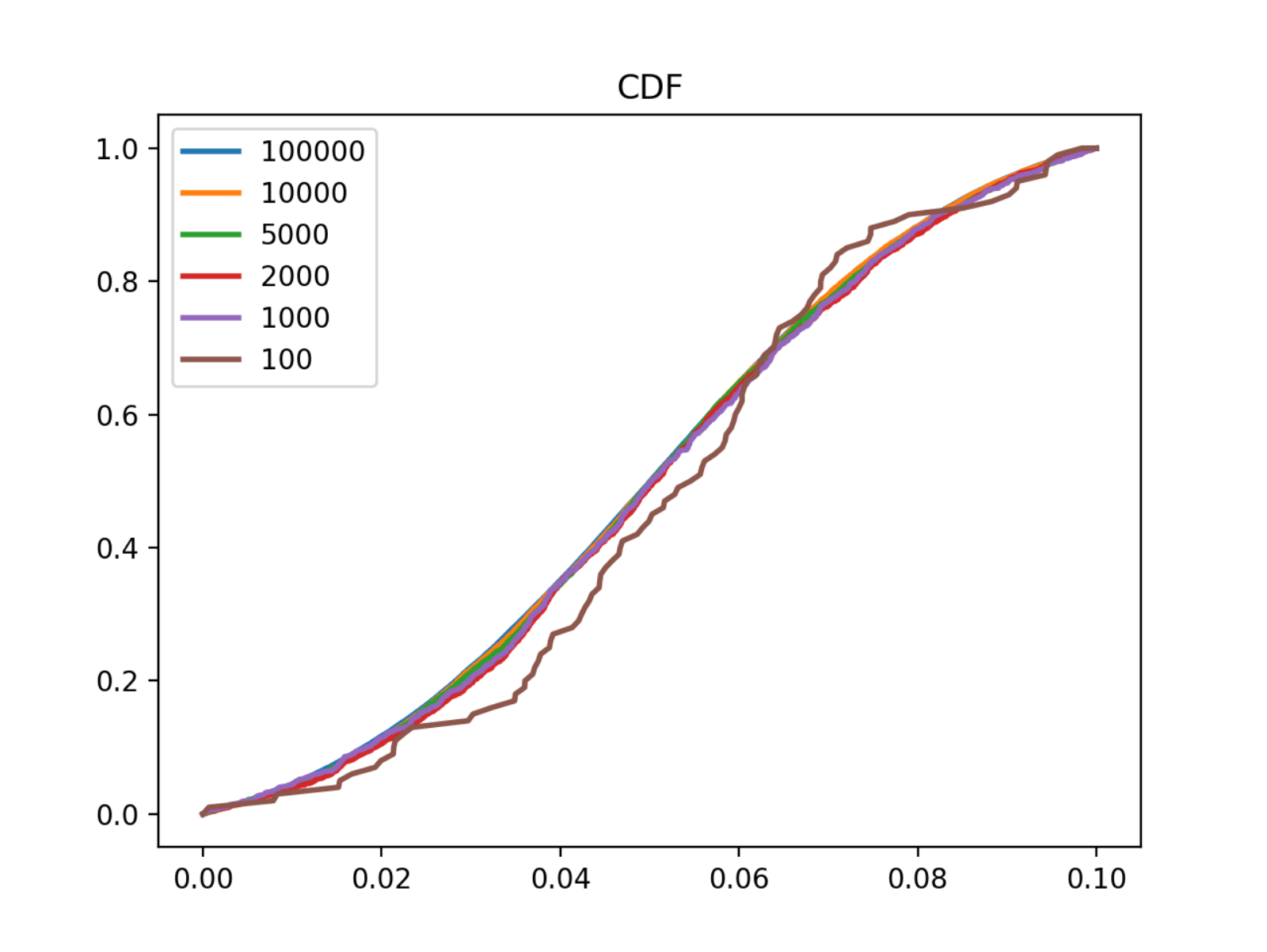}
  }
  ~
  \subfigure[$\epsilon = 0.01$, $v = 0.05$]{
    \label{fig:v005_eps001_icdf_rt}
    \includegraphics[width=0.30\textwidth]{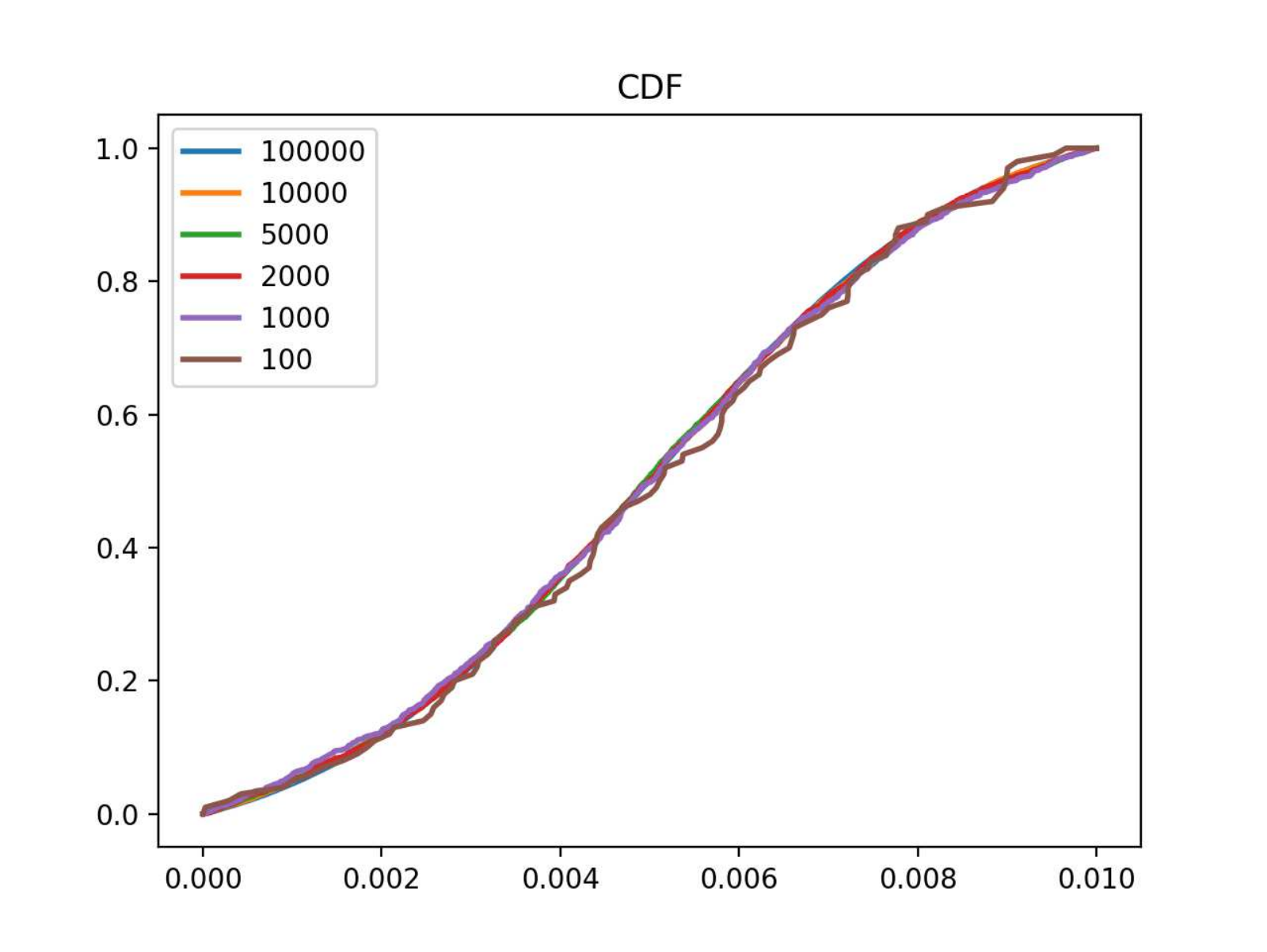}
  }
  ~
  \subfigure[$\epsilon = 0.001$, $v = 0.05$]{
    \label{fig:v005_eps0001_icdf_rt}
    \includegraphics[width=0.30\textwidth]{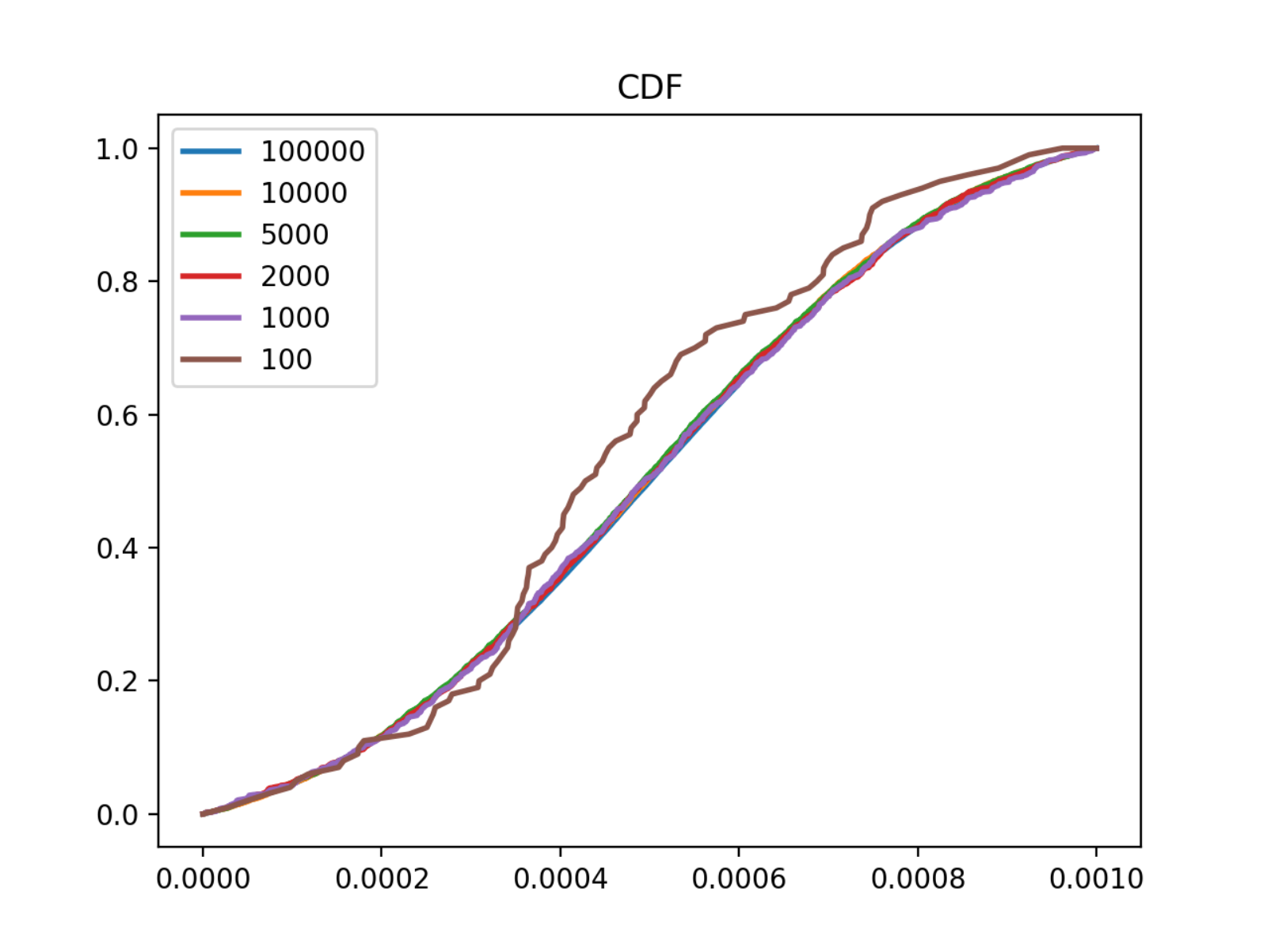}
  }
  \caption{Cumulative density functions of $\delta \sim G(0,\epsilon,\epsilon/2,\sqrt{\epsilon^{2}/12})$ for selected $\epsilon$ and $v$, calculated at different numbers of realizations $N$ of $\delta$. Shown are $N = 100,1000,2000,5000,10000,100000$. Sampling larger $N$ provides distributions that better approximate the theoretical \texttt{cdf}.}\label{fig:icdf_n}
  \end{minipage}
\end{figure}

The code in Figure \ref{code:cdf} details how $\bar{z_{\ast}}$, $\sigma_{z_{\ast}}$, $\bar{\delta}$, $\sigma_{\delta}$, \texttt{cdf}($z_{\ast}$), and \texttt{cdf}($\delta$) are calculated.
First, the results calculated using the Monte Carlo sampling in Figure \ref{code:sampling} are loaded from disk. Then, for each $N$ in $[100,1000,2000,5000,10000,100000]$, $N$ individual $z_{\ast}$ are randomly selected for calculations of $\bar{z_{\ast}}$, $\sigma_{z_{\ast}}$, $\bar{\delta}$, $\sigma_{\delta}$, and \texttt{cdf}($z_{\ast}$) (if \texttt{inputs} is False) or \texttt{cdf}($\delta$) (if \texttt{inputs} is True).  The use of \texttt{numpy.random.choice} ensures $N$ unique $z_{\ast}$ and $\delta$ are selected for the calculations of mean, standard deviation, and \texttt{cdf}.
The strategy used of writing intermediate results to disk is meant to reduce computational cost, as the expensive Monte Carlo sampling of the \texttt{solve} function will not need to be repeated for each new evaluation of $\bar{z_{\ast}}$, $\sigma_{z_{\ast}}$, $\bar{\delta}$, $\sigma_{\delta}$, \texttt{cdf}($z_{\ast}$) and \texttt{cdf}($\delta$).
The calculated mean locations $\bar{z_{\ast}}$ and corresponding standard deviations $\sigma_{z_{\ast}}$, due to uniform random perturbations $\delta \sim U(0,\epsilon)$ on the boundary condition, are presented in Table \ref{tab:multiN} for $v = 0.05$, $\epsilon = 0.1$, and several values of $N$. Results for $\delta \sim G(0,\epsilon,\epsilon/2,\sqrt{\epsilon^{2}/12})$, are also presented in the same table. 
Looking at Table \ref{tab:multiN}, one can infer that roughly $N=10000$ or greater should be used to calculate $\bar{z_{\ast}}$ and $\sigma{z_{\ast}}$.
The quality of the results due to increasing $N$ samplings of $\delta \sim U(0,\epsilon)$ can be visualized in plots of \texttt{cdf}($z_{\ast}$) and \texttt{cdf}($\delta$), in Figure \ref{fig:cdf} and Figure \ref{fig:icdf}, respectively.
Similarly, plots for $\delta \sim G(0,\epsilon,\epsilon/2,\sqrt{\epsilon^{2}/12})$ are shown in Figure \ref{fig:cdf_n} and Figure \ref{fig:icdf_n}.  The representative quality of the sampled distribution for large $N$ is seen by \texttt{cdf}($\delta$) closely approximating the theoretical \texttt{cdf}.
While Figure \ref{fig:cdf_n} shows that the basic shape of the \texttt{cdf}($z_{\ast}$) sharpens for smaller $v$ and doesn't change dramatically for any choice of $\epsilon$, it is apparent that the location of $\bar{z_{\ast}}$ does change with $\epsilon$.  A comparison can be more quantitatively made by examining the moments of $z_{\ast}$, as is presented in Table \ref{tab:multiepsilon}.
Results for $\bar{z_{\ast}}$, $\sigma_{z_{\ast}}$, $\bar{\delta}$, and $\sigma_{\delta}$ for $N=100000$ Monte Carlo realizations of $\delta$, for all combinations of $v$ and $\epsilon$ considered, for both $\delta \sim U(0,\epsilon)$ and $\delta \sim G(0,\epsilon,\epsilon/2,\sqrt{\epsilon^{2}/12})$.
It should be noted that the moments $\bar{\delta}$, and $\sigma_{\delta}$ for both the selected truncated Gaussian and uniform distributions of $\delta$ are constant with respect to $v$, and are directly scaled with $\epsilon$.
Also note that the moments of $\delta \sim U(0,\epsilon)$ and $\delta \sim G(0,\epsilon,\epsilon/2,\sqrt{\epsilon^{2}/12})$ are not a perfect match; the mean of $\delta$ is roughly identical, however $\sigma_{\delta}$ for each distribution is only similar.
The results show that $\bar{z_{\ast}}$ decreases with decreasing $\epsilon$, where the decrease in $\epsilon$ is more pronounced for larger $v$. The behavior of $\sigma_{z_{\ast}}$ mirrors that of $\bar{z_{\ast}}$, but on a smaller scale. Note that small perturbations in $\epsilon$ (or $v$) can lead to significant differences in the stochastic response in the output. Looking back at results from Xiu \cite{Xiu:2004}, we note that results for $\bar{z_{\ast}}$ with $\delta \sim U(0,\epsilon)$ are very consistent, while the results for $\sigma_{z_{\ast}}$ are less so.

\begin{figure}[htbp!]
  \centering
  \subfigure[$\epsilon = 0.1$, $v = 0.1$]{
    \label{fig:v01_eps01_ux}
    \includegraphics[width=0.3\textwidth]{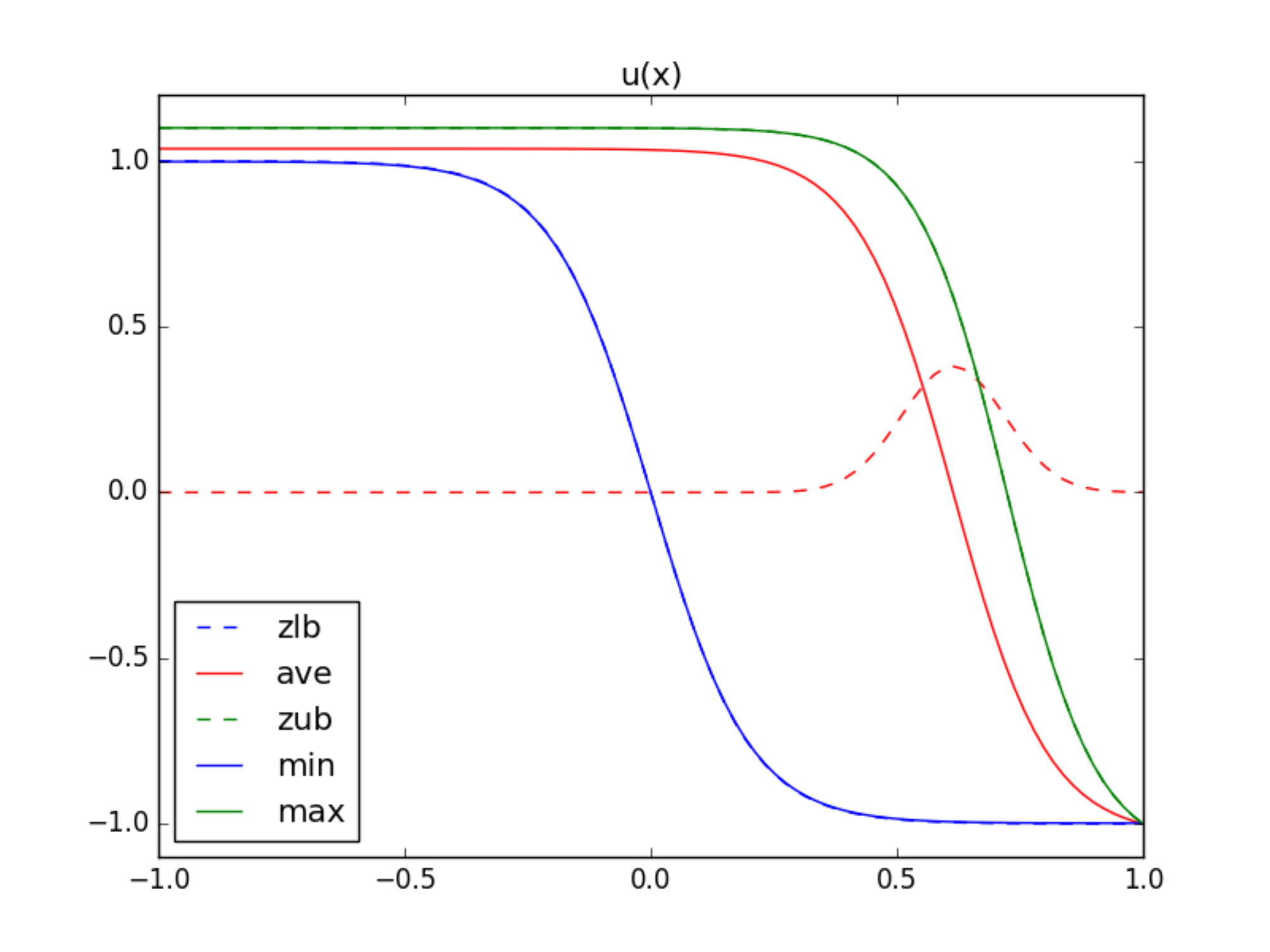}
  }
  ~
  \subfigure[$\epsilon = 0.01$, $v = 0.1$]{
    \label{fig:v01_eps001_ux}
    \includegraphics[width=0.3\textwidth]{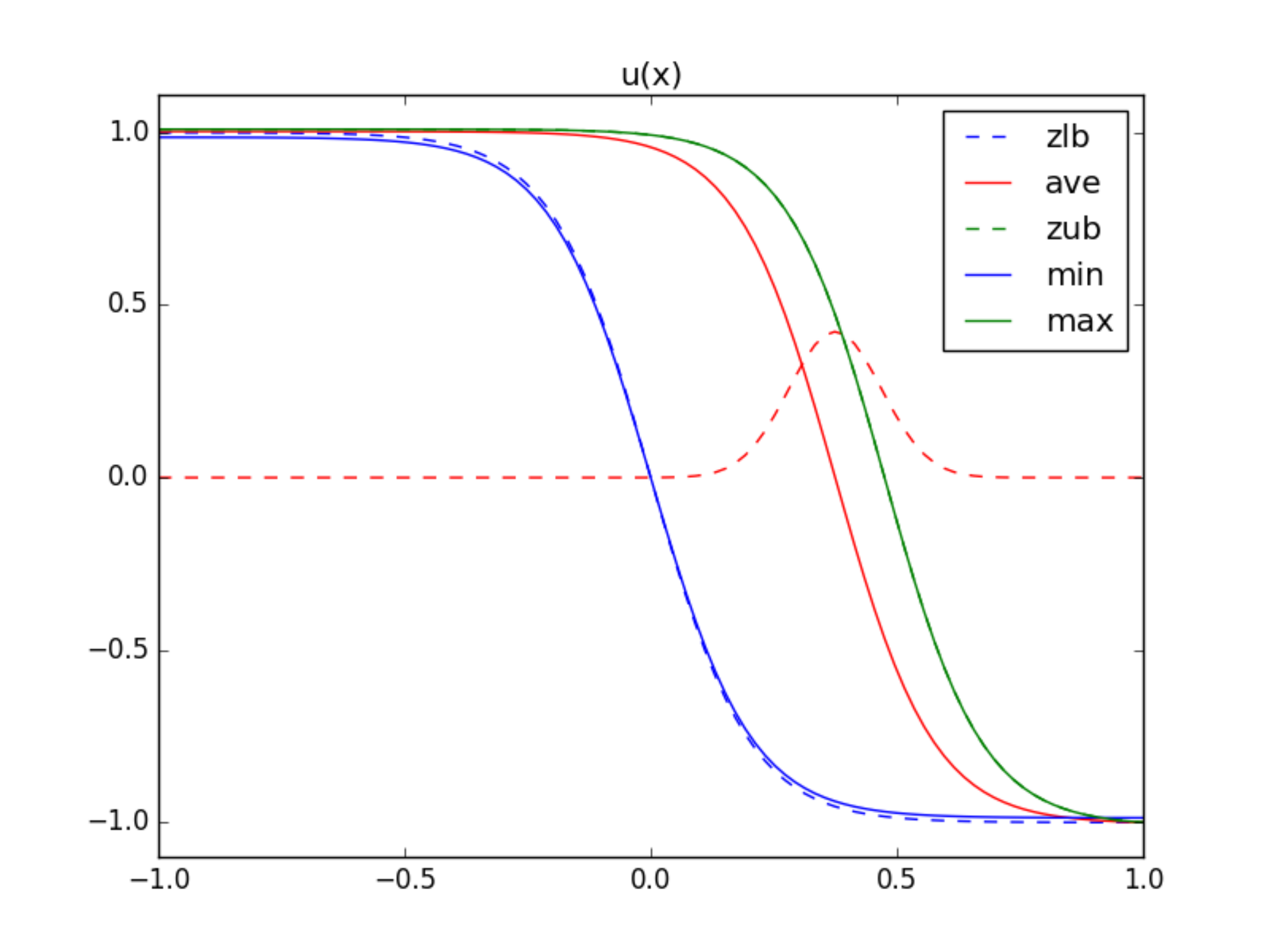}
  }
  ~
  \subfigure[$\epsilon = 0.001$, $v = 0.1$]{
    \label{fig:v01_eps0001_ux}
    \includegraphics[width=0.3\textwidth]{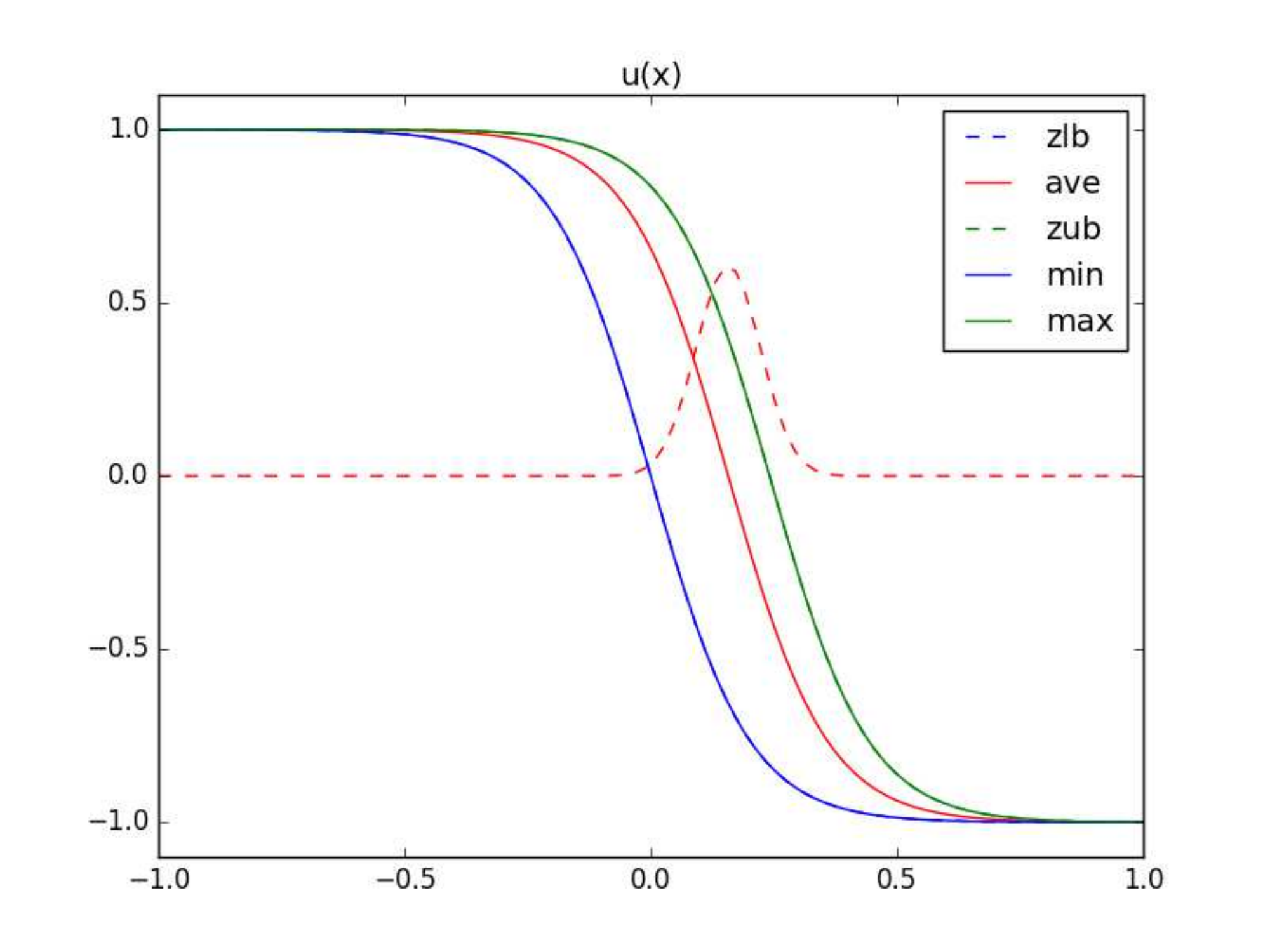}
  }
  \subfigure[$\epsilon = 0.1$, $v = 0.05$]{
    \label{fig:v005_eps01_ux}
    \includegraphics[width=0.3\textwidth]{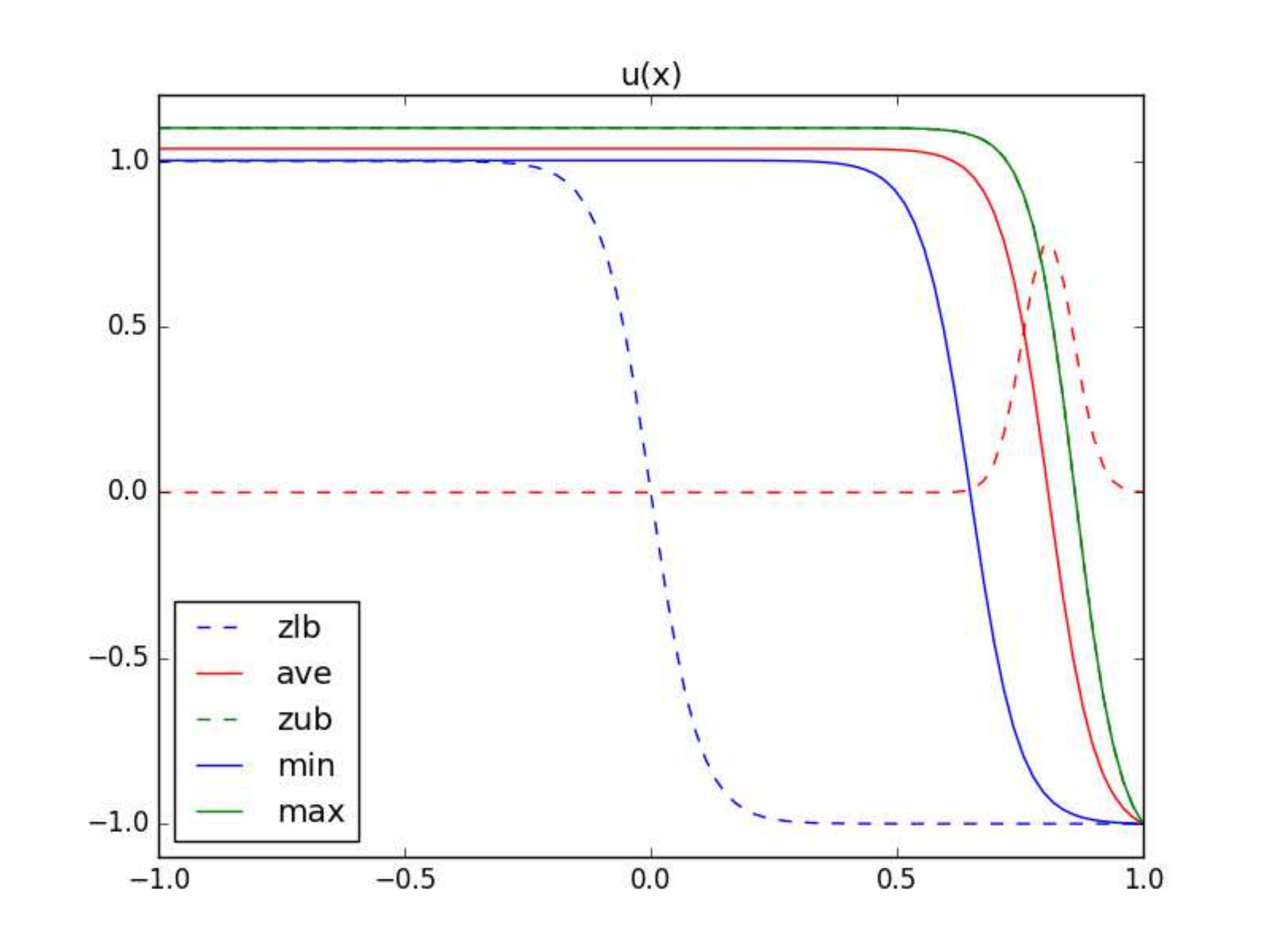}
  }
  ~
  \subfigure[$\epsilon = 0.01$, $v = 0.05$]{
    \label{fig:v005_eps001_ux}
    \includegraphics[width=0.3\textwidth]{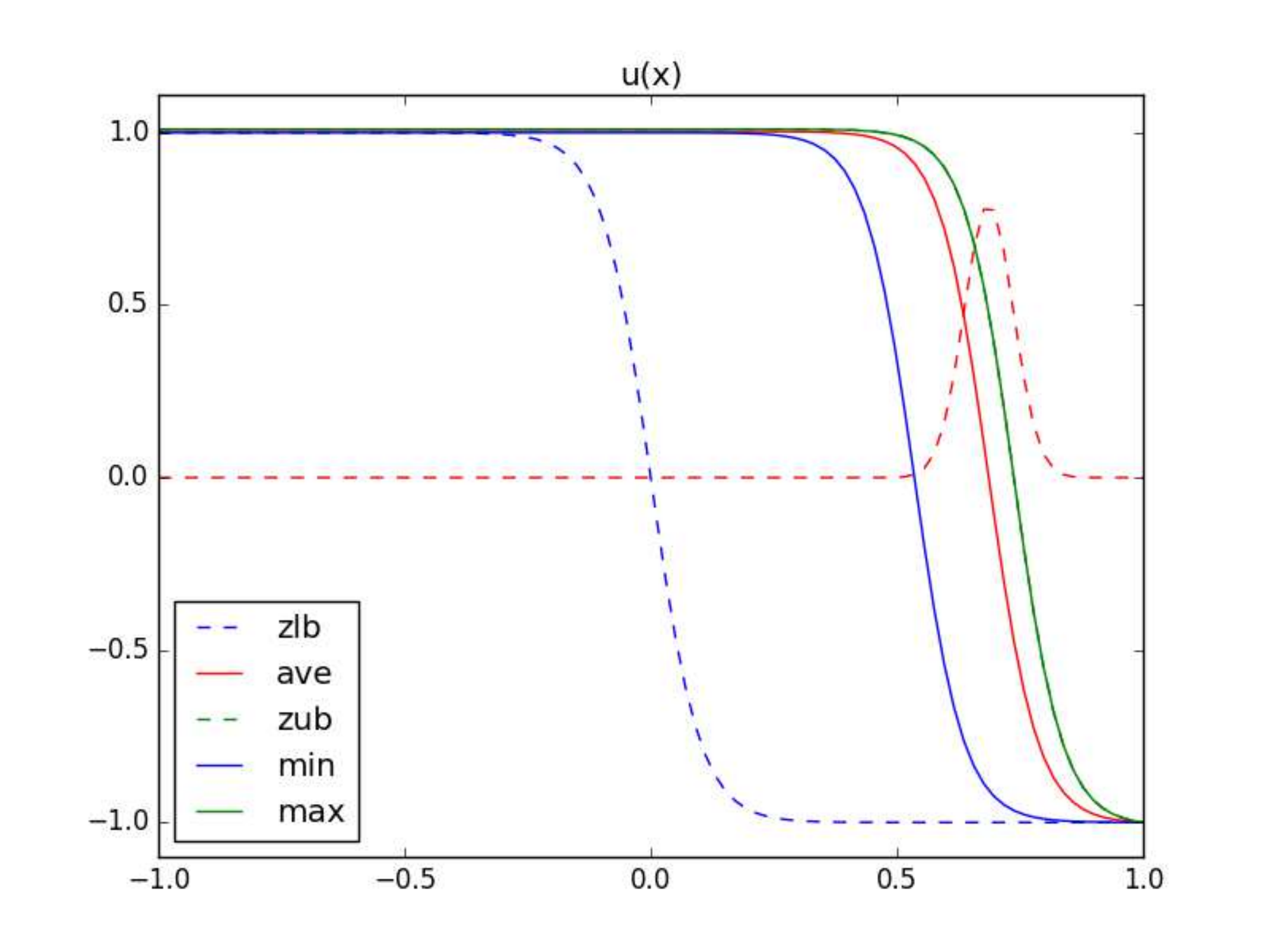}
  }
  ~
  \subfigure[$\epsilon = 0.001$, $v = 0.05$]{
    \label{fig:v005_eps0001_ux}
    \includegraphics[width=0.3\textwidth]{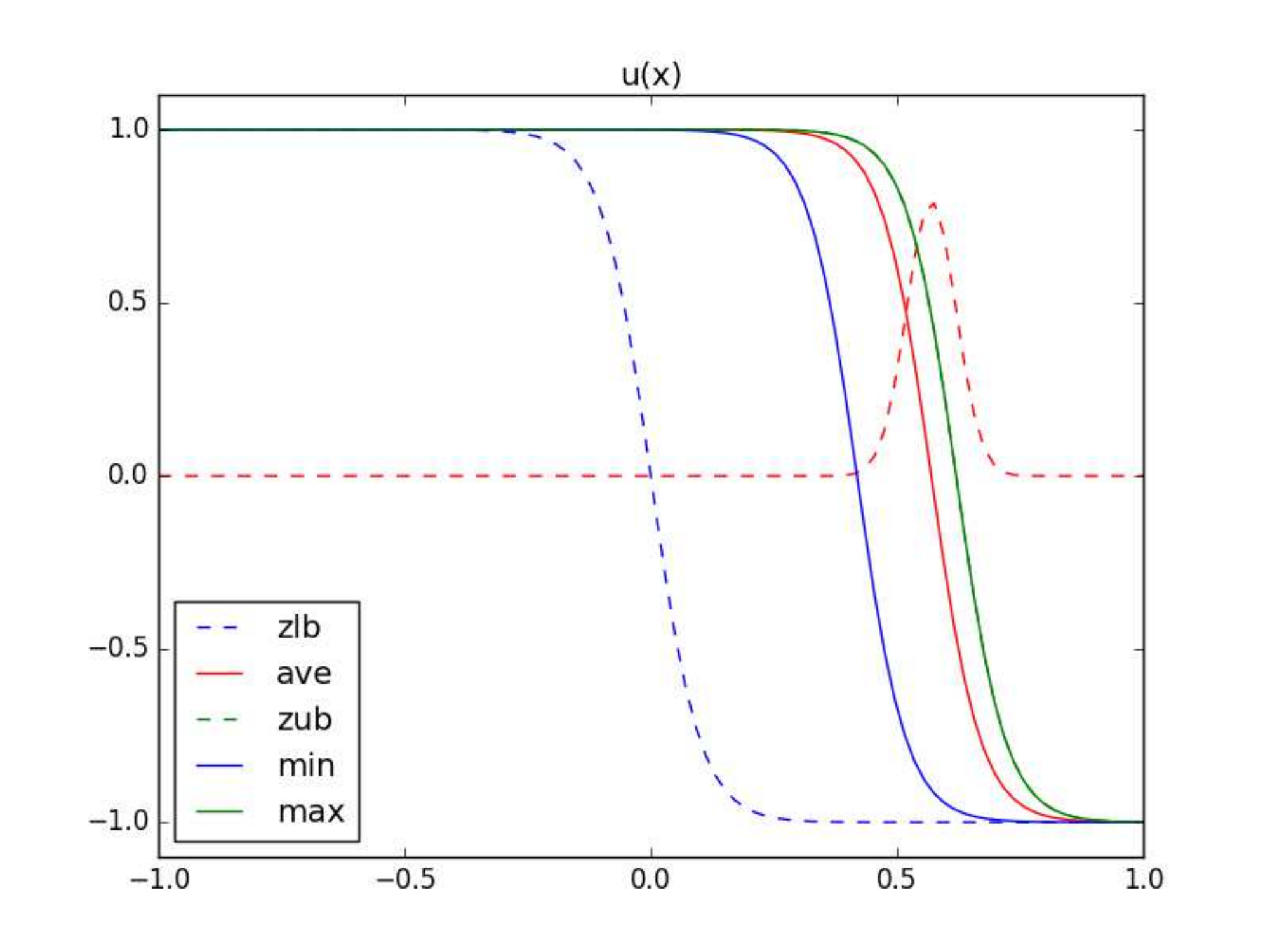}
  }
  \caption{Stochastic solutions for $\delta \sim U(0,\epsilon)$ and $v$. The upper and lower boundaries, $zub$ and $zlb$, are the deterministic solutions corresponding to the limiting values of the random inputs for $\delta$ (i.e. $[0, \epsilon]$). Also shown are $max$ and $min$, calculated using $z$ which are at most three standard deviations away from $\bar{z_{\ast}}$ (i.e. $ave$).}\label{fig:averages}
\end{figure}

To get an idea of the dynamics of $u(z)$ in (\ref{eq:burgers_steady}), we find stochastic solutions of Burgers' equation with a left boundary wall at $\delta \sim U(0, \epsilon)$ for a given $v$ and $\epsilon$.  The code in Figure \ref{code:averages} additionally provides the dynamics at an estimate of where the tails of the distribution for $z_{\ast}$ might be located -- it finds the solutions of $u(z)$ at $\delta(v,$ \texttt{max}$(\bar{z_{\ast}} - 3*\sigma_{z_{\ast}}, 0))$ and $\delta(v,$ \texttt{min}$(\bar{z_{\ast}} + 3*\sigma_{z_{\ast}}, \epsilon))$. Results are shown in Figure \ref{fig:averages}.

% END COMMENT

% results/plots of Monte Carlo sampling futures of z

% results/plots of Monte Carlo sampling of \delta

%\COMMENT{ % monte carlo sampling of delta, calculating \bar{z_{\ast}}
\begin{figure}[htbp!]
  \centering
  \begin{minipage}{1.00\textwidth}
% \begin{\outputtextsize}
\begin{verbatim}
### mc_supersensitive.py :: MC sampling of \delta, calculating \bar{z_{\ast}} ###
N = 100000; uniform = True
tol = 1e-09  # fit tolerance
M = int(round(1.1*N)) # number of samples N, with a padding

import numpy as np; import pickle; import rtnorm
from exact_supersensitive import solve
dist = np.random.uniform if uniform else rtnorm._rtnorm

if __name__ == '__main__':

  for v in (0.05, 0.1):
    for eps in (1e-1, 1e-2, 1e-3):
      from pathos.pools import ProcessPool
      p = ProcessPool()

      deltas = dist(0, eps, size=M) # get realizations of delta
      results = p.uimap(solve, [v]*len(deltas), deltas) # get futures of z
      z = np.empty(M, dtype=[('z',float),('delta',float),('fit',float)])
      i = j = 0
      for d,x,y in results:
        if i is N:
          break
        if y > tol:
          j += 1
          print 'miss {} at {} with {} from {}'.format(j,i,y,d)
          continue
        z[i] = x[0],d,y
        i += 1
      print "buffer, miss:", (M-N, j)
      p.close()

      z = np.sort(z[:N], order='z') # sort by the N 'best fit'
      fname = 'burgers_MC_{}_deltas_v_{}_eps_{}.pkl'.format(N,v,eps)
      z.dump(fname)
      with open(fname, 'ab') as f:
        pickle.dump((v,eps), f)
      p.join()
      zeros = lambda x: len(x) - np.count_nonzero(x)
      print "zeros (delta, z):", (zeros(z['delta']), zeros(z['z']))
      p.restart()
  p.close()
  p.join()

\end{verbatim}
% \end{\outputtextsize}
  \end{minipage}
  \caption{Monte Carlo sampling of $\delta$, where $z_{\ast}$ is calculated with \texttt{solve} in parallel using an iterative map \texttt{uimap} from \texttt{pathos}. After obtaining $N$ results within the selected tolerance $tol$, tuples of ($z_{\ast}$, $\delta$, \texttt{fit}) are written to disk. Failures to \texttt{solve} within $tol$ are dropped.}\label{code:sampling}
\end{figure}
%}%\END COMMENT

%\COMMENT{ % provide sampling from a truncated Gaussian distribution
\begin{figure}[htbp!]
  \centering
  \begin{minipage}{1.00\textwidth}
% \begin{\outputtextsize}
\begin{verbatim}
### rtnorm.py :: random sampling from a truncated Gaussian distribution ###
from scipy.stats import truncnorm

def rtnorm(a, b, loc=0., scale=1., size=None):
  """random sampling from a truncated Gaussian distribution"""
  mu, sigma = float(loc), float(scale)
  a, b = float(a), float(b)
  if not mu == 0. or not sigma == 1.:
    a, b = (a-mu) / sigma, (b-mu) / sigma
    r = rtstdnorm(a, b, size=size)
    return r * sigma + mu
  return rtstdnorm(a, b, size=size)

def _rtnorm(a, b, size=None):
  "rtnorm distribution with first two moments of a uniform distribution"
  return rtnorm(a, b, .5*(a+b), ((b-a)**2/12.)**.5, size)

def rtstdnorm(a, b, size=None):
  "rtnorm distribution with mean of zero and variance of one"
  if size is None:
    return truncnorm(a,b).rvs(1)[0]
  return truncnorm(a,b).rvs(size)

\end{verbatim}
% \end{\outputtextsize}
  \end{minipage}
  \caption{Functions to provide random sampling from a truncated Gaussian distribution in the interval $[a, b]$. \texttt{\_rtnorm} is a special case of \texttt{rtnorm}, providing sampling from a truncated normal distribution with \texttt{loc} and \texttt{scale} equal to the first two moments of a uniform normal distribution defined in the interval $[a, b]$.}\label{code:rtnorm}
\end{figure}
%}%\END COMMENT

%\COMMENT{ % compute mean, std, PDF|CDF from Monte Carlo samples
\begin{figure}[htbp!]
  \centering
  \begin{minipage}{1.00\textwidth}
% \begin{\outputtextsize}
\begin{verbatim}
### cdf_supersensitive.py :: compute mean, std, and cdf ###
v, eps = 0.1, 0.1
N = [100,1000,2000,5000,10000,100000]
inputs = False

import numpy as np; import pickle
from mystic.math.measures import mean, std
import matplotlib.pyplot as plt

def cdf(z, lo=None, hi=None):
  z = np.array(z)
  lo, hi = z.min() if lo is None else lo, z.max() if hi is None else hi
  z = np.sort(z[np.logical_and(z >= lo, z <= hi)])
  y = (1. + np.arange(len(z)))/len(z)
  if lo < z.min():
    z,y = np.insert(z, 0, lo), np.insert(y, 0, 0)  
  if hi > z.max():
    return np.append(z, hi), np.append(y, 1)
  return z,y

if __name__ == '__main__':

  fname = 'burgers_MC_{}_deltas_v_{}_eps_{}.pkl'.format(N[-1],v,eps)
  with open(fname, 'rb') as f:
    z = pickle.load(f)
  z,d = z['z'],z['delta']

  if not inputs:
    plt.xlim(0,1)
  for M in reversed(N): # select M samples
    z = np.random.choice(z, M, replace=False)
    d = np.random.choice(d, M, replace=False)
    print "v, eps, N:", (v, eps, M)
    print "mean, std:", mean(z), std(z) # mean and std of futures
    print "'' for delta:", mean(d), std(d)
    x,y = cdf(d if inputs else z, lo=0); title='CDF' # calculate the CDF
    plt.plot(x,y, linewidth=2)
  plt.title(title)
  plt.show()

\end{verbatim}
% \end{\outputtextsize}
  \end{minipage}
  \caption{Compute \texttt{mean}, \texttt{std}, and \texttt{cdf} for Monte Carlo futures of $z$ (if \texttt{inputs = False}) or \texttt{cdf} of $\delta$ (if \texttt{inputs = True}). The use of \texttt{random.choice} ensures we have random sampling with no duplication. Results for the moments of $z_{\ast}$ are summarized in Tables \ref{tab:multiN} and \ref{tab:multiepsilon}, while plots of the \texttt{cdf} can be found in Figures \ref{fig:cdf}, \ref{fig:icdf}, \ref{fig:cdf_n}, and \ref{fig:icdf_n}.}\label{code:cdf}
\end{figure}

%\END COMMENT

%TABLE: Mean locations $\bar{z_{\ast}}$ of the transition layer and the corresponding standard deviations $\sigma_{z_{\ast}}$ subject to perturbation $\delta \sim U(0,\epsilon)$ of the boundary condition.
%       Contents: z and \sigma_{z_{\ast}} at \epsilon = 0.1 and v = 0.05 for several N (Table V in Xiu)
%       Also mean(delta),std(delta)? Plot PDF? (Output of pdfcdf_supersensitive.py)
\begin{table}[htbp!]
  \centering
  \begin{tabular}{c c c c c c c c}
    \hline
    $\delta$ & $N$ & $100$ & $1000$ & $2000$ & $5000$ & $10000$ & $100000$ \\
    \hline
   %$\bar{z_{\ast}}$ & 0.812405962248 & 0.809010038689 & 0.807908658202 & 0.808533446971 & 0.807786790916 & 0.807264192348 \\
   %$\sigma_{z_{\ast}}$ & 0.0473087338875 & 0.0482140926789 & 0.051597770357 & 0.051291624943 & 0.0521025021647 & 0.052791379663 \\
   %$U$ & $\bar{z_{\ast}}$ & 0.81240596 & 0.80901004 & 0.80790866 & 0.80853345 & 0.80778679 & 0.80726419 \\
   %$ $ & $\sigma_{z_{\ast}}$ & 0.04730873 & 0.04821409 & 0.05159777 & 0.05129162 & 0.05210250 & 0.05279138 \\
    $U$ & $\bar{z_{\ast}}$ & 0.81258074 &  0.80846617 & 0.80955572 & 0.80819600 & 0.80813561 & 0.80743291 \\
    $ $ & $\sigma_{z_{\ast}}$ & 0.04115351 & 0.04967964 & 0.04936012 & 0.05234097 & 0.05225653 & 0.05250859 \\
    \hline
   %$G$ & $\bar{z_{\ast}}$ & 0.81246440 & 0.81519252 & 0.81354562 & 0.81396660 & 0.81418272 & 0.81400354 \\
   %$ $ & $\sigma_{z_{\ast}}$ & 0.03928706 & 0.03675762 & 0.03906461 & 0.03882715 & 0.03931032 & 0.03893414 \\
    $G$ & $\bar{z_{\ast}}$ & 0.81899998 & 0.81578016 & 0.81479423 & 0.81442699 & 0.81452547 & 0.81414045 \\
    $ $ & $\sigma_{z_{\ast}}$ & 0.02995636 & 0.03700579 & 0.03694921 & 0.03770155 & 0.03818269 & 0.03847056 \\
    \hline
  \end{tabular}
  \caption{The mean locations $\bar{z_{\ast}}$ of the transition layer and the corresponding standard deviations $\sigma_{z_{\ast}}$ subject to uniform random perturbation $\delta \sim U(0,\epsilon)$ (or truncated Gaussian random perturbation $\delta \sim G(0,\epsilon,\epsilon/2,\sqrt{\epsilon^{2}/12})$) on the boundary condition. In every case, Monte Carlo futures of $z$ were calculated with $v=0.05$ and $\epsilon=0.1$. Results approach asymptotic for large $N$. In Table \ref{tab:multiepsilon}, we examine $N=100000$ for various $v$ and $\epsilon$.}\label{tab:multiN}
\end{table}

% END COMMENT

% results/plots of Monte Carlo sampling futures of z (for rtnorm)

% results/plots of Monte Carlo sampling of \delta (for rtnorm)

%TABLE: Mean locations $\bar{z_{\ast}}$ of the transition layer and the corresponding standard deviations $\sigma_{z_{\ast}}$ subject to perturbation $\delta \sim U(0,\epsilon)$ of the boundary condition.
%       Contents: z and \sigma_{z_{\ast}} for \epsilon and v at N=100000 (Table IV in Xiu)
%       Also see plots of MC PoF (v01.eps and v005.eps)
%       (Output of pof_supersensitive.py)
\begin{table}[htbp!]
  \centering
  \begin{tabular}{c c c c c c c c c}
    \hline
    \multicolumn{1}{c}{$ $} & \multicolumn{1}{c}{$ $} & \multicolumn{3}{c}{$v=0.05$} & \multicolumn{1}{c}{$ $} & \multicolumn{3}{c}{$v=0.1$} \\
    \cline{3-5} \cline{7-9}
    $\delta$ & $\epsilon$ & $0.1$ & $0.01$ & $0.001$ & $ $ & $0.1$ & $0.01$ & $0.001$ \\ 
    \hline
   %$\bar{z_{\ast}}$ & 0.807264192348 & 0.686467241189 & 0.569867148211 & $ $ & 0.614071648556 & 0.374636048043 & 0.157444166571 \\
   %$\sigma_{z_{\ast}}$ & 0.052791379663 & 0.0503092101224 & 0.0502886375178 & $ $ & 0.105005892629 & 0.0943062706661 & 0.0657727983528 \\
   %$U$ & $\bar{z_{\ast}}$ & 0.80726419 & 0.68646724 & 0.56986715 & $ $ & 0.61407165 & 0.37463605 & 0.15744417 \\
   %$ $ & $\sigma_{z_{\ast}}$ & 0.05279138 & 0.05030921 & 0.05028864 & $ $ & 0.10500589 & 0.09430627 & 0.06577280 \\
   %$ $ & $\bar{\delta}$ & 0.04997934 & 0.00500006 & 0.00049834 & $ $ & 0.04977898 & 0.00500095 & 0.00049940 \\
   %$ $ & $\sigma_{\delta}$ & 0.02889793 & 0.00288701 & 0.00028875 & $ $ & 0.02884908 & 0.00288415 & 0.00028903 \\
   %$\bar{\delta}$ & 0.0499793422551 & 0.00500005599738 & 0.000498341365321 & $ $ & 0.0497789769242 & 0.00500095159544 & 0.00049940186669 \\
    $U$ & $\bar{z_{\ast}}$ & 0.80743291 & 0.68614016 & 0.57015331 & $ $ & 0.61442027 & 0.37461811 & 0.15771333 \\
    $ $ & $\sigma_{z_{\ast}}$ & 0.05250859 & 0.05086699 & 0.05024926 & $ $ & 0.10501165 & 0.09423549 & 0.06572515 \\
    $ $ & $\bar{\delta}$ & 0.05001091 & 0.00498799 & 0.00050024 & $ $ & 0.04994279 & 0.00499975 & 0.00050054 \\
    $ $ & $\sigma_{\delta}$ & 0.02881884 & 0.00288971 & 0.00028850 & $ $ & 0.02891048 & 0.00288486 & 0.00028875 \\
   %$\bar{\delta}$ & 0.0499793422551 & 0.00500005599738 & 0.000498341365321 & $ $ & 0.0497789769242 & 0.00500095159544 & 0.00049940186669 \\
   %$\sigma_{\delta}$ & 0.0288979275828 & 0.00288701400155 & 0.000288754351005 & $ $ & 0.0288490846615 & 0.00288415044489 & 0.000289038659948 \\
    \hline
   %$G$ & $\bar{z_{\ast}}$ & 0.81400354 & 0.69283637 & 0.57672381 & $ $ & 0.62847655 & 0.38668929 & 0.16336534 \\
   %$ $ & $\sigma_{z_{\ast}}$ & 0.03893414 & 0.03677701 & 0.03628689 & $ $ & 0.07618020 & 0.06970244 & 0.05202887 \\
   %$ $ & $\bar{\delta}$ & 0.05002038 & 0.00498961 & 0.00050032 & $ $ & 0.05007693 & 0.00499849 & 0.00049933 \\
   %$ $ & $\sigma_{\delta}$ & 0.02347176 & 0.00234851 & 0.00023536 & $ $ & 0.02347501 & 0.00234817 & 0.00023476 \\
    $G$ & $\bar{z_{\ast}}$ & 0.81414045 & 0.69277367 & 0.57648146 & $ $ & 0.62783214 & 0.38670812 & 0.16337040 \\
    $ $ & $\sigma_{z_{\ast}}$ & 0.03847056 & 0.03692673 & 0.03692400 & $ $ & 0.07657660 & 0.06974101 & 0.05219059 \\
    $ $ & $\bar{\delta}$ & 0.05005416 & 0.00498816 & 0.00049950 & $ $ & 0.04989020 & 0.00499956 & 0.00049957 \\
    $ $ & $\sigma_{\delta}$ & 0.02347194 & 0.00235580 & 0.00023565 & $ $ & 0.02352199 & 0.00234733 & 0.00023499 \\
    \hline
  \end{tabular}
  \caption{The mean locations $\bar{z_{\ast}}$ of the transition layer and the corresponding standard deviations $\sigma_{z_{\ast}}$ subject to uniform random perturbation $\delta \sim U(0,\epsilon)$ (or truncated Gaussian random perturbation $\delta \sim G(0,\epsilon,\epsilon/2,\sqrt{\epsilon^{2}/12})$) on the boundary condition. Similarly, $\bar{\delta}$ and $\sigma_{\delta}$ are calculated for $\delta$. $N=100000$ Monte Carlo realizations of $\delta$ were used in each case. The sensitivity of $z_{\ast}$ to $v$ and $\epsilon$ can also be seen in Figures \ref{fig:cdf} and \ref{fig:cdf_n}.}\label{tab:multiepsilon}
\end{table}

% END COMMENT

\COMMENT{ % Monte Carlo sampling of deltas, calculate futures of z
%TABLE: Number of "misses" (trials where solve does not converge to tol). Also included are the average convergence tol for the missed trials and the time to sample N points..
\begin{verbatim}
miss 1 at 56 with 0.0227541158362 from 0.0579323860128
miss 2 at 420 with 1.08271734267e-06 from 0.082008021016
miss 3 at 476 with 0.0362331948064 from 0.0362406524749
miss 4 at 522 with 0.0296308387651 from 0.0401611944574
miss 5 at 625 with 0.00747596653335 from 0.082584101347
miss 6 at 637 with 0.00616529352484 from 0.00603517309019
...
miss 773 at 108684 with 2.67659212491e-09 from 0.000938867664105
miss 774 at 108744 with 0.00527182991681 from 0.000328841659553
miss 775 at 108978 with 0.00297695675124 from 0.000796011693718
miss 776 at 108984 with 4.7394108347e-08 from 0.000565791469853
miss 777 at 109006 with 4.53690456314e-05 from 0.000592830896188
miss 778 at 109120 with 0.00322534850217 from 4.27416905986e-05
buffer, miss: (10000, 778)
started: 2016-09-08 11:26:31.157277
ended: 2016-09-08 14:03:04.713190
took: 2:36:33.555913
wrote results to burgers_MC_100000_deltas_v_0.1_eps_0.001.pkl
zeros (delta, z): (0, 0)
\end{verbatim}
}%\END COMMENT
\begin{table}[htbp!]
  \centering
  \begin{tabular}{c c c c c c c c c}
    \hline
    \multicolumn{1}{c}{$ $} & \multicolumn{1}{c}{$ $} & \multicolumn{3}{c}{$v=0.05$} & \multicolumn{1}{c}{$ $} & \multicolumn{3}{c}{$v=0.1$} \\
    \cline{3-5} \cline{7-9}
    $\delta$ & $\epsilon$ & $0.1$ & $0.01$ & $0.001$ & $ $ & $0.1$ & $0.01$ & $0.001$ \\ 
    \hline
    $U$ & $misses$ & 3 & 3 & 4 & $ $ & 0 & 0 & 0 \\
    $ $ & $tol$ & 6.7e-6  & 2.6e-8 & 1.4e-8 & $ $ & - & - & - \\
    $ $ & $time$ & 2:08:44 & 2:09:42 & 2:03:50 & $ $ & 2:11:18 & 2:00:09 & 2:04:54 \\
    \hline
   $G$ & $misses$ & 2 & 2 & 5 & $ $ & 0 & 0 & 1 \\
   $ $ & $tol$ & 2e-5 & 3.2e-8 & 7.6e-9 & $ $ & - & - & 5e-9 \\
   $ $ & $time$ & 2:03:02 & 2:02:04 & 1:55:16 & $ $ & 2:00:57 & 1:59:25 & 7:11:50 \\
    \hline
  \end{tabular}
  \caption{Details on the number of times \texttt{solve} fails to converge within \texttt{tol = 1e-9} (i.e. ``$miss$'') for \texttt{M = 110000} potential samplings of $\delta$. The number of $misses$, average $tol$ for a failed \texttt{solve}, and $time$ for collecting $N = 100000$ samples is provided.}
 \label{tab:misses}
\end{table}

%\COMMENT{ % calculate min/average/max exact solution of Burger's equation
\begin{figure}[htbp!]
  \centering
  \begin{minipage}{1.00\textwidth}
% \begin{\outputtextsize}
\begin{verbatim}
### u_supersensitive.py :: min/ave/max for solutions of Burgers' equation ###
import matplotlib.pyplot as plt
import numpy as np; import scipy.stats as ss
from mystic.solvers import fmin_powell
from exact_supersensitive import solve

def _delta(v, z):  # calculate delta from v and z_{\ast}
    # zdiff is abs(solved z_{\ast} - given z_{\ast}) for v and delta (i.e. x[0])
    zdiff = lambda x,v,z: abs(solve(v, x[0])[1][0] - z)
    return fmin_powell(zdiff, [0.0], args=(v, z), disp=0, full_output=1)[0]

def _u(v, delta):  # generate u(x) for given v,delta
    def _zA(v, delta):  # calculate z_{\ast},A from v and delta
        z, A = solve(v, delta)[1]
        return z, abs(A)
    v = float(v)
    z,A = _zA(v, delta)
    # build a function u(x) for fixed v,A,z_{\ast}
    return lambda x: -A*np.tanh(0.5*(A/v)*(x - z))

v, eps = 0.1, 0.1
z, std = 0.614420266306, 0.105011650866
x = np.linspace(-1,1,100); n = 3

umin = _u(v, 0.0)  # calculating u at lower limit for delta
plt.plot(x, umin(x), label='zlb', color='b', linestyle='dashed')
lb, ub = -1.1, 1.1
uave = _u(v, _delta(v, z))  # calculating u at z_{\ast}
plt.plot(x, uave(x), label='ave', color='r')
ub = max(ub, max(uave(x))+.1)
umax = _u(v, eps)  # calculating u at upper limit for delta
plt.plot(x, umax(x), label='zub', color='g', linestyle='dashed')
ub = max(ub, max(umax(x))+.1)

# calculate lower and upper limits based on z_{\ast} +/- n*std
zlb, zub = max(z-n*std,0.0), min(z+n*std,eps)
uzlb = umin if zlb is 0.0 else _u(v, _delta(v, zlb))  # u at min: n*std
uzub = umax if zub is eps else _u(v, _delta(v, zub))  # u at max: n*std
plt.plot(x, .1*ss.norm(z,std).pdf(x), color='r', linestyle='dashed')
plt.plot(x, uzlb(x), label='min', color='b')
plt.plot(x, uzub(x), label='max', color='g')
plt.xlim(-1, 1); plt.ylim(lb, ub)
plt.legend(); plt.title('u(x)')
plt.show()

\end{verbatim}
% \end{\outputtextsize}
  \end{minipage}
  \caption{Calculation of $u(z)$ due to direct solution of Burger's equation at 
the lower and upper limit of $\delta$ and at $\delta(v, \bar{z_{\ast}})$.  We also solve for $u(z)$ at $\delta(v,$ \texttt{max}$(\bar{z_{\ast}} - 3*\sigma_{z_{\ast}}, 0))$ and $\delta(v,$ \texttt{min}$(\bar{z_{\ast}} + 3*\sigma_{z_{\ast}}, \epsilon))$. The values of \texttt{z} and \texttt{std} are copied from Table \ref{tab:multiepsilon}.}
  \label{code:averages}
\end{figure}
%}%\END COMMENT

% results/plots of min/average/max for exact solution

% \newpage

\section{Using Sampling to Estimate Probability of Success}
The probability that some condition is true can be determined in a very straightforward way by sampling.  By applying sampling on an indicator function that defines the criteria for ``success'', we can calculate the probability of success $P($\texttt{success}$)$. We define
\begin{equation}
\label{eqn:success} 
  \tt{success} \defeq
  \begin{cases}
    \it{z_{\ast}} > \frac{\it{x}\bar{\it{z_{\ast}}}}{\text{100}}, \text{with} \\
    \it{x} = \text{100} + \it{dx}, \it{dx} \in \lbrack \text{0,15} \rbrack
  \end{cases}
\end{equation}
(i.e. when the shock occurs at a distance $z_{\ast}$ that is greater than the product of average distance $\bar{z_{\ast}}$ and scaling factor $\frac{x}{100}$). The scaling factor is introduced to demonstrate the decrease in probability of success as the target distance for the shock event moves further away from the origin.

The code in Figure \ref{code:pof} details how the indicator function \texttt{success} is passed the samples generated by the code in Figure \ref{code:sampling}, and is used to produce $P($\texttt{success}$)$. Monte Carlo estimates for $P($\texttt{success}$)$ can be found in Figure \ref{fig:sampling} for both $\delta \sim U(0,\epsilon)$ and $\delta \sim G(0,\epsilon,\epsilon/2,\sqrt{\epsilon^{2}/12})$.  While the code in Figure \ref{code:pof} uses \texttt{v = 0.1} and \texttt{eps = 0.1}, simply editing those two variables will calculate $P($\texttt{success}$)$ for other values of $v$ and $\epsilon$.
As see in the latter figure, the selection of a different distribution for $\delta$ does have an effect on the results, with the probability of success for $\delta \sim U(0,\epsilon)$ found to be slightly greater than for $\delta \sim G(0,\epsilon,\epsilon/2,\sqrt{\epsilon^{2}/12})$.  Additional calculation time required, beyond what is listed in Table \ref{tab:misses}, is negligible.

%\begin{compactitem}
% \item Discuss code for Monte Carlo calculation of $P(z_{\ast} > \frac{x\bar{z_{\ast}}}{100})$.
% \item Discuss timings and compuational complexity for the above.
%\end{compactitem}

%\COMMENT{ % calculate P(z_{\ast} > \bar{z_{\ast}}) with Monte Carlo
\begin{figure}[htbp!]
  \centering
  \begin{minipage}{1.00\textwidth}
% \begin{\outputtextsize}
\begin{verbatim}
### pof_supersensitive.py :: calculate P(success) with Monte Carlo ###
v, eps = 0.1, 0.1
M = 100000  # number of samples
pcnt = range(0,16)

import numpy as np; import pickle
from mystic.math.measures import mean, std

fname = 'burgers_MC_{}_deltas_v_{}_eps_{}.pkl'.format(M,v,eps)
with open(fname, 'rb') as f:
  z = pickle.load(f)
z = z['z']

z_, s_z = mean(z), std(z)
print "v, eps, N:", (v, eps, M)
print "mean, std:", z_, s_z  # sampled mean and std of z_{\ast}

success = lambda z,zave: z > zave  # define "success" indicator function

for p in pcnt:
  dz = round(1 + 0.01*p, 2)
  prob = success(z, z_{\ast}dz).sum()/float(len(z))  # probability of success
  print "P(z > %r*bar{z}) = %r" % (dz,prob)

\end{verbatim}
% \end{\outputtextsize}
  \end{minipage}
  \caption{Monte Carlo estimate for probability of success $P(z_{\ast} > \frac{x\bar{z_{\ast}}}{100})$ at $\delta \sim U(0,\epsilon)$ (or $\delta \sim G(0,\epsilon,\epsilon/2,\sqrt{\epsilon^{2}/12})$) for $\epsilon = 0.1$ and $v = 0.1$. We define \texttt{success} as where $z_{\ast} > \frac{x\bar{z_{\ast}}}{100}$, with $x = 100 + dx$, $dx \in [0,15]$ as in (\ref{eqn:success}). See Figure \ref{fig:MCdist} for the resulting distribution when the sampling of $\delta$ and subsequent calculation of $P($\texttt{success}$)$ is repeated $100000$ times.}
  \label{code:pof}
\end{figure}
%}%\END COMMENT

% \newpage

\section{Calculating Bounds on Probability of Success}
\begin{figure}[htbp!]
  \centering
  \subfigure[uniform; $v = 0.1$]{
    \label{fig:v01}
    \includegraphics[width=0.45\textwidth]{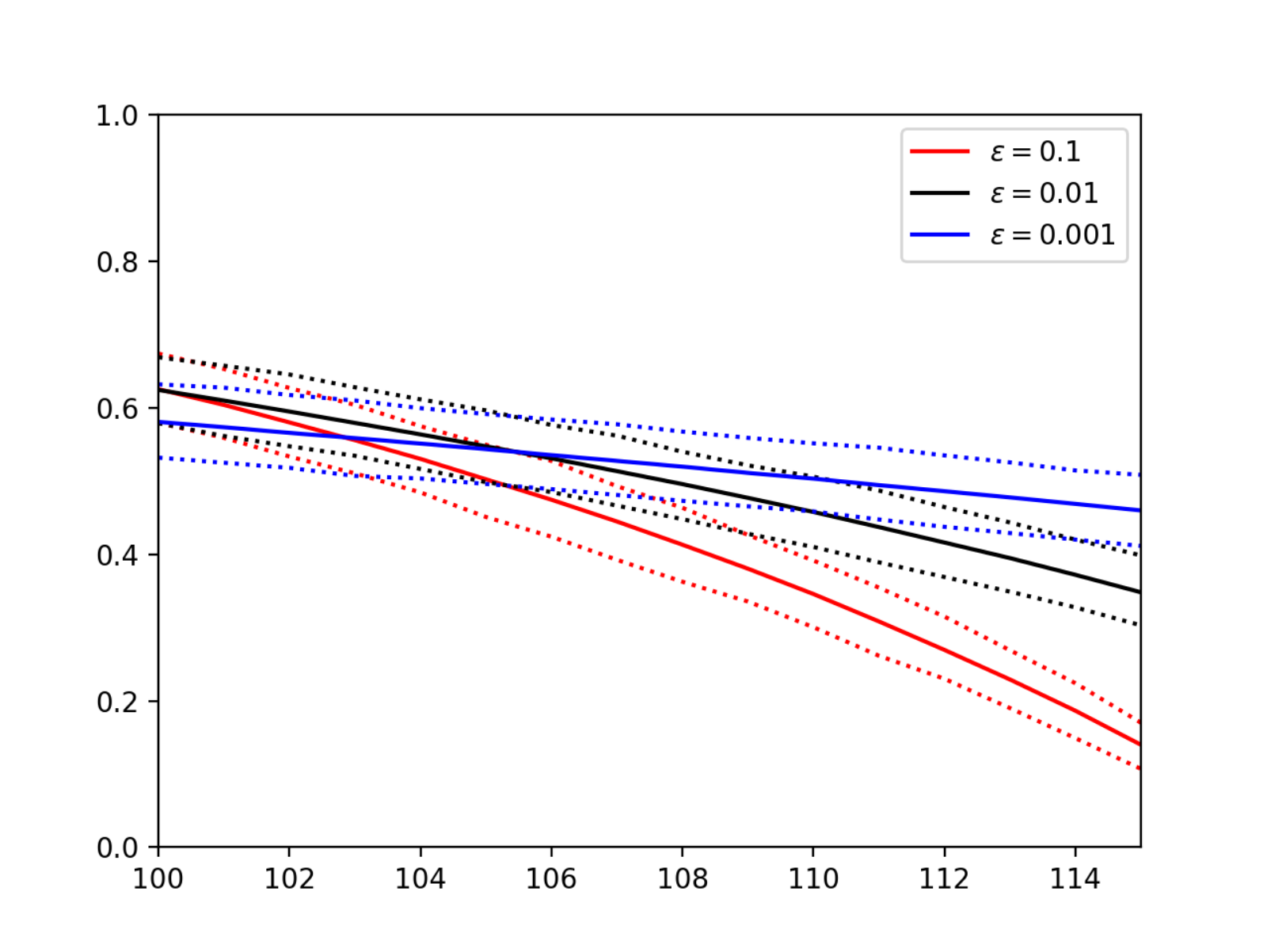}
  }
  ~
  \subfigure[rtnorm; $v = 0.1$]{
    \label{fig:v01r}
    \includegraphics[width=0.45\textwidth]{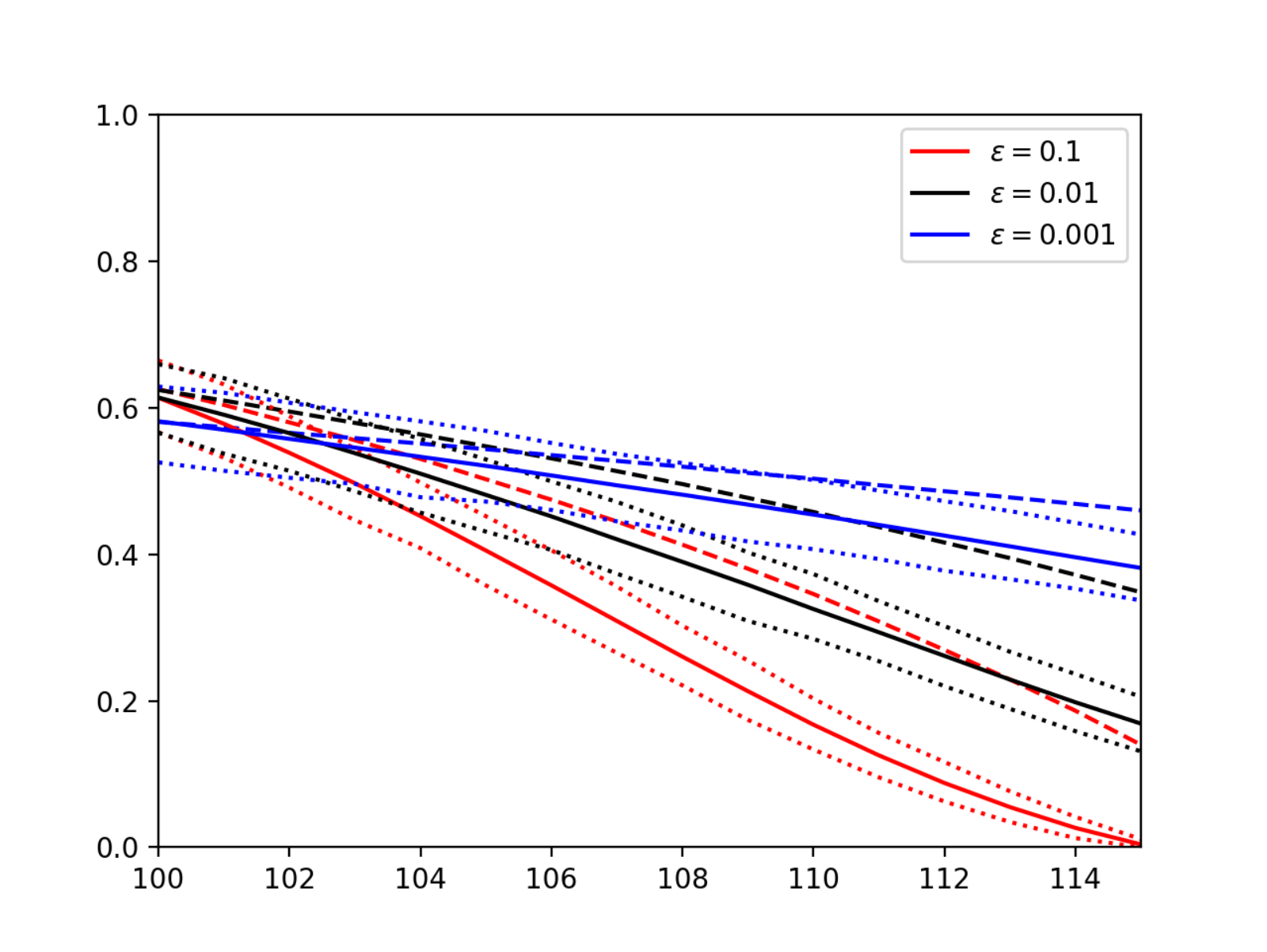}
  }
  \subfigure[uniform; $v = 0.05$]{
    \label{fig:v005}
    \includegraphics[width=0.45\textwidth]{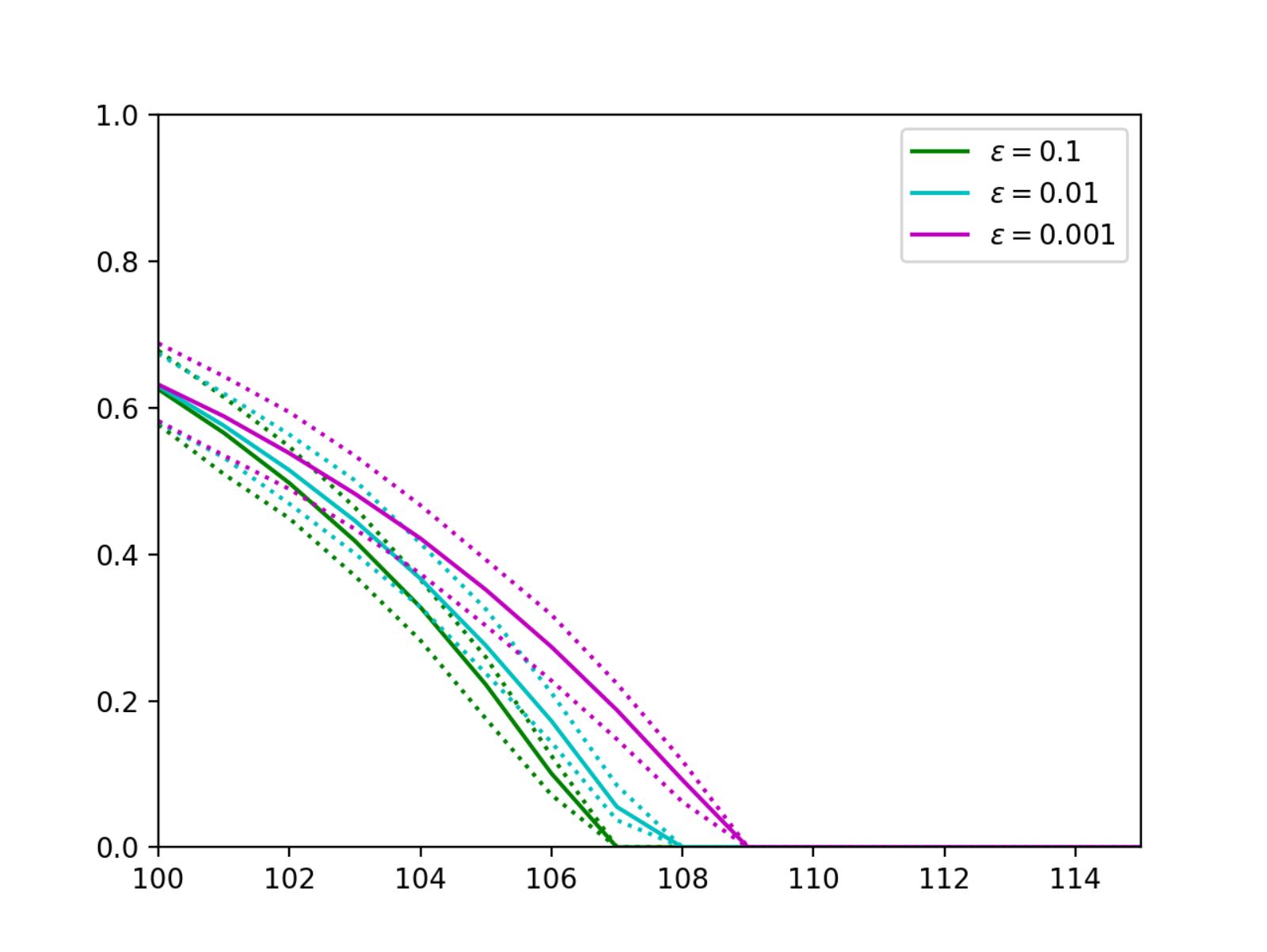}
  }
  ~
  \subfigure[rtnorm; $v = 0.05$]{
    \label{fig:v005rr}
    \includegraphics[width=0.45\textwidth]{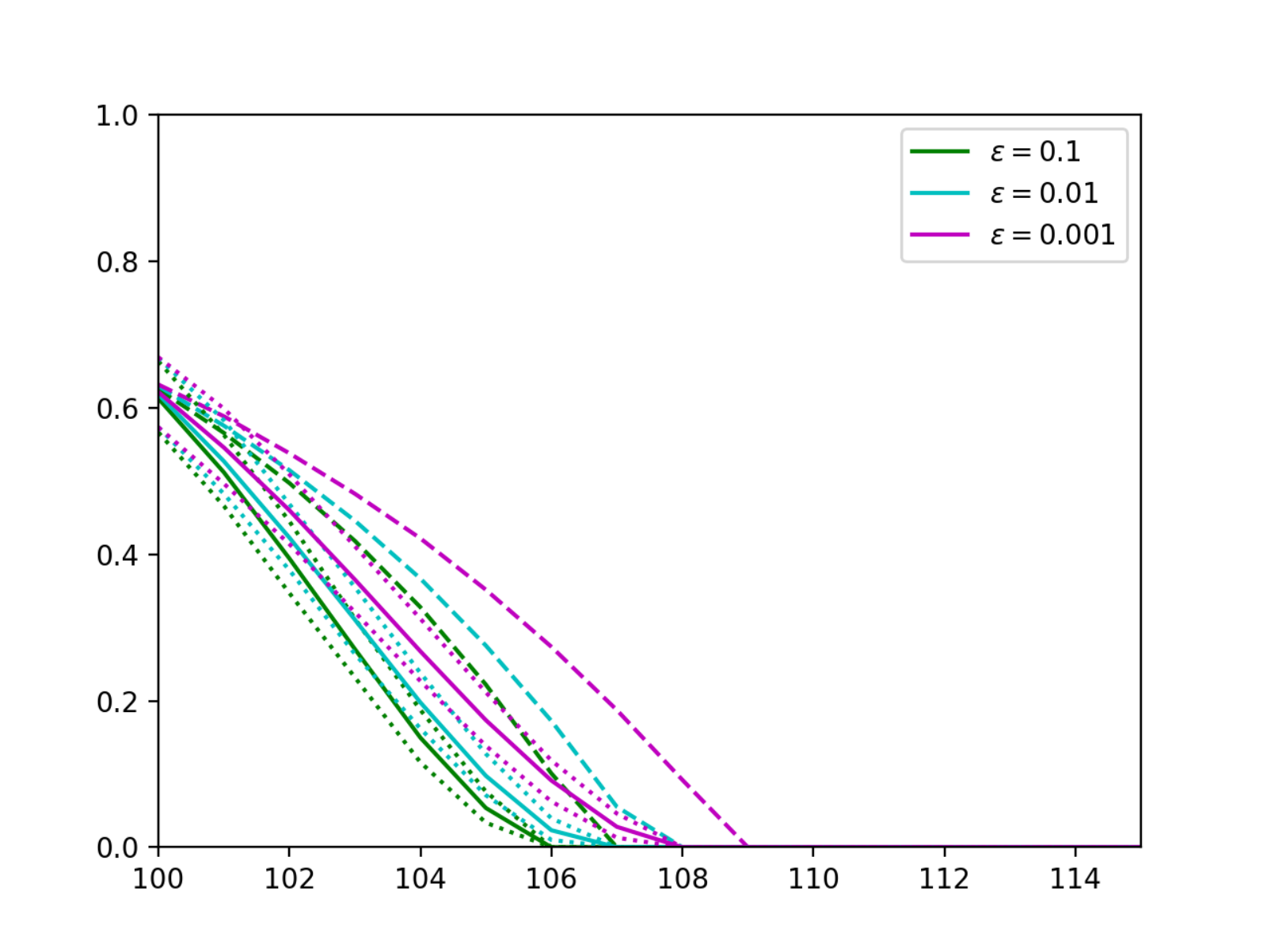}
  }
  \caption{Monte Carlo estimate for probability of success $P(z_{\ast} > \frac{x\bar{z_{\ast}}}{100})$, with bounds also calculated by Monte Carlo sampling. Plots are shown in each subfigure for three values of $\epsilon$, where $\delta \sim U(0,\epsilon)$ or $\delta \sim G(0,\epsilon,\epsilon/2,\sqrt{\epsilon^{2}/12})$. The Monte Carlo bounds are estimated by repeatedly ($100000$ times) performing $N = 100000$ independent samplings of $\delta$ and calculating $P($\texttt{success}$)$ for each iteration -- then selecting the minimum and maximum from the resulting distribution of $P($\texttt{success}$)$. The average $P($\texttt{success}$)$ for $\delta \sim U(0,\epsilon)$ are shown in all plots, for ease of comparison between results using the two different distributions of $\delta$.}
 \label{fig:sampling}
\end{figure}

%NOTE: see https://en.wikipedia.org/wiki/Probability_bounds_analysis ?
The value of $P(\tt{success})$ is dependent on sampling $\delta \sim U(0,\epsilon)$ (or $\sim G(0,\epsilon,\epsilon/2,\sqrt{\epsilon^{2}/12})$); thus, any new sampling of $\delta$ will likely yield slightly different results for $P(\tt{success})$.  A natural question is: \textit{How different?}
This question can be directly addressed if the distribution of calculated $P($\texttt{success}$)$ is known.
The most common way of resolving the output distribution for a given function, given random variable inputs, is (maybe unsurprisingly) Monte Carlo sampling.
Thus, by running the code in Figures \ref{code:sampling} and \ref{code:pof} repeatedly, we can begin to resolve what the output distribution of $P($\texttt{success}$)$ looks like for different samplings of $\delta$ (Figure \ref{fig:MCdist}).  A simple way to identify lower and upper bounds on $P($\texttt{success}$)$ is to use the minimum and maximum values found in the sampled distribution of $P($\texttt{success}$)$.  In Figure \ref{fig:sampling}, we see the sampled bounds obtained thusly for $100000$ independent calculations of $P($\texttt{success}$)$, each with $N = 100000$ independent samplings of $\delta$.
The results show an obvious trend of decreasing $P(\tt{success})$ for increasing $\epsilon$. It also appears that $P(\tt{success})$ decreases with decreasing $v$, to the point where $P($\texttt{success}$)$ is $0$. Additionally, when one compares sampling of $\delta \sim U(0,\epsilon)$ versus $\delta \sim G(0,\epsilon,\epsilon/2,\sqrt{\epsilon^{2}/12})$, sampling from the uniform distribution yields slightly larger $P($\texttt{success}$)$. 

In engineering, finance, and many fields of science, there is a goal of minimizing risk or optimizing design for a given system under uncertainty.  Minimization of risk, or optimization of design, can be performed without full knowledge of the output distribution of the system, if the most likely and extremal output values of the system are known.  Thus, in addition to the most likely value (as found in the previous section), we are also interested in determining the greatest upper bound and the least lower bound, with certainty.  Knowing the greatest upper bound and least lower bound for probability of success enables us to make guaranteed decisions about the system. If we don't know the extrema, then we fall into the realm of having some measure of confidence in our predictions (as opposed to being able to guarantee with certainty).
With respect to the system defined by (\ref{eq:burgers_steady}), have we thus far found the \textit{optimal} lower and upper bounds for $P($\texttt{success}$)$?  No, likely not.  Given more time, and more samplings of $\delta$, it is likely that eventually a more extremal value of $P($\texttt{success}$)$ will be found. Given the current computational cost of $2$ hr $\times$ $100000$ ($22.8$ years, without hierarchical parallel computing -- or alternately, roughly $35$ days on 256 compute nodes with 8 cores each), it's unlikely that better bounds will be found in a reasonable time and reasonable use of computing resources. Ultimately, to find the optimal lower and upper bounds, we need a change in approach.
%\begin{compactitem}
% \item Discuss MC calculation of bounds on $P(z_{\ast} > \frac{x\bar{z_{\ast}}}{100})$.
% \item MC bounds for different distributions of $\delta$ (non-uniform). %(#3)
% \item Discuss timings and compuational complexity for the above.
% \item Discuss use of std and other estimations of bounds.
%\end{compactitem}

% \newpage

\input ./ouq.tex % \section{An OUQ Refresher}

\section{Rigorous Upper Bounds on Probability of Success}
The OUQ problem at hand is to find the global maximum of the probability function $P($\texttt{success}$)$, where 
\texttt{success} is the success/failure criterion on solutions of (\ref{eq:burgers_steady}) defined in (\ref{eqn:success}).
Substituting into (\ref{eq:ext_E_qoi}), we have $q(x,y) = \tt{success}$, $g(X) = \tt{model}$,
\begin{subequations}
	\label{eq:qoi}
	\begin{align}
		\underline{Q}(\mathcal{A}) & \defeq \inf_{\mu \in \mathcal{A}} \P_{\mu}[\tt{success}], 
		\quad\text{and}\quad \\
		\overline{Q}(\mathcal{A}) & \defeq \sup_{\mu \in \mathcal{A}} \P_{\mu}[\tt{success}].
	\end{align}
\end{subequations}
Unlike sampling methods, OUQ doesn't require the distribution of the input random variables to be specified; we instead represent the unknown probability distribution by a product measure (composed of discrete support points, each with an accompanying weight), and any constraining information
% on the moments (e.g.\ \texttt{mean}, \texttt{std}) of the random variables
is specified in the set $\mathcal{A}$.
Specifically, we will define $\mathcal{A}$ as constraints on the feasible set of solutions, where:\begin{equation}
\label{eqn:mathcala} 
\mathcal{A} = \left\{ (g, \mu) \,\middle|\,
    \begin{matrix}
        g = \tt{model} \,:\, \mu \in [\tt{lb}, \tt{ub}] \to \mathbb{R}, \\
	\mu = \sum_{i = 0}^{3} w_{i} \delta_{i}, \\
        \sum_{i = 0}^{3} w_{i} = 0, \\
        \mathbb{E}_{\mu}[g] = \tt{z\_mean}, \\
        \bar{\mu} = \tt{d\_mean}
    \end{matrix} \right\}
\end{equation}
This imposes a mean constraint on $\delta$, a mean constraint on $z_{\ast}$, and
normalizes the weights. We will use three support points to represent the
product measure in the problem.
To solve this OUQ problem, we will need to write the code
for the bounds, constraints, and objective function -{}- then
we will plug the code into a global optimizer.

The optimization is performed by leveraging the \texttt{product\_measure} class in the \texttt{mystic} software.  The code in Figures \ref{code:ouqvar}, \ref{code:ouqfail}, \ref{code:ouqconstrain}, and \ref{code:ouqpof} detail how to use \texttt{mystic} to find the upper bound on
$P(\tt{success})$ at $\delta \sim U(0,\epsilon)$ with $v = 0.1$ and $\epsilon = 0.1$ under the constraints defined in $\mathcal{A}$.
While the code in Figure \ref{code:ouqvar} is specifically for \texttt{v = 0.1} and \texttt{eps = 0.1}, editing those two variables -- and looking up the corresponding $\bar{z_{\ast}}$, $\sigma_{z_{\ast}}$, $\bar{\delta}$, and $\sigma_{\delta}$ in Table \ref{tab:multiepsilon} -- will enable calculations of the upper bound of $P($\texttt{success}$)$ for other values of $v$ and $\epsilon$. Similarly, making the modification \texttt{MINMAX = 1} will produce the lower bound instead of the upper bound.
In general, the code in Figure \ref{code:ouqconstrain} imposes the constraints, $\mathcal{A}$, on the product measure used in the optimization problem. Specifically \texttt{\_constrain} imposes a mean constraint for $\delta$ on the product measure \texttt{c}, and \texttt{constrain} imposes a mean constraint for $z_{\ast}$ on \texttt{c} and additionally normalizes the measure's weights. These constraints are used by \texttt{solver} (in Figure \ref{code:ouqpof}) to specify the feasible set of solutions that the optimizer can search.
In terms of implementation, simple constraints on input parameter ranges are imposed with the \texttt{SetStrictRanges} method, and all other constraints are imposed with the \texttt{SetConstraints} method. 
If we add or subtract new constraints in $\mathcal{A}$, we then will have to modify the corresponding code in Figure \ref{code:ouqconstrain}.

\begin{figure}[htbp!]
  \centering
  \subfigure[$v = 0.1$]{
    \label{fig:v01}
    \includegraphics[width=0.45\textwidth]{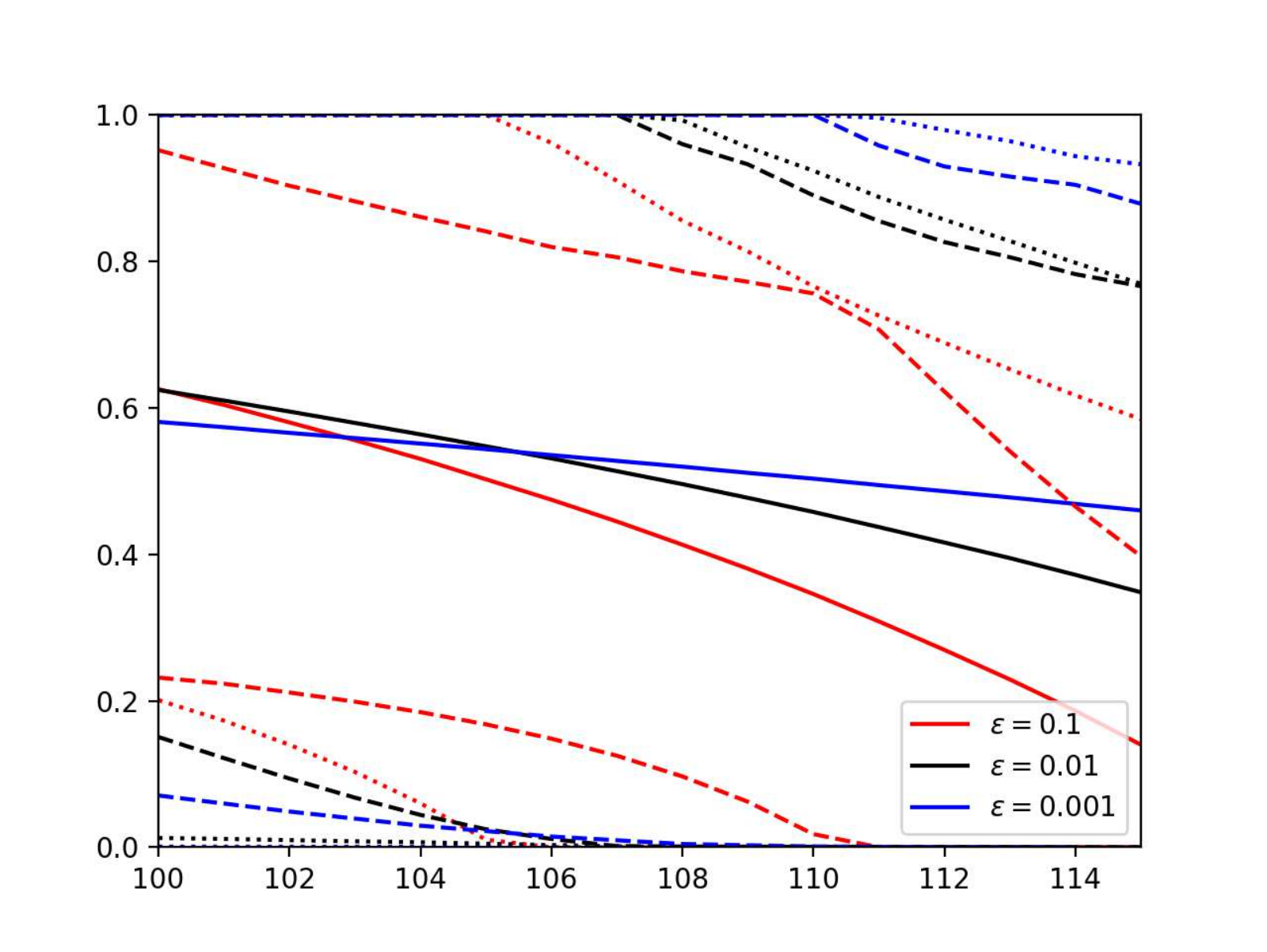}
  }
  ~
  \subfigure[$v = 0.05$]{
    \label{fig:v005}
    \includegraphics[width=0.45\textwidth]{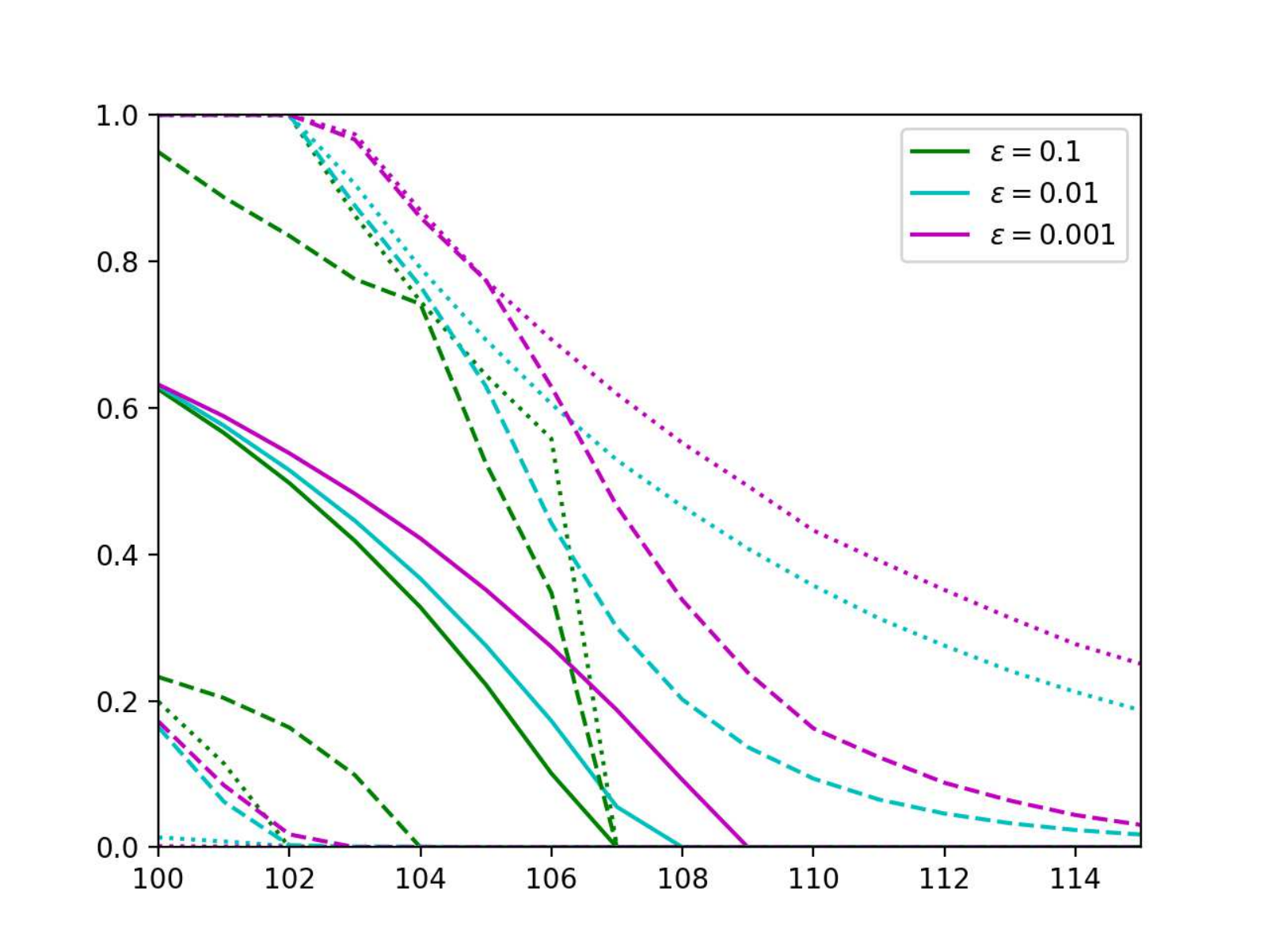}
  }
  \caption{Monte Carlo estimate for probability of success $P(z_{\ast} > \frac{x\bar{z_{\ast}}}{100})$, with bounds calculated by OUQ. Plots are shown in each subfigure for three values of $\epsilon$, where $\delta \sim U(0,\epsilon)$. OUQ bounds are calculated with a mean constraint on $\delta$ (dotted lines), and with a mean constraint and a variance constraint on $\delta$ (dashed lines). The presence of new information (i.e. the variance constraint) can be seen to tighten the bounds. As expected, $P($\texttt{success}$)$ decreases with increasing $x$.}
  \label{fig:OUQ}
\end{figure}

Results are plotted in Figure \ref{fig:OUQ} for:\begin{equation}
\label{eqn:A_meanD_varD} 
\mathcal{A} = \left\{ (g, \mu) \,\middle|\,
    \begin{matrix}
        g = \tt{model} \,:\, \mu \in [\tt{lb}, \tt{ub}] \to \mathbb{R}, \\
	\mu = \sum_{i = 0}^{3} w_{i} \delta_{i}, \\
        \sum_{i = 0}^{3} w_{i} = 0, \\
        \sigma_{g} = \tt{d\_std}, \\
        \bar{\mu} = \tt{d\_mean}
    \end{matrix} \right\}
\end{equation} and:\begin{equation}
\label{eqn:A_meanD} 
\mathcal{A} = \left\{ (g, \mu) \,\middle|\,
    \begin{matrix}
        g = \tt{model} \,:\, \mu \in [\tt{lb}, \tt{ub}] \to \mathbb{R}, \\
	\mu = \sum_{i = 0}^{3} w_{i} \delta_{i}, \\
        \sum_{i = 0}^{3} w_{i} = 0, \\
        \bar{\mu} = \tt{d\_mean}
    \end{matrix} \right\}
\end{equation} for $\epsilon \in 0.1, 0.01, 0.001$ and $v \in 0.05, 0.1$.
From the figure, we can see that increasing $\epsilon$ tightens the bounds.
This is in part due to the selection of \texttt{d\_mean} and \texttt{d\_std}
from Table \ref{tab:multiepsilon}, where the values for each $\epsilon$ and
$v$ are different. Perhaps a better, or just an alternate, comparison would be
to use a fixed \texttt{d\_mean} and \texttt{d\_std} for all optimizations on
the set of $\epsilon$ and $v$ above.
It's also notable that when comparing the constraints defined in (\ref{eqn:A_meanD_varD}) versus (\ref{eqn:A_meanD}), we see that the presence of new information (the variance constraint) causes a tightening of the bounds.

\begin{figure}[htbp!]
  \centering
  \subfigure[$x = 100 + dx, dx \in \lbrack 0,15 \rbrack$]{
    \label{fig:constraints}
    \includegraphics[width=0.45\textwidth]{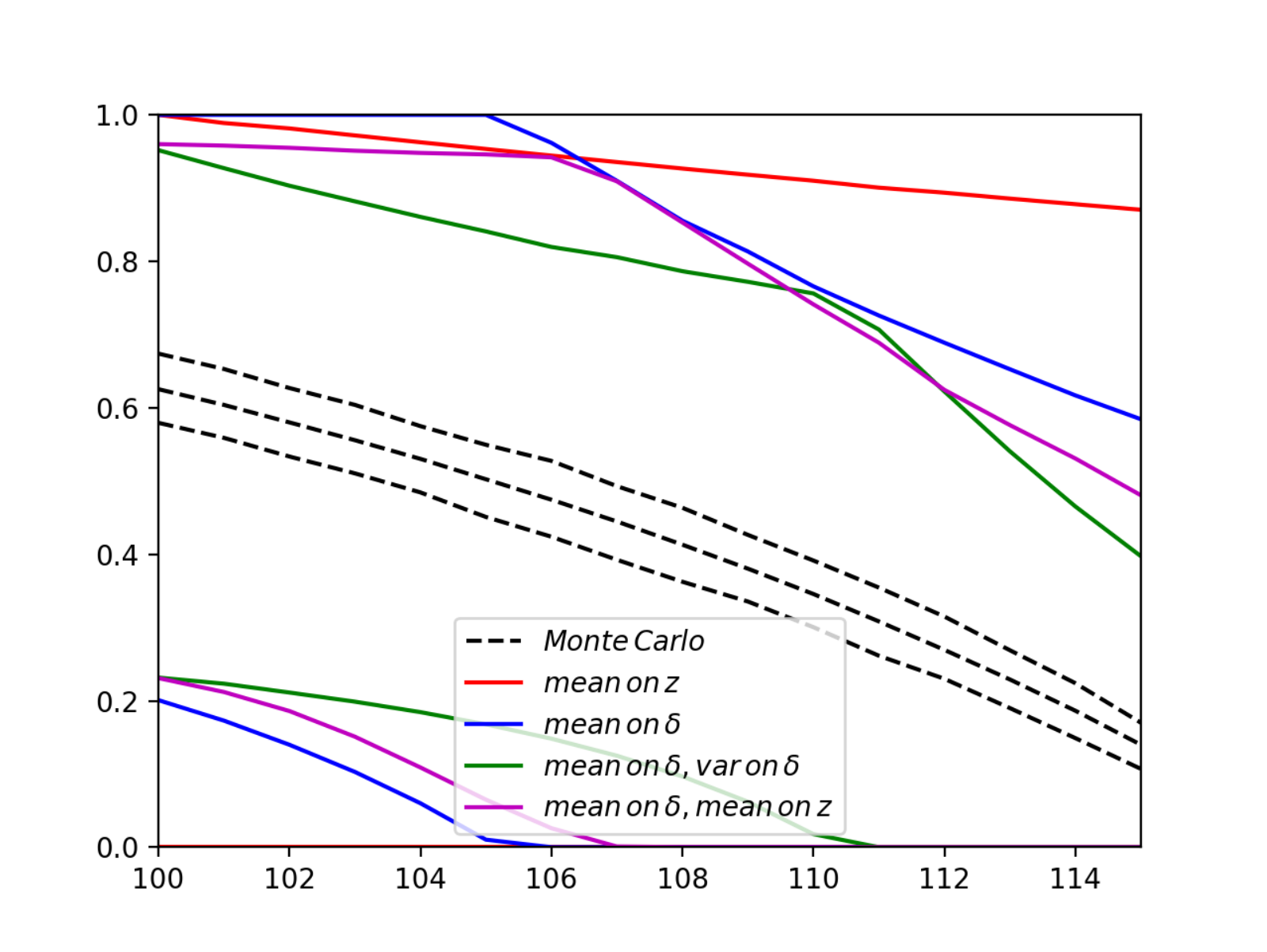}
  }
  ~
  \subfigure[Monte Carlo; $P(z_{\ast} > \bar{z_{\ast}})$]{
    \label{fig:MCdist}
    \includegraphics[width=0.45\textwidth]{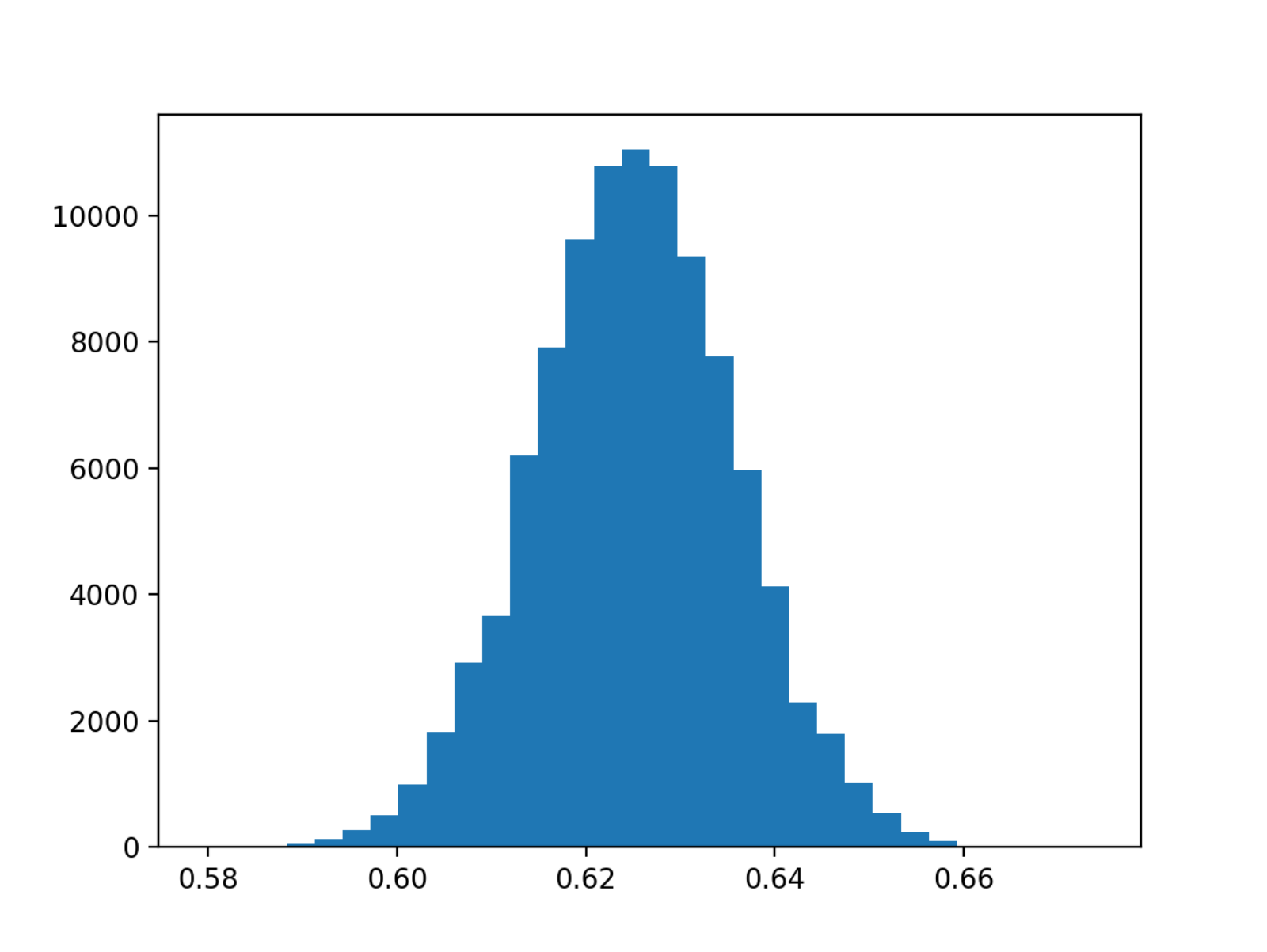}
  }
  \caption{Monte Carlo estimate for probability of success $P(z_{\ast} > \frac{x\bar{z_{\ast}}}{100})$ at $\delta \sim U(0,0.1)$ and $v = 0.1$, with bounds calculated by OUQ.
OUQ bounds are calculated with a mean constraint on $\delta$, a mean constraint and a variance constraint on $\delta$, a mean constraint on $z$, and a mean constraint on $\delta$ and a mean constraint on $z$. The effect of having different information constraints can be seen on the calculated lower and upper bounds. More specifically, the presence of additional information can be seen to generally tighten the bounds. Notice the striking difference between the bounds calculated with OUQ and Monte Carlo sampling. The plot of the distribution of $100000$ Monte Carlo sampled values of $P(z_{\ast} > \frac{x\bar{z_{\ast}}}{100})$ emphasize how much more capable OUQ is at describing behavior governed by rare events.}
  \label{fig:punchline}
\end{figure}

In Figure \ref{fig:constraints}, we further explore the notion that adding new constraining information to $\mathcal{A}$ tightens the bounds. We examine the OUQ bounds calculated in four cases: (\ref{eqn:mathcala}), (\ref{eqn:A_meanD_varD}), (\ref{eqn:A_meanD}), and\begin{equation}
\label{eqn:A_meanZ} 
\mathcal{A} = \left\{ (g, \mu) \,\middle|\,
    \begin{matrix}
        g = \tt{model} \,:\, \mu \in [\tt{lb}, \tt{ub}] \to \mathbb{R}, \\
        \mu = \sum_{i = 0}^{3} w_{i} \delta_{i}, \\
        \sum_{i = 0}^{3} w_{i} = 0, \\
        \mathbb{E}_{\mu}[g] = \tt{z\_mean}, \\
    \end{matrix} \right\}
\end{equation}
and also compare to the bounds calculated with Monte Carlo sampling.
We see very clearly when comparing results for (\ref{eqn:mathcala}) and (\ref{eqn:A_meanZ}) that the presence of a new constraint tightens the bounds.  Similarly for (\ref{eqn:A_meanD_varD}) and (\ref{eqn:A_meanD}). Additionally, the difference in the bounds calculated with OUQ compared to the bounds found using Monte Carlo sampling is striking.  The Monte Carlo bounds are much tighter -- however, are not the \textit{optimal} upper and lower bounds on $P($\texttt{success}$)$ -- they are merely estimates. The bounds calculated from OUQ are rigorous \textit{optimal} upper and lower bounds.  Upon reflection, this should be expected, as the OUQ optimization drives the solver to the extrema of $P($\texttt{success}$)$, while Monte Carlo sampling will struggle to select values at the tails of the distribution.

The computational cost of OUQ versus Monte Carlo should also be noted. Not only does OUQ provide better bounds, but it does so at a greatly reduced computational cost. The cost of an OUQ calculation is strongly dependent on the form of \texttt{constrain} (see Table \ref{tab:timing}). If an equation exists that maps the space to the feasible set, then the optimization typically is faster than without the constraint; however, for numerically imposed constraints (e.g.\ a nested optimization), the cost tends to be much larger. Constraints on the output of $g(X)$ also tend to be expensive, as imposing a constraint on the output will require several evaluations of $g(X)$ each iteration, in addition to the required number of evaluations by the optimization algorithm itself.
The timings in Table \ref{tab:timing} are for the \texttt{DifferentialEvolution2} solver, which is a robust global optimizer, but can be quite slow to converge.Time to solution could be improved by judicious tuning of the optimization settings (e.g.\ \texttt{npop}, \texttt{ngen}, \texttt{crossover}, etc), or use of an ensemble (e.g.\ \texttt{lattice}) of fast solvers as used in Figure \ref{code:exactly}.
Timings are presented for a 2.7 GHz Macbook with an Intel Core i7 processor and 16 GB 1600 MHz DDR3 memory, running python 2.7.14 and \texttt{mystic} 0.3.2.
Notably, the OUQ calculations can be performed in a reasonable time on a relatively good laptop, while the estimation of upper and lower bounds using Monte Carlo sampling requires institutional computing resources. 

One will notice that the OUQ section did not specify whether
$\delta \sim U(0,\epsilon)$ or 
$\delta \sim G(0,\epsilon,\epsilon/2,\sqrt{\epsilon^{2}/12})$ was used.
This is because definitions of $\mathcal{A}$ specified constraints on the moments of $\delta$ and $z_{\ast}$, and did not specify a particular input or output distribution. While it is possible to specify the distribution of the input (or output) random variables in $\mathcal{A}$, it is in reality unlikely that this information is actually known -- more commonly, only the moment information is known about the input and output variables due to measurements taken on the system. OUQ is built to handle these kinds of physical constraints explicitly, while Monte Carlo cannot.

%\begin{compactitem}
% \item Discss problem formulation, including code, for min and max.
% \item Describe how A maps onto the code in the figures.
% \item Discuss potential effect on bounds of choosing the same \texttt{mean(d)}, \texttt{var(d)} constraint values (as opposed to values calculated from each MC run) for all $\epsilon$. %e.g.\ for same constraining value, 0.001 should be wider than 0.1
% \item Discss problem variants for the different information constraints.
% \item Discss computational complexity and cost.
%end{compactitem}

%\newpage
%\COMMENT{ % find OUQ bounds for Burgers' equation with perturbed boundary
\begin{figure}[htbp!]
  \centering
  \begin{minipage}{1.00\textwidth}
% \begin{\outputtextsize}
\begin{verbatim}
### OUQ_supersensitive.py :: calculate bounds for Burgers' equation ###
from numpy import inf
from mystic.math.measures import split_param
from mystic.math.discrete import product_measure
from mystic.math.stats import meanconf
from mystic.math import almostEqual
from mystic.solvers import DifferentialEvolutionSolver2
from mystic.termination import Or, VTR, ChangeOverGeneration as COG
from mystic.monitors import VerboseMonitor
from exact_supersensitive import solve

v, eps = 0.1, 0.1
MINMAX = -1  # {'maximize':-1, 'minimize':1}
nx = 3  # number of support points for delta
npts = (nx,)  # dimensionality of support for the random variables
pcnt = 0.00  # percent increase for z_mean
N = 100000  # number of realizations of delta
z_mean = 0.614420266306; z_std = 0.105011650866
d_mean = 0.0499427947961; d_std = 0.0289104773543

# bounds on weights and positions for support points for delta
w_lower = [0.0]; w_upper = [1.0]; x_lower = [0.0]; x_upper = [eps]
lb = (nx * w_lower) + (nx * x_lower); ub = (nx * w_upper) + (nx * x_upper)

# define constraints on mean and std of inputs and outputs
z_range = meanconf(z_std,N); d_range = meanconf(d_std,N)
target = (z_mean, d_mean,); error = (z_range, d_range,)

\end{verbatim}
  \end{minipage}
  \caption{
Global variables and imports required for the
calculation of lower and upper bounds of the
probability of success $P(z_{\ast} > \frac{x\bar{z_{\ast}}}{100})$
at $\delta \sim U(0,0.1)$ and $v = 0.1$.
The values of \texttt{z\_mean}, \texttt{z\_std}, \texttt{d\_mean}, and \texttt{d\_std} are taken from Table
\ref{tab:multiepsilon}.
\texttt{MINMAX}=-1 indicates an upper bound calculation is to be performed.
\texttt{pcnt} is percent increase from $\bar{z_{\ast}}$ (i.e. $\frac{x}{100} - 1$).
\texttt{z\_mean} and \texttt{d\_mean} are used as \texttt{target} values for the
moment constraints. An \texttt{error} of up to the 95\% confidence interval
(\texttt{z\_range} and \texttt{d\_range}) is allowed in the mean constraints
on $z$ and $\delta$.
\texttt{lb} and \texttt{ub} define the lower and upper limits for the weights and positions for the support points for $\delta$. The \texttt{len(lb)} is $6$ because we are using \texttt{nx = 3} support points. The variables and imports defined here are used in the code in Figures \ref{code:ouqfail}, \ref{code:ouqconstrain}, and \ref{code:ouqpof}.
}
  \label{code:ouqvar}
\end{figure}

\begin{figure}[htbp!]
  \centering
  \begin{minipage}{1.00\textwidth}
% \begin{\outputtextsize}
\begin{verbatim}
### OUQ_supersensitive.py (cont'd) :: calculate bounds for Burgers' equation ###
# solve Burgers' equation for z given "parameter vector" [v,delta]
zsolve = lambda rv: solve(*rv)[1][0]

# the model function (for fixed v)
model = lambda rv: zsolve((v,)+rv)

# "success" and "failure" indicator functions
success = lambda z,zave: z > zave
failure = lambda rv: not success(model(rv), (1 + pcnt)*z_mean)

\end{verbatim}
  \end{minipage}
  \caption{
Define model function and failure function for the
calculation of lower and upper bounds of the
probability of success $P(z_{\ast} > \frac{x\bar{z_{\ast}}}{100})$
at $\delta \sim U(0,0.1)$ and $v = 0.1$.
\texttt{failure} is defined as when \texttt{(1+pcnt)*z\_mean} is greater than
$z_{\ast}$, where $z_{\ast}$ is determined by the analytical solution to (\ref{eq:burgers_steady}) (using \texttt{solve}).
The \texttt{model} function is used in both the constraints (Figure \ref{code:ouqconstrain}) and the objective function (Figure \ref{code:ouqpof}),
and solves $z_{\ast}$ for
given inputs $[v,\delta]$, where $v=0.1$ and $\delta$ is a random variable. 
}
  \label{code:ouqfail}
\end{figure}

\begin{figure}[htbp!]
  \centering
  \begin{minipage}{1.00\textwidth}
% \begin{\outputtextsize}
%%% product_measure,almostEqual
%%% npts,model,target,error,x_lb,x_ub
\begin{verbatim}
### OUQ_supersensitive.py (cont'd) :: calculate bounds for Burgers' equation ###
x_lb = split_param(lb, npts)[-1]  # lower bounds on positions only
x_ub = split_param(ub, npts)[-1]  # upper bounds on positions only

def constraints(rv):
    c = product_measure().load(rv, npts)
    # impose norm on each discrete measure
    for measure in c:
        if not almostEqual(float(measure.mass), 1.0):
            measure.normalize()
    # impose expectation value and other constraints on product measure
    E = float(c.expect(model))
    if E > (target[0] + error[0]) or E < (target[0] - error[0]):
        c.set_expect((target[0],error[0]), model, (x_lb,x_ub), _constraints)
    return c.flatten()  # extract parameter vector of weights and positions

def _constraints(c):
    E = float(c[0].mean)
    if E > (target[1] + error[1]) or E < (target[1] - error[1]):
        c[0].mean = target[1]
    return c

\end{verbatim}
  \end{minipage}
  \caption{
Define constraints function(s) for the
calculation of lower and upper bounds of the
probability of success $P(z_{\ast} > \frac{x\bar{z_{\ast}}}{100})$
at $\delta \sim U(0,0.1)$ and $v = 0.1$.
The function \texttt{constraints} imposes constraints on the parameter vector \texttt{rv} by normalizing the weights of each discrete measure, and imposing
mean constraints on $z_{\ast}$ and $\delta$ (in \texttt{\_constraints}).
The constraints defined here are used by the optimizer in Figure \ref{code:ouqpof}.
}
  \label{code:ouqconstrain}
\end{figure}

\begin{figure}[htbp!]
  \centering
  \begin{minipage}{1.00\textwidth}
% \begin{\outputtextsize}
%%% product_measure,almostEqual,inf
%%% npts,model,target,error,failure,MINMAX
%%% lb,ub,npop,maxiter,maxfun,npts,constraints,
%%% convergence_tol,ngen,crossover,scaling,objective,MINMAX
\begin{verbatim}
### OUQ_supersensitive.py (cont'd) :: calculate bounds for Burgers' equation ###
def objective(rv):
    c = product_measure().load(rv, npts)
    E = float(c[0].mean)
    if E > (target[1] + error[1]) or E < (target[1] - error[1]):
        return inf
    E = float(c.expect(model))
    if E > (target[0] + error[0]) or E < (target[0] - error[0]):
        return inf
    return MINMAX * c.pof(failure)

# solver parameters
npop, maxiter, maxfun = 10, 1000, 1e+6
convergence_tol = 1e-6; ngen = 10; crossover = 0.9; scaling = 0.4

# configure solver and find extremum of objective
solver = DifferentialEvolutionSolver2(len(lb),npop)
solver.SetRandomInitialPoints(min=lb,max=ub)
solver.SetStrictRanges(min=lb,max=ub)
solver.SetEvaluationLimits(maxiter,maxfun)
solver.SetGenerationMonitor(VerboseMonitor(1, 1, npts=npts))
solver.SetConstraints(constraints)
solver.SetTermination(Or(COG(convergence_tol, ngen), VTR(convergence_tol, -1.0)))
solver.Solve(objective, CrossProbability=crossover, ScalingFactor=scaling)
print "results: %s" % MINMAX * solver.bestEnergy
print "# evals: %s" % solver.evaluations
solver.SaveSolver('solver.pkl')

\end{verbatim}
% \end{\outputtextsize}
  \end{minipage}
  \caption{
Define an objective function, then configure the optimizer and solve for lower and upper bounds of the probability of success $P(z_{\ast} > \frac{x\bar{z_{\ast}}}{100})$
at $\delta \sim U(0,0.1)$ and $v = 0.1$.
\texttt{objective} essentially constructs a \texttt{product\_measure} from the parameter vector, and then calculates the probability of failure. \texttt{objective} also includes two safety checks that throw out any solutions in the unlikely case that the constraints solver \texttt{constraints} fails to impose the constraints.
The \texttt{solver} uses a Differential Evolution algorithm to find the minimum of \texttt{objective} with respect to bounds \texttt{lb} and \texttt{ub} and constraints \texttt{constraints}. Custom termination conditions are used to signal when \texttt{solver} should stop, while \texttt{SaveSolver} preserves the final state of \texttt{solver} to a file.
Results are shown in Figures \ref{fig:OUQ} and \ref{fig:constraints}.
}
  \label{code:ouqpof}
\end{figure}
%}%\END COMMENT

% show/plot results of MC PoF, OUQ PoF bounds

\begin{table}[htbp!]
  \centering
  \begin{tabular}{c c c c c}
    \hline
    $ $ & $\bar{\delta}$ & $\bar{\delta}\, \sigma_{\delta}$ & $\bar{z_{\ast}}$ & $\bar{\delta}\, \bar{z_{\ast}}$ \\
    \hline
    \texttt{upper} & \texttt{0:03:02} & \texttt{0:20:08} & \texttt{1 day, 22:47:18} & \texttt{18:41:05} \\
    \texttt{lower} & \texttt{0:00:03} & \texttt{0:01:39} & \texttt{0:00:00} & \texttt{1 day, 09:54:17} \\
    \hline
  \end{tabular}
  \caption{Sample timings for OUQ optimizations defined by (\ref{eq:qoi} - \ref{eqn:A_meanZ}), where the name of the column is an indication of the unique constraint applied (e.g.\ $\bar{\delta}$ is a mean constraint on $\delta$).
The vast difference in timing is due to both the topological nature of the $P(\tt{success})$ for the chosen values of $v$ and $\epsilon$, and the nature of the constraints in $\mathcal{A}$. Note that the optimizations leverage $100000$ Monte Carlo realizations of $\delta$ from Figure \ref{code:sampling}. The values in the table correspond to a single point on Figure \ref{fig:OUQ}.}
  \label{tab:timing}
\end{table}

\COMMENT{ % timings for OUQ POF
\begin{compactitem}
  \item 110000 Monte Carlo realizations of \delta
took: 2:36:33.555913, where a buffer of 10000 was used
but only 778 did not make the fit tolerance of 1e-09.
  \item OUQ mean(d) max: (1.08\%) took: 0:03:02.461392  
  \item OUQ mean(d) min: (1.08\%) took: 0:00:02.806286
  \item OUQ mean(d)var(d) max: (1.08\%) took: 0:20:08.086603
                                   (3N) took: 0:00:44.717957
  \item OUQ mean(d)var(d) min: (1.08\%) took: 0:01:38.979158
                                   (3N) took: 0:22:20.824934
  \item OUQ mean(z) max: (1.08\%) took: 1 day, 22:47:17.974492
  \item OUQ mean(z) min: (1.08\%) took: 0:00:00.000000
  \item OUQ mean(d)mean(z) max: (1.08\%) took: 9:32:11.946637 + 9:08:53.188908
  \item OUQ mean(d)mean(z) min: (1.08\%) took: 12:02:06.115143 + 21:52:10.772389
\end{compactitem}
} % END COMMENT

% \newpage

\section{Discussion and Future Work}
In this manuscript, we demonstrate utility of the OUQ approach to understanding
the behavior of a system that is governed by partial differential equations
(more specifically, by Burgers' equation).
In particular, we solve the problem of predicting shock location when
we know only bounds on viscosity and on the initial conditions.
We calculate bounds on $P(\tt{success})$ given the uncertainty in $\delta$, using both a standard Monte Carlo sampling approach, and OUQ.
The results highlight the stark contrast in the ability of each approach to elucidate the rare-event behavior that occurs at the tails of the unknown distribution of $\delta$. OUQ uses numerical optimization to discover the governing behavior at the tails of the distribution for input random variables. In contrast, using Monte Carlo sampling to find accurate extremal behavior in many cases is computationally infeasible. OUQ finds \textit{optimal} upper and lower bounds on the quantity of interest with regard to the constraining information $\mathcal{A}$ in the system, while Monte Carlo can only provide an estimate of the upper and lower bounds -- and as mentioned above, requires assumptions on the form of the distribution for the input random variables.
In short, OUQ provides more accurate bounds at a lower computational cost,
and additionally 
can take advantage of solution-constraining information that Monte Carlo
cannot. Since OUQ requires fewer assumptions on the form of the inputs,
the predicted bounds are more rigorous than those obtained with Monte Carlo.

OUQ is a very capable but relatively new approach, and has not yet been applied in many fields of science and engineering.
This manuscript is the first application of OUQ (that we know of) to a complex system governed by partial differential equations.  In our example, we had moment constraints on the input and output variables; however, OUQ is a general formalism, and has been shown to handle other types of constraining physical and statistical information (e.g.\ legacy data and associated uncertainties, uniqueness, further moments) \cite{OSSMO:2011, Sullivan:2013, Kamgaouq:2014}.
Adding further constraining information often has the drawback of making an OUQ calculation more expensive (as complex nested numerical optimizations may be needed to produce the space of feasible solutions). The cost of the optimization is roughly due to the most expensive inner optimization times the number of iterations. Hence, there is a trade-off of adding information and the potential computational cost/complexity \cite{MOSSO:2010}.

In this manuscript, we used OUQ to solve for the optimal probability of success. We direct the reader to earlier work (especially \cite{OSSMO:2011}) to see how to use OUQ for model certification and validation, design of experiments, and calculations of optimal model, model error, and other quantities relevant to risk and statistics under uncertainty. We will close by summarizing a few of the aforementioned applications of OUQ.
For example, finding the optimal model in OUQ is essentially finding the functional form that minimizes the worst upper and lower bounds on the model error for all possible models.
%The notion of optimal experiment design in OUQ depends on the question that one wants answered.  For example, a calculation of the best next experiment is done by finding the constraint (and hence required measurement or data point) to add to $\mathcal{A}$ that most increases the likelihood of success.
Experiment design with OUQ starts by selecting a quantity of interest to maximize or minimize, such as $P(\tt{success})$ or model error, and a base set $\mathcal{A}$ of assumptions, measurements, and data; then, as new information is added to $\mathcal{A}$, the change in bounds for the quantity of interest is calculated. The piece of information that most dramatically tightens the bounds indicates the best next experiment to perform.
The addition of new information to $\mathcal{A}$ can be manual, or optimizer-driven, with a OUQ calculation of bounds performed for each new $\mathcal{A}$.
It should be noted that many OUQ problem formulations require one or more outer optimization loops, and one or more inner optimization loops, to be able to impose the constraints in $\mathcal{A}$ -- thus OUQ can be significantly more computationally expensive than other methods that leverage strong approximations.
Since this manuscript is one of the first applications of OUQ to a complex physical system, we intend to follow this study with other examples of using OUQ for complex systems -- including examples of experiment design, and in the calculation of the optimal model, for a complex physical system.

% \newpage %XXX: to push bib to separate page

\section{Acknowledgements}

This work was supported by the U.S. Department of Energy, Office of
Science, Office of Advanced Scientific Computing Research,
specifically through the Mathematical Multifaceted Integrated
Capability Centers (MMICCS) program, under contract number
DE-SC0019303, and by the Uncertainty Quantification Foundation
under the Statistical Learning program. The Uncertainty Quantification
Foundation is a non-profit dedicated to the advancement of predictive
science through research, education, and the development and
dissemination of advanced technologies.

\bibliographystyle{plain}
\bibliography{./refs}

\end{document}

%% file: oed.tex
\section{OUQ for Model Validation and Experiment Design}

%The Problem:
The need to make possibly critical decisions about complex systems with limited information is common to nearly all areas of science and industry. The design problems used to formulate statistical estimators are typically not well-posed, and lack a well-defined notion of solution. As a result, different teams with the same goals and data will find vastly different statistical estimators, solutions, and notions of uncertainties. We claim that the recently developed OUQ theory, combined with recent software enabling fast global solutions of constrained non-convex optimization problems, provides a methodology for rigorous model certification, validation, and optimal design under uncertainty.
This OUQ-based methodology is especially pertinent in the context of beamline science, where data analysis generally (1) lacks full automation, (2) requires an expert, and (3) commonly utilizes phenomenological profile functions with loosely-defined statements of validity. The lack of automated analysis and rigorously validated models create a bottleneck in the scientific throughput of beamline science. For example, instrument scientists generally use very conservative back-of-the-envelope estimates to ensure sufficient beamtime is allocated to produce high-quality data. Consequently, a significant portion of the beamtime that could be used for new science is wasted.

%Proposed innovation:
The validation of a statistical estimators against real-world data generally is a laborious manual process. Significant human effort is spent cleaning and reducing data, as well as selecting, training, and validating the statistical estimators for each special-case. High-performance computing (HPC), however, offers an alternative approach: posing the construction of an optimal statistical estimators as an at-scale optimization problem. In this new paradigm, human intellectual effort goes into the design of an over-arching optimization problem, instead of being spent on manual derivation and validation of a statistical estimators. We propose a methodology for optimal design/planning under uncertainty, where our tools enable the determination of the impact of any new piece of information on all possible valid outcomes. Our vision is that HPC combined with OUQ will stimulate the emergence of a new paradigm where the parameter space to conduct key experiments is numerically designed and optimized.
We leverage the \texttt{mystic} optimization framework \cite{McKernsHungAivazis:2009, MSSFA:2011}
to impose constraining information from different sources on the search space of an optimizer, where the solution produced respects all constraints from the different data/models/etc.
Once all available information has been encoded as constraints, the numerical transform generated from these constraints guarantees that any solution is valid with respect to the constraints.
We can then use this methodology, for example, in the validation of phenomenological models against high-fidelity materials theory codes, or in the calculation of the minimum beamtime required to ensure data has met a high-quality threshold, or in other design of experiments questions.

% END COMMENT

%% file: ouq.tex
\section{An OUQ Refresher}

Rigorous quantification of the effects of epistemic and aleatoric
uncertainty is an increasingly important component of research studies
and policy decisions in science, engineering, and finance.  In the
presence of incomplete imperfect knowledge (sometimes called
\emph{epistemic uncertainty}) about the objects involved, and
especially in a high-consequence decision-making context, it makes
sense to adopt a posture of \emph{healthy conservatism}, i.e.\ to
determine the best and worst outcomes consistent with all available
knowledge.  This posture naturally leads to uncertainty quantification
(UQ) being posed as an optimization problem.  Such optimization
problems are typically high dimensional, and hence can be difficult
and expensive to solve computationally (depending on the nature of the
constraining information).

In \cite{OSSMO:2011, Sullivan:2013}, a theoretical framework is outlined for \emph{optimal uncertainty quantification} (OUQ), namely the calculation of optimal lower and upper bounds on probabilistic output quantities of interest, given quantitative information about (underdetermined) input probability distributions and response functions.  In their computational formulation \cite{MOSSO:2010, MSSFA:2011}, whose derivation we outline below, OUQ problems require optimization over discrete (finite support) probability distributions of the form
\begin{equation}
	\label{eq:reduced_measure_in_multiindex_form}
	\mu = \sum_{\bs{i} = \bs{0}}^{\bs{M}} w_{\bs{i}} \delta_{\bs{x}_{\bs{i}}} \text{,}
\end{equation}
where $\bs{i} = \bs{0}, \dotsc, \bs{M}$ is a finite range of indices, the $w_{\bs{i}}$ are non-negative weights that sum to $1$, and the $\bs{x}_{\bs{i}}$ are points in some input parameter space $\mathcal{X}$;  $\delta_{a}$ denotes the Dirac measure (unit point mass) located at a point $a \in \mathcal{X}$, i.e., for $E \subseteq \mathcal{X}$,
\[
	\delta_{a}(E) \defeq 
	\begin{cases}
		1\text{,} & \text{if $a \in E$,} \\
		0\text{,} & \text{if $a \notin E$.}
	\end{cases}
\]
Note that $\delta_{x_{i}}$ signifies a unit point mass, and not $\delta$ the position of the left boundary wall; however, in the next section, we will see for the example in this paper that $\delta_{x_{i}} = \delta$.

Many UQ problems such as certification, prediction, reliability estimation, and risk analysis can be posed as the calculation or estimation of an expected value, i.e.\ an integral, although this expectation (integral) may depend in intricate ways upon various probability measures, parameters, and models.  This point of view on UQ is similar to that of \cite{BP:1996}, in which formulations of many problem objectives in reliability are represented in a unified framework, and the decision-theoretic point of view of \cite{S:1995}.  In the presentation below, an important distinction is made between the ``real'' values of objects of interest, which are decorated with daggers (e.g.\ $g^{\dagger}$ and $\mu^{\dagger}$), versus possible models or other representatives for those objects, which are not so decorated.

The system of interest is a measurable \emph{response function} $g^{\dagger} \colon \mathcal{X} \to \mathcal{Y}$ that maps a measurable space $\mathcal{X}$ of \emph{inputs} into a measurable space $\mathcal{Y}$ of outputs.  The inputs of this response function are distributed according to a probability measure $\mu^{\dagger}$ on $\mathcal{X}$;  $\mathcal{P}(\mathcal{X})$ denotes the set of all probability measures on $\mathcal{X}$.  The UQ objective is to determine or estimate the expected value under $\mu^{\dagger}$ of some measurable \emph{quantity of interest} $q \colon \mathcal{X} \times \mathcal{Y} \to \R$, i.e.\
\begin{equation}
	\label{eq:E_qoi}
	\E_{X \sim \mu^{\dagger}}[q(X, g^{\dagger}(X))].
\end{equation}
The probability measure $\mu^{\dagger}$ can be interpreted in either a frequentist or subjectivist (Bayesian) manner, or even just as an abstract probability measure.  A typical example is that the event $[ g^{\dagger}(X) \in E ]$, for some measurable set $E \subseteq \mathcal{Y}$, constitutes some undesirable ``failure'' outcome, and it is desired to know the $\mu^{\dagger}$ probability of failure, in which case $q$ is the indicator function
\[
	q(x, y) \defeq
	\begin{cases}
		1 \text{,} & \text{if $y \in E$,} \\
		0 \text{,} & \text{if $y \notin E$.}
	\end{cases}
\]

In practice, the real response function and input distribution pair $(g^{\dagger}, \mu^{\dagger})$ are not known precisely.  In such a situation, it is not possible to calculate \eqref{eq:E_qoi} even by approximate methods such as Monte Carlo or other sampling techniques for the simple reason that one does not know \emph{which} probability distribution to sample, and it may be inappropriate to simply assume that a chosen model pair $(g^{\mathrm{m}}, \mu^{\mathrm{m}})$ is $(g^{\dagger}, \mu^{\dagger})$.  However, it may be known (perhaps with some degree of statistical confidence) that $(g^{\dagger}, \mu^{\dagger}) \in \mathcal{A}$ for some collection $\mathcal{A}$ of pairs of functions $g \colon \mathcal{X} \to \mathcal{Y}$ and probability measures $\mu \in \mathcal{P}(\mathcal{X})$.  If knowledge about which pairs $(g, \mu) \in \mathcal{A}$ are more likely than others to be $(g^{\dagger}, \mu^{\dagger})$ can be encapsulated in a probability measure $\pi \in \mathcal{P}(\mathcal{A})$ --- what a Bayesian probabilist would call a \emph{prior} --- then, instead of \eqref{eq:E_qoi}, it makes sense to calculate or estimate
\begin{equation}
	\label{eq:EE_qoi}
	\E_{(g, \mu) \sim \pi} \Bigl[ \E_{X \sim \mu}[q(X, g(X))] \Bigr].
\end{equation}
(A Bayesian probabilist would also incorporate additional data by conditioning to obtain the \emph{posterior} expected value of $q$.)

However, in many situations, either due to lack of knowledge or being in a high-consequence regime, it may be either impossible or undesirable to specify such a $\pi$.  In such situations, it makes sense to adopt a posture of \emph{healthy conservatism}, i.e.\ to determine the best and worst outcomes consistent with the available knowledge.  Hence, instead of \eqref{eq:E_qoi} or \eqref{eq:EE_qoi}, it makes sense to calculate or estimate
\begin{subequations}
	\label{eq:ext_E_qoi}
	\begin{align}
		\underline{Q}(\mathcal{A}) & \defeq \inf_{(g, \mu) \in \mathcal{A}} \E_{X \sim \mu}[q(X, g(X))]
		\quad\text{and}\quad \\
		\overline{Q}(\mathcal{A}) & \defeq \sup_{(g, \mu) \in \mathcal{A}} \E_{X \sim \mu}[q(X, g(X))].
	\end{align}
\end{subequations}
If the probability distributions $\mu$ are interpreted in a Bayesian sense, then this point of view is essentially that of the robust Bayesian paradigm \cite{Berger:1994} with the addition of uncertainty about the forward model(s) $g$.  Within the operations research and decision theory communities, similar questions have been considered under the name of distributionally robust optimization \cite{DY:2010, GS:2010, S:1995}.  Distributional robustness for polynomial chaos methods has been considered in \cite{NX:2011}.   Our interest lies in providing a UQ analysis for \eqref{eq:E_qoi} by the efficient calculation of the extreme values \eqref{eq:ext_E_qoi}.

%\COMMENT { % more OUQ details
In order to compute the extreme values of the optimization problems \eqref{eq:ext_E_qoi}, an essential step is finding finite-dimensional problems that are equivalent to (i.e.\ have the same extreme values as) the problems \eqref{eq:ext_E_qoi}.  This will be achieved by reducing the search over all measures to a search over discrete measures of the form \eqref{eq:reduced_measure_in_multiindex_form}.
A strong analogy can be made here with finite-dimensional linear programming:  to find the extreme value of a linear functional on a polytope, it is sufficient to search over the extreme points of the polytope;  the extremal scenarios of $\mathcal{A}$ turn out to consist of discrete functions and probability measures that are themselves far more singular than would ``typically'' be encountered ``in reality'' but nonetheless encode the full range of possible outcomes in much the same way as a polytope is the convex hull of its ``atypical'' extreme points.

One general setting in which a finite-dimensional reduction can be effected is that in which, for each candidate response function $g \colon \mathcal{X} \to \mathcal{Y}$, the set of input probability distributions $\mu \in \mathcal{P}(\mathcal{X})$ that are admissible in the sense that $(g, \mu) \in \mathcal{A}$ is a (possibly empty) generalized moment class.  More precisely, assume that it is known that the $\mu^{\dagger}$-distributed input random variable $X$ has $K$ independent components $(X_{0}, \dotsc, X_{K - 1})$, with each $X_{k}$ taking values in a Radon space\footnote{This technical requirement is not a serious restriction in practice, since it is satisfied by most common parameter and function spaces.  A \emph{Radon space} is a topological space on which every Borel probability measure $\mu$ is \emph{inner regular} in the sense that, for every measurable set $E$, $\mu(E) = \sup \{ \mu(K) \mid K \subseteq E \text{ is compact} \}$.  A simple example of a non-Radon space is the unit interval $[0, 1]$ with the lower limit topology.} $\mathcal{X}_{k}$;  this is the same as saying that $\mu^{\dagger}$ is a product of marginal probability measures $\mu^{\dagger}_{k}$ on each $\mathcal{X}_{k}$.  By a ``generalized moment class'', we mean that interval bounds are given for the expected values of finitely many\footnote{This is a ``philosophically reasonable'' position to take, since one can verify finitely many such inequalities in finite time.} test functions $\varphi$ against either the joint distribution $\mu$ or the marginal distributions $\mu_{k}$.  This setting encompasses a wide spectrum of possible dependence structures for the components of $X$, all the way from independence, through partial correlation (an inequality constraint on $\E_{\mu}[X_{i} X_{j}]$), to complete dependence ($X_{i}$ and $X_{j}$ are treated as a single random variable $(X_{i}, X_{j})$ with arbitrary joint distribution).  This setting also allows for coupling of the constraints on $g$ and those on $\mu$ (e.g.~by a constraint on $\E_{\mu}[g]$).

To express the previous paragraph more mathematically, we assume that our information about reality $(g^{\dagger}, \mu^{\dagger})$ is that it lies in the set $\mathcal{A}$ defined by
\begin{equation}
	\label{eq:A}
	\mathcal{A} 
	\defeq
	\left\{ (g, \mu) \,\middle|\, \begin{array}{c}
		\text{$g \colon \mathcal{X} = \mathcal{X}_{0} \times \dots \times \mathcal{X}_{K - 1} \to \mathcal{Y}$ is measurable,} \\
		\text{$\mu = \mu_{0} \otimes \dots \otimes \mu_{K - 1}$ is a product measure on $\mathcal{X}$,} \\
		\text{$\langle$conditions that constrain $g$ pointwise$\rangle$} \\
		\text{$\E_{\mu}[\varphi_{j}] \leq 0$ for $j = 1, \dotsc, N$,} \\
		\text{$\E_{\mu_{k}}[\varphi_{k, j_{k}}] \leq 0$ for $k = 0, \dotsc, K - 1$, $j_{k} = 1, \dotsc, N_{k}$}
	\end{array} \right\}
\end{equation}
for some known measurable functions $\varphi_{j} \colon \mathcal{X} \to \R$ and $\varphi_{k, j_{k}} \colon \mathcal{X}_{k} \to \R$.  In this case, the following reduction theorem holds:

\begin{thm}[{\cite[\S4]{OSSMO:2011}}]
	\label{thm:reduction}
	Suppose that $\mathcal{A}$ is of the form \eqref{eq:A}.  Then
	\begin{equation}
		\label{eq:equivalence}
		\underline{Q}(\mathcal{A}) = \underline{Q}(\mathcal{A}_{\Delta})
		\quad\text{and}\quad
		\overline{Q}(\mathcal{A}) = \overline{Q}(\mathcal{A}_{\Delta}),
	\end{equation}
	where
	\begin{equation}
		\label{eq:ADelta}
		\mathcal{A}_{\Delta}
		\defeq
		\left\{ (g, \mu) \in \mathcal{A} \,\middle|\, \begin{array}{c}
			\text{for $k = 0, \dotsc, K - 1$,} \\
			\text{$\mu_{k} = \sum_{i_{k} = 0}^{N + N_{k}} w_{k, i_{k}} \delta_{x_{k, i_{k}}}$} \\
			\text{for some $x_{k, 1}, x_{k, 2}, \dotsc, x_{k_{N + N_{k}}} \in \mathcal{X}_{k}$} \\
			\text{and $w_{k, 1}, w_{k, 2}, \dotsc, w_{k_{N + N_{k}}} \geq 0$} \\
			\text{with $w_{k, 1} + w_{k, 2} + \dots + w_{k_{N + N_{k}}} = 1$}
		\end{array} \right\}.
	\end{equation}
\end{thm}

Informally, Theorem~\ref{thm:reduction} says that if all one knows about the random variable $X = (X_{0}, \dotsc, X_{K - 1})$ is that its components are independent, together with inequalities on $N$ generalized moments of $X$ and $N_{k}$ generalized moments of each $X_{k}$, then for the purposes of solving \eqref{eq:ext_E_qoi} it is legitimate to consider each $X_{k}$ to be a \emph{discrete} random variable that takes at most $N + N_{k} + 1$ distinct values $x_{k, 0}, x_{k, 1}, \dotsc, x_{k, N + N_{k}}$;  those values $x_{k, i_{k}} \in \mathcal{X}_{k}$ and their corresponding probabilities $w_{k, i_{k}} \geq 0$ are the optimization variables.

For the sake of conciseness and to reduce the number of subscripts
required, multi-index notation will be used in what follows to express
the product probability measures $\mu$ of the form
\[
	\mu = \bigotimes_{k = 0}^{K - 1} \sum_{i_{k} = 0}^{N + N_{k}} w_{k, i_{k}} \delta_{x_{k, i_{k}}}
\]
that arise in the finite-dimensional reduced feasible set $\mathcal{A}_{\Delta}$ of equation \eqref{eq:ADelta}.  Write $\bs{i} \defeq (i_{0}, \dotsc, i_{K - 1}) \in \N_{0}^{K}$ for a multi-index, let $\bs{0} \defeq (0, \dotsc, 0)$, and let
\[
	\bs{M} \defeq (M_{0}, \dotsc, M_{K - 1}) \defeq (N + N_{0}, \dotsc, N + N_{K - 1}).
\]
Let $\# \bs{M} \defeq \prod_{k = 1}^{K} (M_{k} + 1)$.  With this notation, the $\# \bs{M}$ support points of the measure $\mu$, indexed by $\bs{i} = \bs{0}, \dotsc, \bs{M}$, will be written as
\[
	\bs{x}_{\bs{i}} \defeq (x_{1, i_{1}}, x_{2, i_{2}}, \dotsc, x_{K, i_{K}}) \in \mathcal{X}
\]
and the corresponding weights as
\[
	w_{\bs{i}} \defeq w_{1, i_{1}} w_{2, i_{2}} \dots w_{K, i_{K}} \geq 0,
\]
so that
\begin{equation}
	\label{eq:product_and_sum}
	\mu = \bigotimes_{k = 0}^{K - 1} \sum_{j_{k} = 0}^{N + N_{k}} w_{k, j_{k}} \delta_{x_{k, j_{k}}} = \sum_{\bs{i} = \bs{0}}^{\bs{M}} w_{\bs{i}} \delta_{\bs{x}_{\bs{i}}}.
\end{equation}
It follows from \eqref{eq:product_and_sum} that, for any integrand $f \colon \mathcal{X} \to \R$, the expected value of $f$ under such a discrete measure $\mu$ is the finite sum
\begin{equation}
	\label{eq:discrete_expect}
	\E_{\mu}[f] = \sum_{\bs{i} = \bs{0}}^{\bs{M}} w_{\bs{i}} f(\bs{x}_{\bs{i}})
\end{equation}
(It is worth noting in passing that conversion from product to sum representation and back as in \eqref{eq:product_and_sum} is an essential task in the numerical implementation of these UQ problems, because the product representation captures the independence structure of the problem, whereas the sum representation is best suited to integration (expectation) as in \eqref{eq:discrete_expect}.)

% \textbf{DO WE NEED THE REDUCTION OVER $g$ AS WELL?  IF NOT, THEN DITCH THE NEXT TWO PARAGRAPHS.}

Furthermore, not only is the search over $\mu$ effectively finite-dimensional, as guaranteed by Theorem~\ref{thm:reduction}, but so too is the search over $g$:  since integration against a measure requires knowledge of the integrand \emph{only} at the support points of the measure, only the $\# \bs{M}$ values $y_{\bs{i}} \defeq g(\bs{x}_{\bs{i}})$ of $g$ at the support points $\{ \bs{x}_{\bs{i}} \mid \bs{i} = \bs{0}, \dotsc, \bs{M} \}$ of $\mu$ need to be known.  So, for example, if $g^{\dagger}$ is known, then it is only necessary to evaluate it on the finite support of $\mu$.  Another interesting situation of this type is considered in \cite{Sullivan:2013}, in which $g^{\dagger}$ is not known exactly, but is known via legacy data at some points of $\mathcal{X}$ and is also known to satisfy a Lipschitz condition --- in which case the space of admissible $g$ is infinite-dimensional before reduction to the support of $\mu$, but the finite-dimensional collection of admissible values $( y_{\bs{0}}, \dotsc, y_{\bs{M}} )$ has a polytope-like structure.

Theorem~\ref{thm:reduction}, formulae \eqref{eq:product_and_sum}--\eqref{eq:discrete_expect}, and the remarks of the previous paragraph imply that $\overline{Q}(\mathcal{A})$ is found by solving the following finite-dimensional maximization problem (and $\underline{Q}(\mathcal{A})$ by the corresponding minimization problem):
\begin{align}
	\label{eq:ouq_explicit}
	\text{maximize:}\quad &  \sum_{\bs{i} = \bs{0}}^{\bs{M}} w_{\bs{i}} q(\bs{x}_{\bs{i}}, y_{\bs{i}}) ; \\
	\notag
	\text{among:}\quad
	& y_{\bs{i}} \in \mathcal{Y} \text{ for $\bs{i} = \bs{0}, \dotsc, \bs{M}$,} \\
	\notag
	& w_{k, i_{k}} \in [0, 1] \text{ for $k = 0, \dotsc, K - 1$ and $i_{k} = 0, \dotsc, M_{k}$,} \\
	\notag
	& x_{k, i_{k}} \in \mathcal{X}_{k} \text{ for $k = 0, \dotsc, K - 1$ and $i_{k} = 0, \dotsc, M_{k}$;} \\
	\notag
	\text{subject to:}\quad
	& y_{\bs{i}} = g(\bs{x}_{\bs{i}}) \text{ for some $\mathcal{A}$-admissible $g \colon \mathcal{X} \to \mathcal{Y}$,} \\
	\notag
	& \sum_{\bs{i} = \bs{0}}^{\bs{M}} w_{\bs{i}} \varphi_{j}(\bs{x}_{\bs{i}}) \leq 0 \text{ for $j = 1, \dotsc, N$,} \\
	\notag
	& \sum_{i_{k} = 0}^{M_{k}} w_{k, i_{k}} \varphi_{k, j_{k}}(x_{k, i_{k}}) \leq 0 \text{ for $k = 0, \dotsc, K - 1$ and $j_{k} = 1, \dotsc, N_{k}$,} \\
	\notag
	& \sum_{i_{k} = 0}^{M_{k}} w_{k, i_{k}} = 1 \text{ for $k = 0, \dotsc, K - 1$.}
	\notag
\end{align}

Generically, the reduced OUQ problem \eqref{eq:ouq_explicit} is non-convex, although there are special cases that can be treated using the tools of convex optimization and duality \cite{BP:2005, DY:2010, S:1995, VBC:2007}.  Therefore, numerical methods for global optimization must be employed to solve \eqref{eq:ouq_explicit}.  Unsurprisingly, the numerical solution of \eqref{eq:ouq_explicit} is much more computationally intensive when $\# \bs{M}$ is large --- the so-called \emph{curse of dimension}.
%} % END COMMENT 